%Paper: dg-ga/9501002
%From: Johan Rade <rade@leland.Stanford.EDU>
%Date: Sat, 7 Jan 1995 17:46:43 -0800
%Date (revised): Sat, 7 Jan 1995 17:58:07 -0800
%Date (revised): Mon, 9 Jan 1995 15:39:51 -0800
%Date (revised): Mon, 9 Jan 1995 16:12:10 -0800
%Date (revised): Wed, 25 Jan 1995 03:38:52 -0800
%Date (revised): Wed, 25 Jan 1995 04:36:42 -0800

%%%%%%%%%%%%%%%%%%%%%%%%%%%%%%%%%%%%%%%%%%%%%%%%%%%%%%%%%%%%%%%%%%%%%%%%%%%%
%%%%%%%%%%%%%%%%%%%%%%%%%%%%%%%%%%%%%%%%%%%%%%%%%%%%%%%%%%%%%%%%%%%%%%%%%%%%

%font.tex

%
%         \bb = blackboard bold --- math only
%         \euf = fraktur (gothic) --- math only
%         \bfsl = bold slanted --- text only
%         \bfit = bold italic --- text only
%         \seventeenpoint = switch to seventeen point size
%         \boldseventeenpoint = switch to bold seventeen point
%         \twelvepoint = switch to twelve point size
%         \boldtwelvepoint = switch to bold twelve point
%         \tenpoint = switch to ten point size (default)
%         \boldtenpoint = switch to bold ten point
%         \ninepoint = switch to nine point size
%         \eightpoint = switch to eight point size
%
%
%         Font families msa and msb are loaded and make all AMS TeX
%         characters defineable. See sample definitions at the end.
%

% Text Fonts

\font\twelverm     = cmr12
\font\tenrm        = cmr10
\font\ninerm       = cmr9
\font\eightrm      = cmr8
\font\sevenrm      = cmr7
\font\sixrm        = cmr6
\font\fiverm       = cmr5

\font\twelveit     = cmti12
\font\tenit        = cmti10
\font\nineit       = cmti9

\font\twelvesl     = cmsl10 at 12pt
\font\tensl        = cmsl10
\font\ninesl       = cmsl9

\font\twelvebf     = cmbx12
\font\tenbf        = cmbx10
\font\ninebf       = cmbx9
\font\eightbf      = cmbx8
\font\sevenbf      = cmbx7
\font\sixbf        = cmbx6
\font\fivebf       = cmbx5

\font\twelvebfit     = cmbxti10 at 12pt
\font\tenbfit        = cmbxti10
\font\ninebfit       = cmbxti10 at 9pt

\font\twelvebfsl     = cmbxsl10 at 12pt
\font\tenbfsl        = cmbxsl10
\font\ninebfsl       = cmbxsl10 at 9pt

\font\twelvett    = cmtt12
\font\tentt       = cmtt10
\font\ninett      = cmtt9

% Math Fonts

\font\twelvemit     = cmmi12
\font\tenmit        = cmmi10
\font\ninemit       = cmmi9
\font\eightmit      = cmmi8
\font\sevenmit      = cmmi7
\font\sixmit        = cmmi6
\font\fivemit       = cmmi5

\font\twelvebfmit    = cmmib10 at 12pt
\font\tenbfmit       = cmmib10

\font\eightbfmit     = cmmib8
\font\sevenbfmit     = cmmib7
\font\sixbfmit       = cmmib6
\font\fivebfmit      = cmmib5

\font\twelvesy    = cmsy10 at 12pt
\font\tensy       = cmsy10
\font\ninesy      = cmsy9
\font\eightsy     = cmsy8
\font\sevensy     = cmsy7
\font\sixsy       = cmsy6
\font\fivesy      = cmsy5

\font\twelvebfsy    = cmbsy10 at 12pt
\font\tenbfsy       = cmbsy10

\font\eightbfsy     = cmbsy8
\font\sevenbfsy     = cmbsy7
\font\sixbfsy       = cmbsy6
\font\fivebfsy      = cmbsy5

\font\tenex      = cmex10

% AMS TeX Math Fonts

\font\twelvemsa    = msam10 at 12pt
\font\tenmsa       = msam10
\font\ninemsa      = msam9
\font\eightmsa     = msam8
\font\sevenmsa     = msam7
\font\sixmsa       = msam6
\font\fivemsa      = msam5

\font\twelvemsb    = msbm10 at 12pt
\font\tenmsb       = msbm10
\font\ninemsb      = msbm9
\font\eightmsb     = msbm8
\font\sevenmsb     = msbm7
\font\sixmsb       = msbm6
\font\fivemsb      = msbm5

\font\twelveeuf    = eufm10 at 12pt
\font\teneuf       = eufm10
\font\nineeuf      = eufm9
\font\eighteuf     = eufm8
\font\seveneuf     = eufm7
\font\sixeuf       = eufm6
\font\fiveeuf      = eufm5

\font\teneufb       = eufb10

\font\twelveeus    = eusm10 at 12pt
\font\teneus       = eusm10
\font\nineeus      = eusm9
\font\eighteus     = eusm8
\font\seveneus     = eusm7
\font\sixeus       = eusm6
\font\fiveeus      = eusm5

\font\teneusb       = eusb10

% LaTeX fonts

\font\twelvelasy     = lasy10 at 12pt
\font\tenlasy        = lasy10
\font\ninelasy       = lasy9
\font\eightlasy      = lasy8
\font\sevenlasy      = lasy7
\font\sixlasy        = lasy6
\font\fivelasy       = lasy5

\font\tenblasy        = lasyb10

% Various fonts

%\font\tenfi      = cmfi10

%%%%%%%%%%%%%%%%%

\def\hexdigit#1{\ifcase#1%
  0\or1\or2\or3\or4\or5\or6\or7\or
  8\or9\or A\or B\or C\or D\or E\or F\fi}

\newfam\msafam
\newfam\msbfam
\newfam\euffam
\newfam\eusfam

\edef\hexmsa{\hexdigit\msafam}
\edef\hexmsb{\hexdigit\msbfam}
\def\bb{\fam\msbfam}
\def\euf{\fam\euffam}
\def\eus{\fam\eusfam}

\def\defaultfont{\rm\fam-1{}}

\newfam\lasyfam
\edef\hexlasy{\hexdigit\lasyfam}

%%%%%%%%%%%%%%%%%%%

\textfont\itfam        = \nullfont
\textfont\slfam        = \nullfont
\textfont3             = \tenex
\scriptfont3           = \tenex
\scriptscriptfont3     = \tenex

\def\twelvepoint{%  Text: 12  Math: 12 - 8 - 6
  \textfont0                 = \twelverm
  \scriptfont0               = \eightrm
  \scriptscriptfont0         = \sixrm
  \textfont1                 = \twelvemit
  \scriptfont1               = \eightmit
  \scriptscriptfont1         = \sixmit
  \textfont2                 = \twelvesy
  \scriptfont2               = \eightsy
  \scriptscriptfont2         = \sixsy
  \textfont\bffam            = \twelvebf
  \scriptfont\bffam          = \eightbf
  \scriptscriptfont\bffam    = \sixbf
  \textfont\msafam           = \twelvemsa
  \scriptfont\msafam         = \eightmsa
  \scriptscriptfont\msafam   = \sixmsa
  \textfont\msbfam           = \twelvemsb
  \scriptfont\msbfam         = \eightmsb
  \scriptscriptfont\msbfam   = \sixmsb
  \textfont\euffam           = \twelveeuf
  \scriptfont\euffam         = \eighteuf
  \scriptscriptfont\euffam   = \sixeuf
  \textfont\eusfam           = \twelveeus
  \scriptfont\eusfam         = \eighteus
  \scriptscriptfont\eusfam   = \sixeus
  \textfont\lasyfam          = \twelvelasy
  \scriptfont\lasyfam        = \eightlasy
  \scriptscriptfont\lasyfam  = \sixlasy
  \def\rm{\twelverm\fam0}%
  \def\it{\twelveit\fam1}%
  \def\sl{\twelvesl}%
  \def\bf{\twelvebf\fam\bffam}%
  \def\bfit{\twelvebfit}%
  \def\bfsl{\twelvebfsl}%
  \def\tt{\twelvett}%
  \twelvebaselines\normalbaselines\defaultfont}

\def\boldtwelvepoint{%  Text: 12  Math: 12 - 8 - 6
  \twelvepoint
  \textfont0                 = \twelvebf
  \scriptfont0               = \eightbf
  \scriptscriptfont0         = \sixbf
  \textfont1                 = \twelvebfmit
  \scriptfont1               = \eightbfmit
  \scriptscriptfont1         = \sixbfmit
  \textfont2                 = \twelvebfsy
  \scriptfont2               = \eightbfsy
  \scriptscriptfont2         = \sixbfsy
  \def\rm{\twelvebf\fam\bffam}%
  \def\it{\twelvebfit\fam1}%
  \def\sl{\twelvebfsl}%
  \defaultfont}

\def\tenpoint{%  Text: 10  Math: 10 - 7 - 5
  \textfont0                 =\tenrm
  \scriptfont0               =\sevenrm
  \scriptscriptfont0         =\fiverm
  \textfont1                 =\tenmit
  \scriptfont1               =\sevenmit
  \scriptscriptfont1         =\fivemit
  \textfont2                 =\tensy
  \scriptfont2               =\sevensy
  \scriptscriptfont2         =\fivesy
  \textfont\bffam            =\tenbf
  \scriptfont\bffam          =\sevenbf
  \scriptscriptfont\bffam    =\fivebf
  \textfont\msafam           =\tenmsa
  \scriptfont\msafam         =\sevenmsa
  \scriptscriptfont\msafam   =\fivemsa
  \textfont\msbfam           =\tenmsb
  \scriptfont\msbfam         =\sevenmsb
  \scriptscriptfont\msbfam   =\fivemsb
  \textfont\euffam           =\teneuf
  \scriptfont\euffam         =\seveneuf
  \scriptscriptfont\euffam   =\fiveeuf
  \textfont\eusfam           =\teneus
  \scriptfont\eusfam         =\seveneus
  \scriptscriptfont\eusfam   =\fiveeus
  \textfont\lasyfam          =\tenlasy
  \scriptfont\lasyfam        =\sevenlasy
  \scriptscriptfont\lasyfam  =\fivelasy
  \def\rm{\tenrm\fam0}%
  \def\it{\tenit\fam1}%
  \def\sl{\tensl}%
  \def\bf{\tenbf\fam\bffam}%
  \def\bfit{\tenbfit}%
  \def\bfsl{\tenbfsl}%
  \def\tt{\tentt}%
  \tenbaselines\normalbaselines\defaultfont}

\def\boldtenpoint{%  Text: 10  Math: 10 - 7 - 5
  \tenpoint
  \textfont0                 =\tenbf
  \scriptfont0               =\sevenbf
  \scriptscriptfont0         =\fivebf
  \textfont1                 =\tenbfmit
  \scriptfont1               =\sevenbfmit
  \scriptscriptfont1         =\fivebfmit
  \textfont2                 =\tenbfsy
  \scriptfont2               =\sevenbfsy
  \scriptscriptfont2         =\fivebfsy
  \textfont\euffam           =\teneufb
  \textfont\eusfam           =\teneusb
  \textfont\lasyfam          =\tenblasy
  \def\rm{\tenbf\fam\bffam}%
  \def\it{\tenbfit\fam1}%
  \def\sl{\tenbfsl}%
  \defaultfont}

\def\ninepoint{%  Text: 9  Math: 9 - 6 - 5
  \textfont0                 = \ninerm
  \scriptfont0               = \sixrm
  \scriptscriptfont0         = \fiverm
  \textfont1                 = \ninemit
  \scriptfont1               = \sixmit
  \scriptscriptfont1         = \fivemit
  \textfont2                 = \ninesy
  \scriptfont2               = \sixsy
  \scriptscriptfont2         = \fivesy
  \textfont\bffam            = \ninebf
  \scriptfont\bffam          = \sixbf
  \scriptscriptfont\bffam    = \fivebf
  \textfont\msafam           = \ninemsa
  \scriptfont\msafam         = \sixmsa
  \scriptscriptfont\msafam   = \fivemsa
  \textfont\msbfam           = \ninemsb
  \scriptfont\msbfam         = \sixmsb
  \scriptscriptfont\msbfam   = \fivemsb
  \textfont\euffam           = \nineeuf
  \scriptfont\euffam         = \sixeuf
  \scriptscriptfont\euffam   = \fiveeuf
  \textfont\eusfam           = \nineeus
  \scriptfont\eusfam         = \sixeus
  \scriptscriptfont\eusfam   = \fiveeus
  \textfont\lasyfam          = \ninelasy
  \scriptfont\lasyfam        = \sixlasy
  \scriptscriptfont\lasyfam  = \fivelasy
  \def\rm{\ninerm\fam0}%
  \def\it{\nineit\fam1}%
  \def\sl{\ninesl}%
  \def\bf{\ninebf\fam\bffam}%
  \def\bfit{\ninebfit}%
  \def\bfsl{\ninebfsl}%
  \def\tt{\ninett}%
  \ninebaselines\normalbaselines\defaultfont}

\def\twelvebaselines{\normalbaselineskip=14.4pt%
  \normallineskiplimit=0pt \normallineskip=1.2pt}
\def\tenbaselines{\normalbaselineskip=12pt%
  \normallineskiplimit=0pt \normallineskip=1pt%
  \setbox\strutbox=\hbox{\vrule height8.5pt depth3.5pt width0pt}}
\def\ninebaselines{\normalbaselineskip=10.8pt%
  \normallineskiplimit=0pt \normallineskip=0.9pt}

\tenpoint

% Examples of character definitions:

\mathchardef\boxplus=        "2\hexmsa01
\mathchardef\boxtimes=       "2\hexmsa02
\mathchardef\square=         "0\hexmsa03
\mathchardef\rightleftarrow= "3\hexmsa1D
\mathchardef\rightsquigarrow="3\hexmsa20
\mathchardef\therefore=      "3\hexmsa29
\mathchardef\because=        "3\hexmsa2A
\mathchardef\varpropto=      "3\hexmsa5F
\mathchardef\Subset=         "3\hexmsa62
\mathchardef\pitchfork=      "3\hexmsa74

\mathchardef\ltimes=         "2\hexmsb6E
\mathchardef\rtimes=         "2\hexmsb6F
\mathchardef\eth=            "0\hexmsb67
\mathchardef\varkappa=       "0\hexmsb7B

\mathchardef\Diamond=        "2\hexlasy33

%%%%%%%%%%%%%%%%%%%%%%%%%%%%%%%%%%%%%%%%%%%%%%%%%%%%%%%%%%%%%%%%%%%%%%%%%%%%
%%%%%%%%%%%%%%%%%%%%%%%%%%%%%%%%%%%%%%%%%%%%%%%%%%%%%%%%%%%%%%%%%%%%%%%%%%%%

%macros.tex

\hyphenation{mono-pole mono-poles}

\def\cases#1{\left\{\,\vcenter{\openup\jot
  \ialign{$\ds##\hfil$&\quad##\hfil\crcr#1\crcr}}\right.}

\def\eqalign#1{\null\,\vcenter{\openup\jot
  \ialign{&\strut\hfil$\ds##$%
          &$\ds{}##$\hfil\cr#1\crcr}}\,}

% Note the & after \ialign{. This allows \eqalign with any number of entries
% on each line.

\mathchardef\oldepsilon="010F
\mathchardef\epsilon="0122
\mathchardef\oldphi="011E
\mathchardef\phi="0127
\mathcode`\*="0203   % makes * an ord-atom rather than a bin-atom

%See The TeX book, p.179

\def\itimes{\mathbin{\raise1.5pt\hbox{\vrule width6pt height.4pt depth0pt
  \kern-.1pt\vrule width.4pt height2pt depth0pt}}}

\def\diracc#1#2{\ooalign{\hfill$#1\mskip1mu/$\hfill\cr
  $#1\raise-.2ex\hbox{$#1#2$}$}}
\def\dbar{{\overline\partial}}

\def\IC{{\bb C}}  \def\IE{{\bb E}}  \def\IF{{\bb F}}

\def\IR{{\bb R}}  \def\IZ{{\bb Z}}

    \def\eusC{{\eus C}}
\def\eusD{{\eus D}}    
  \def\eusH{{\eus H}}  
\def\eusJ{{\eus J}}    
\def\eusM{{\eus M}}    \def\eusO{{\eus O}}
\def\eusP{{\eus P}}

\chardef\norwegiano="1C
\chardef\norwegianO="1F
\def\breveaccent#1{{\accent"15 #1}}
\def\cedilla#1{\setbox0=\hbox{#1}\ifdim\ht0=1ex \accent'30 #1%
  \else{\ooalign{\hidewidth\char'30\hidewidth\crcr\unhbox0}}\fi}
\def\eufu{{\euf u}}
\def\Chern{{\rm c}}

\def\o{\relax\ifmmode{\euf o}\else\norwegiano\fi}

\def\u{\relax\ifmmode\let\next=\eufu\else\let\next=\breveaccent\fi\next}

\def\O{\relax\ifmmode{\rm O}\else\norwegianO\fi}

\def\U{{\rm U}}

\def\c{{\ifmmode\let\next=\Chern\else\let\next=\cedilla\fi\next}}

\def\dist{\mathop{\rm dist}\nolimits}

\def\index{\mathop{\rm index}\nolimits}

\def\Hom{\mathop{\rm Hom}\nolimits}

\def\range{\mathop{\rm range}}

\def\comdia#1{{\advance\normallineskip by3pt%
               \advance\normallineskiplimit by3pt %
               \matrix{#1}}}
\def\mapright#1{\mathrel{\smash{\mathop{\longrightarrow}\limits^{#1}}}}
\def\mapdown#1{\Big\downarrow\rlap{$\vcenter{\hbox{$\scriptstyle#1$}}$}}

\def\plus{{\sss+}}
\def\minus{{\sss-}}

\def\ds{\displaystyle}
\def\ts{\textstyle}
\def\ss{\scriptstyle}

\def\iint{\int\!\!\int}

\def\dingbat{\smallskip\line{\hfill\vbox{\hrule width1cm}\hskip.5cm%
  \lower2pt\hbox{$\circ$}\hskip.5cm\vbox{\hrule width1cm}\hfill }\smallskip}

\def\today{\ifcase\month\or January\or February\or March\or April\or May%
  \or June\or July\or August\or September\or October\or November\or December%
  \fi\space\number\day,\space\number\year}

\font\manfont=manfnt
\def\strutdepth{\dp\strutbox}
\def\marginaldanger{\strut\vadjust{\kern-\strutdepth\specialdanger}}
\def\specialdanger{\vtop to\strutdepth{\baselineskip\strutdepth
   \vss\llap{\manfont\char"7E\qquad\qquad}\null}}

%%%%%%%%%%%%%%%%%%%%%%%%%%%%%%%%%%%%%%%%%%%%%%%%%%%%%%%%%%%%%%%%%%%%%%%%%%%%
%%%%%%%%%%%%%%%%%%%%%%%%%%%%%%%%%%%%%%%%%%%%%%%%%%%%%%%%%%%%%%%%%%%%%%%%%%%%

%jrpaperstyle.tex

%\ifundefined{AmSTeX}\else\error{AMSTeX}\fi

%parameters
\magnification\magstephalf
\font\titlefont=cmbx12 at 16pt
\hoffset=.5truein
\hsize=5.7truein
\voffset=.4375truein
\vsize=8.25truein
\parskip=6pt plus2pt minus2pt
\smallskipamount=6pt plus2pt minus2pt
\medskipamount=12pt plus4pt minus4pt
\bigskipamount=18pt plus6pt minus6pt
\abovedisplayskip=12pt plus4pt minus4pt
\belowdisplayskip=12pt plus4pt minus4pt
\def\ninebaselines{\normalbaselineskip=12pt%
  \normallineskiplimit=1pt \normallineskip=1pt}
\def\tenbaselines{\normalbaselineskip=14pt plus1pt minus1pt%
  \normallineskiplimit=1pt \normallineskip=1pt}
\tenpoint
\newdimen\titlebaselineskip \titlebaselineskip=24pt
\newdimen\authorbaselineskip \authorbaselineskip=18pt
\def\smallrule{\centerline{\vbox{\hrule width 1in}}}

\newif\ifshowlabel\showlabelfalse

% page layout
\headline={\ifnum\pageno>1\ifodd\pageno\hfil{\ninerm\folio}%
  \else{\ninerm\folio}\hfil\fi\else\hfil\fi}
\footline={\hfil}
\def\makeheadline{\vbox to0pt{\vskip-34pt%
   \line{\vbox to8.5pt{}\the\headline}\vss}\nointerlineskip}

% title
\def\titlelines#1\par{\leavevmode\bigskip\vskip\parskip\nointerlineskip
  \centerline{\vbox{\baselineskip=\titlebaselineskip\titlefont
  \halign{\hfil##\hfil\cr#1\crcr}}}\par}
\def\authorlines#1\par{\medskip\vskip\parskip\nointerlineskip
  \centerline{\vbox{\twelvepoint\baselineskip=\authorbaselineskip
  \halign{\hfil##\hfil\cr#1\crcr}}}\par}
\def\titlerule{\bigskip\vskip\parskip\nointerlineskip\smallrule
  \bigskip\vskip\parskip\prevdepth=-1001pt}
\def\abstract#1\par{{\narrower\narrower\noindent#1\par}}
\def\acknowledgements#1\par{{\ninepoint\subsubheading Acknowledgements.\par
   #1\par}}

% body
\newif\iffirstheading\firstheadingtrue
\def\heading#1\par{%
  \iffirstheading\firstheadingfalse
  \else \vskip0pt plus.3\vsize\penalty-250\vskip0pt plus-.3\vsize\fi
  \ifdim\prevdepth=-1001pt\else\bigskip\vskip\parskip\fi
  \centerline{\vbox{\boldtwelvepoint\baselineskip=\authorbaselineskip
  \halign{\hfil##\hfil\cr#1\crcr}}}%
  \smallskip}
\def\subheading#1\par{\removelastskip\smallskip
  \penalty-75\noindent{\boldtenpoint#1}\enspace}
\def\subsubheading#1\par{\penalty-50{\sl#1}\enspace}
\long\def\proclaim#1\par#2\par{\ifdim\lastskip<\smallskipamount
  \removelastskip\smallskip\fi\penalty-100{\bf#1}\enspace
  {\sl#2}\smallskip\penalty100}

\def\demo#1\par{\ifdim\lastskip<\smallskipamount\removelastskip\smallskip
  \penalty-50\fi{\sl#1}\enspace}
\def\subdemo{\subsubheading}
\def\enddemo{\penalty200\hbox{}\nobreak\hfill\space\hbox{$\square$}\smallskip
  \penalty-75}
\def\item#1\par#2\par{{\advance\leftskip by 2\parindent\noindent
  \kern-\parindent\rlap{#1}\kern\parindent#2\par}}

% end
\def\references{\heading\tenpoint\bf References\par\smallskip
  \parskip=2pt plus.7pt minus.7pt\tolerance=1000\overfullrule=0pt
  \frenchspacing\ninepoint}
\def\ref[#1]#2\par {\noindent\kern-2.5cm\hbox to1.5cm{[#1]\hfil}#2%
  \leftskip= 2.5cm\par\leftskip= 0cm}
\def\endrule{\medskip\vskip\parskip\nointerlineskip\smallrule
  \medskip\nointerlineskip}
\def\endlines#1{\vskip\parskip
  \centerline{\ninepoint\vbox{\halign{\hfil##\hfil\cr#1\crcr}}}}

% equation numbers
\def\eqno#1{\leqno{%
  \ifshowlabel\llap{\hbox to.5in{\eightrm\string#1\hss}}\fi
  \hskip\parindent#1}}
\catcode`\@=11
\def\eqalignno#1{\displ@y\tabskip=\centering
  \halign to\displaywidth{\hfil$\@lign\displaystyle{##}$\tabskip=0pt
    &$\@lign\displaystyle{{}##}$\hfil\tabskip=\centering
    &\kern-\displaywidth\rlap{\hskip\parindent$\@lign##$}%
    \tabskip=\displaywidth\crcr#1\crcr}}
\def\displaylinesno#1{\displ@y\tabskip=\centering
  \halign to\displaywidth{\hfil$\@lign\displaystyle{##}$\hfil
  &\kern-\displaywidth\rlap{\hskip\parindent$\@lign##$}%
  \tabskip=\displaywidth\crcr#1\crcr}}
\catcode`\@=12
\newcount\labelcount
\newcount\sectioncount
\labelcount=0
\sectioncount=1
\def\endassign{apa}
\def\nextsection{gnurgelpupp}
\def\beginassign#1{\ifx#1\endassign\relax
  \else\ifx#1\nextsection\advance\sectioncount by1\labelcount=0
  \else\advance\labelcount by1\edef#1{{\rm(\the\sectioncount.\the\labelcount)}}
  \fi\expandafter\beginassign\fi}

%%%%%%%%%%%%%%%%%%%%%%%%%%%%%%%%%%%%%%%%%%%%%%%%%%%%%%%%%%%%%%%%%%%%%%%%%%%%
%%%%%%%%%%%%%%%%%%%%%%%%%%%%%%%%%%%%%%%%%%%%%%%%%%%%%%%%%%%%%%%%%%%%%%%%%%%%

%%%%%%%%%%%%%%%%%%%%%%%%%%%%%%%%
%
%    Products and relations in symplectic Floer homology
%
%               by Marty Betz and Johan Rade
%
%
%    This file is in plain TeX
%
%
%%%%%%%%%%%%%%%%%%%%%%%%%%%%%%%%%

\def\+{\tabalign}  % This redefines \+ to be non-outer,
                   % so it can be used inside an argument of proclaim.

\def\gop{\Diamond}
\def\cop{\mathchar"0\hexmsa7B}

\def\coker{\mathop{\rm coker}}
\def\max{{\rm max}}
\def\fib{{\rm fib}}
\def\domain{\mathop{\rm domain}}
\def\Diff{\mathop{\rm Diff}\nolimits}

\def\rmE{{\rm E}}
\def\rmF{{\rm F}}

\def\bfJ{{\bf J}}
\def\bfT{{\bf T}}
\def\bfg{{\bf g}}
\def\bfc{{\bf c}}
\def\bfj{{\bf j}}
\def\eufo{{\euf o}}
\def\eufJ{{\euf J}}
\def\eufO{{\euf O}}
\def\eufP{{\euf P}}
\def\eufM{{\euf M}}
\def\eufU{{\euf U}}
\def\eufV{{\euf V}}

\font\sans = cmss10
\def\ssJ{\hbox{\sans J}}
\def\ssQ{\hbox{\sans Q}}

\def\eufOtilde{\widetilde{\euf\tilde O}}
\def\eusJhat{\widehat\eusJ}
\def\Sigmabar{{\overline\Sigma}}
\def\alphabar{{\overline\alpha}}

\def\muH{\mu_{\scriptscriptstyle H}}
\def\phiH{\oldphi_{\scriptscriptstyle H}}
\def\plusminus{{\smash{\scriptscriptstyle\pm}}} % Why \smash?
\def\minus{{\smash{\scriptscriptstyle-}}}       % To prevent disaster
\def\plus{{\smash{\scriptscriptstyle+}}}        % in Sect. 3.4.

%%%%%%%%%%%%%%%%%%%%%%%%%%%%%%%%%%%%%%%%%

\beginassign
  \nextsection
  \qbx\qby\qba\qbz\qbb\qbt\qbc\nextsection
  \qca\qcaa\qcb\qcc\qcd\qce\qcf\qcg\qch\qci\qcj\qcx\qck\qcl\qcm\qco
    \qnp\qcq\nextsection
  \nextsection
\endassign

%%%%%%%%%%%%%%%%%%%%%%%%%%%%%%%%%%%%%%%%%%

\titlelines
  Products and Relations\cr
  in Symplectic Floer Homology\cr

\authorlines
  by {\sl Martin Betz} at Austin \cr
  and {\sl Johan R\aa de} at Lund and Stanford \cr

\titlerule

%%%%%%%%%%%%%%%%%%%%%%%%%%%%%%%%%%%%%%%%%%%

\heading \S1. Introduction

This paper gives a detailed and functorial treatment of
products, operations and relations in
Floer homology and Floer cohomology
of monotone symplectic manifolds.
Floer (co)homology groups were
introduced by A.~Floer in a series of papers [F1], [F2], [F3] and [F4].
Basic material on Floer (co)homology can also be found in
[HS], [HZ], [M], [MS1], [S] and [SZ]; see also [Sch1].
Let $M$ be
a monotone symplectic manifold of dimension $2n$.
The Floer homology groups $HF_*(M)$
are given by a chain complex $CF_*(M)$
with  modules generated by the contractible periodic orbits
of a Hamiltonian flow on $M$.
The differential is defined
by counting points in moduli spaces of perturbed pseudoholomorphic curves
parametrized by the cylinder $\IR\times S^1$.
The Floer cohomology groups $HF^*(M)$ are given by the dual complex
$CF^*(M)$.
The Floer (co)homology groups are isomorphic to the ordinary
(co)homology groups $H_*(M)$ and $H^*(M)$,
except that the grading may be reduced modulo an even integer;
see [F3] and [F4] Thm.~5.

The cylinder $\IR\times S^1$ can be viewed as a twice punctured sphere.
In this paper we consider moduli spaces
of perturbed pseudoholomorphic curves parametrized by
Riemann surfaces $\Sigma_{g,k^\plus,k^\minus}$ of genus $g$
with $k=k^\minus+k^\plus$ punctures,
of which  $k^\minus$ punctures are labeled ``negative''
and $k^\plus$ punctures are labeled ``positive''.
These moduli spaces are determined by a choice of conformal structure
on $\Sigma_{g,k^\minus,k^\plus}$ and some additional perturbation data.
These choices are parametrized by a
space $\eusJhat_{g,k^\minus,k^\plus}$.
The construction is invariant under a  framed diffeomorphism group
$\Diff_{g,k}$,
and descends to
$\eufJ_{g,k^\minus,k^\plus}
 =\eusJhat_{g,k^\minus,k^\plus}/\Diff_{g,k}$.
The space $\eusJhat_{g,k^\minus,k^\plus}$ is contractible.
The group $\Diff_{g,k}$ acts freely on a subset
$\eusJhat^*_{g,k^\minus,k^\plus}$
of
$\eusJhat_{g,k^\minus,k^\plus}$
whose complement has infinite codimension.
We let
$\eufJ^*_{g,k^\minus,k^\plus}
 = \eusJ^*_{g,k^\minus,k^\plus} / \Diff_{g,k}$.
This is essentially a classifying space for
$\Diff_{g,k}$.

The singular homology of $\eufJ^*_{g,k^\minus,k^\plus}$ is given by
a chain complex $C_*(\eufJ^*_{g,k^\minus,k^\plus})$
where the modules are generated by maps $\sigma$
from a standard simplex to $\eufJ^*_{g,k^\minus,k^\plus}$.
Such a map also parametrizes moduli spaces $\eufM^d_\sigma$
of dimension $d$.
By counting the number of points in
$\eufM^0_\sigma$ we essentially get a homomorphism
$$Q:C_*(\eufJ^*_{g,k^\minus,k^\plus})
   \to CF^*(M)^{\otimes k^\minus} \otimes CF_*(M)^{\otimes k^\plus}.$$
The product $S_{k^\minus}\times S_{k^\plus}$
of the permutation groups on $k^\minus$ and $k^\plus$ letters
acts on the left hand side
by permuting the punctures,
and on the right hand  side by permuting the factors.
Theorem 4.3.2 and Theorem 4.3.3 essentially state that $Q$
is an $S_{k^\minus}\times S_{k^\plus}$
equivariant chain map of degree $2n(1-g-k^\minus)$.
We prove this by counting the boundary points
in compactifications of 1-dimensional moduli spaces.
Hence there is an induced equivariant homomorphism
$$Q: H_*(\eufJ^*_{g,k^\minus,k^\plus})
   \to H\bigl(CF^*(M)^{\otimes k^\minus} \otimes CF_*(M)^{\otimes
k^\plus}\bigr) .
$$

Let $\Theta_{g,k^\minus,k^\plus}$ denote the canonical generator
of $H_0(\eufJ^*_{g,k^\minus,k^\plus})$.
Then $Q(\Theta_{0,1,1})\in H\bigl(CF^*(M)\otimes CF_*(M)\bigr)$
gives the identity maps $HF_*(M)\to HF_*(M)$
and $HF^*(M)\to HF^*(M)$.
A generalizations of this construction gives
canonical isomorphisms
between the Floer (co)homology groups defined using different
Hamiltonians and different almost complex structures on $M$.
The classes $Q(\Theta_{0,2,0})$ and $Q(\Theta_{0,0,2})$
give symplectic Poincar\'e duality maps.
The classes $Q(\Theta_{0,1,2})$ and $Q(\Theta_{0,2,1})$
give symplectic cup and intersection products.
The classes $Q(\Theta_{0,1,0})$ and $Q(\Theta_{0,0,1})$
give a symplectic unit class in $HF^0(M)$
and a symplectic top class in $HF_{2n}(M)$.

One can obtain the punctured surface
$\Sigma_{g_1+g_2,k^\minus_1+k^\minus_2,k^\plus_1+k^\plus_2}$
by gluing the punctured surface $\Sigma_{g_1,k^\minus_1,k^\plus_1+1}$ to
the punctured surface $\Sigma_{g_2,k^\minus_2+1,k^\plus_2}$.
This gives a map
$$\gop^\ell_{ij}:\eufJ^*_{g_1,k_1^\minus,k_1^\plus+1}
 \times
 \eufJ^*_{g_2,k_2^\minus+1,k_2^\plus}
 \to
  \eufJ^*_{g_1+g_2,k_1^\minus+k_2^\minus,k_1^\plus+k_2^\plus}.
$$
Theorem 4.3.4 states that the diagram
$$\hss\comdia{
    H_*(\eufJ^*_{g_1,k_1^\minus,k_1^\plus+1})
    \otimes
    H_*(\eufJ^*_{g_2,k_2^\minus+1,k_2^\plus})
  & \mapright{Q\otimes Q}
  & {\ninepoint\eqalign{
        & H\bigl(CF^*(M)^{\otimes k_1^\minus}
          \otimes
          CF_*(M)^{\otimes(k^\plus_1+1)}\bigr) \cr
      & \otimes H\bigl(CF^*(M)^{\otimes(k_2^\minus+1)}
          \otimes
          CF_*(M)^{\otimes k_2^\plus}\bigr) \cr
            }} \cr
    \mapdown{\gop^{}_{ij}}
  &
  & \mapdown{\gop^{}_{ij}} \cr
    H_*(\eufJ^*_{g_1+g_2,k_1^\minus+k^\minus_2,k_1^\plus+k_2^\plus})
  & \mapright{Q}
  & H\bigl(CF^*(M)^{\otimes(k_1^\minus+k_2^\minus)}
    \otimes
    CF_*(M)^{\otimes(k^\plus_1+k_2^\plus)}\bigr) \cr
          }
  \hss
$$
commutes.
The homomorphism on the right,
also denoted $\gop^{}_{ij}$,
is induced by the pairing
$CF^*(M)\otimes CF_*(M)\to\IZ$.
In \S4.4 we show that as a consequence,
symplectic Poincar\'e duality is an isomorphism,
the symplectic cup and intersection products are associative,
they are intertwined by symplectic Poincar\'e duality,
the symplectic unit class is an identity element for the
symplectic cup product,
and the symplectic top class is an identity element for the
symplectic intersection product.
As another consequence, in \S4.5 we show that
the usual scheme for identifying Floer (co)homology groups
defined using different Hamiltonians and almost complex structures
is consistent,
and that the homomorphism $Q$ is independent of these choices.

One obtains the punctured surface $\Sigma_{g+1,k^\minus,k^\plus}$
by gluing the punctured surface $\Sigma_{g,k^\minus+1,k^\plus+1}$ to itself.
This gives a map
$$\cop^\ell_{ij}:\eufJ^*_{g,k^\minus+1,k^\plus+1} \to
\eufJ^*_{g+1,k^\minus,k^\plus}.$$
Theorem 4.3.5 states that the diagram
$$\comdia{
    H_*(\eufJ^*_{g,k^\minus+1,k^\plus+1})
  & \mapright{Q}
  & H\bigl(CF^*(M)^{\otimes(k^\minus+1)}
    \otimes
    CF_*(M)^{\otimes(k^\plus+1)}\bigr) \cr
    \mapdown{\cop_{ij}}
  &
  & \mapdown{\cop_{ij}} \cr
    H_*(\eufJ^*_{g+1,k^\minus,k^\plus})
  & \mapright{Q}
  & H\bigl(CF^*(M)^{\otimes k^\minus}
    \otimes
    CF_*(M)^{\otimes k^\plus }\bigr) \cr
          }
$$
commutes.
The homomorphism on the right,
also denoted $\cop_{ij}$,
is induced by the pairing
$CF^*(M)\otimes CF_*(M)\to\IZ$.
As a consequence of this theorem,
the number of isolated perturbed pseudoholomorphic tori in $M$,
with a given conformal structure,
equals the Euler characteristic of $M$.

By Theorem 4.3.6 the homomorphism $Q$ factors through a
homomorphism
$$Q^0 : H_*(\eufJ^{0*}_{g,k^\minus,k^\plus})
   \to H(CF^*(M)^{\otimes k^\minus} \otimes CF_*(M)^{\otimes k^\plus})$$
where $\eufJ^{0*}_{g,k^\minus,k^\plus}$ is the unframed analogue
of $\eufJ^*_{g,k^\minus,k^\plus}$.
In particular,
$Q$ vanishes on the classes given by twists of the framings
and on the classes corresponding to Dehn twists.

These theorems form an attempt to make rigorous the conjectural
Gromov-Witten classes [W], [V], [Ru], [KM]; see also [Fu].
These are homomorphisms
$$\ssQ:H_*(\ssJ_{g,k^\minus,k^\plus})
   \to H(C^*(M)^{\otimes k^\minus} \otimes C_*(M)^{\otimes k^\plus}),$$
with the same gluing properties as our homomorphisms $Q$,
where $\ssJ_{g,k^\minus,k^\plus}$ is the
Mumford-Deligne compactification of the moduli space
of conformal structures on $\Sigmabar_{g,k^\minus,k^\plus}$.

Another rigorous approach to Gromov-Witten classes,
in the genus 0 case,
is the quantum cohomology of Y.~Ruan and G.~Tian [RT1];
see also [MS1].
The symplectic cup product and its relation to
quantum cohomology is discussed in [MS], [P], [RT2] and [Sch2].

In a second paper we will extend our results to
weakly monotone symplectic manifolds,
in the sense of Hofer and Salamon [HS],
by imposing additional transversality conditions
on the simplices spanning $C_*(\eufJ^*_{g,k^\minus,k^\plus})$.

\acknowledgements
The first author would like to thank Ralph Cohen for valuable
advice and encouragement.

%%%%%%%%%%%%%%%%%%%%%%%%%%%%%%%%%%%%%%%%%%%%%%

\heading \S2. Moduli Spaces of Perturbed Pseudoholomorphic Curves

\subheading 2.1. Floer's moduli spaces $\eusM^d(\alpha^\minus,\alpha^\plus)$.

Let $M$ be a compact smooth symplectic manifold of dimension $2n$
with symplectic form $\omega$.
Let $H$ be a smooth function $M\times S^1\to\IR$.
The function $H$ will play the role of a periodic time-dependent Hamiltonian
on $M$.
The periodic orbits for $H$ are the solutions $\alpha:S^1\to M$ to
Hamilton's equations
$${d\alpha\over d\theta} = X_H(\alpha(\theta),\theta) ,
  \eqno\qbx
$$
where $X_H(\cdot,\theta)$ denotes the symplectic gradient of $H(\cdot,\theta)$.
Let $\eusC_H$ denote the set of periodic orbits.
Let $\eusC^0_H$ denote the set of contractible periodic orbits.

The Hamiltonian induces an exact symplectomorphism
$\phiH:M\to M$ by $\phiH(p) = \alpha(2\pi)$,
where $\alpha:[0,2\pi]\to M$ is the unique solution
to \qbx{} with initial data $\alpha(0)=p$.
The periodic orbits correspond to the fixed points
of $\phiH$.

\proclaim Definition 2.1.1.

{
A periodic orbit $\alpha$ is regular
if and only if\/ $1$ is not an eigenvalue of $d_{\alpha(0)}\phiH$.
Equivalently, $\alpha$ is regular if and only if
the graph of $\phiH$ intersects the diagonal transversely
in $M\times M$ at $(\alpha(0),\alpha(0))$.

A Hamiltonian $H$ is regular if and only if all periodic orbits
for $H$ are regular.
}

Let $\eusH_M$ denote the space of smooth time-dependent
Hamiltonians on $M$.

\proclaim Proposition 2.1.2.

The regular time-dependent
Hamiltonians form a dense open subset of $\eusH_M$.

This is proven in [F4] Prop.~2c.1 and [SZ] Thm.~8.1.

An almost complex structure on $M$ is a smooth bundle map
$\bfJ:TM\to TM$ such that $\bfJ^2=-1$.
An almost complex structure $\bfJ$ on $M$ is said to be compatible
with $\omega$ if
$\langle\xi,\zeta\rangle = \omega(\xi,\bfJ\zeta)$ is a Riemannian metric
on $M$.
We denote the space of compatible almost complex structures on $M$
by $\eusJ_M$.
If $\bfJ$ is compatible with $\omega$,
then $X_H=\bfJ\nabla H$.
Hamilton's equations \qbx{} then take the form
$$\bfJ{d\alpha\over d\theta} +\nabla H=0.$$

\proclaim Proposition 2.1.3.

The space $\eusJ_M$ is contractible.

This is well known; see for instance [HZ] p.~15.
It follows from Prop.~2.1.3
that $TM$ has well defined Chern classes,
independently of the choice of $\bfJ$.

\proclaim Assumption 2.1.4.

Throughout the paper we assume that $M$ is monotone.
This means that there exists $k\ge0$ such that
$\langle \omega,a\rangle
 = k \langle \c_1(TM), a\rangle
$
for all $a\in\pi_2(M)$.

In a second paper we will extend the results of this paper to
weakly monotone symplectic manifolds in the sense of [HS].

Given a regular $A=(H,\bfJ)$,
Floer considered smooth maps
$$u:\IR\times S^1\to M$$
that satisfy the equation
$${\partial u\over\partial t} + \bfJ {\partial u\over\partial\theta}
  + \nabla H(u,\theta) = 0,
  \eqno\qby
$$
where $t$ and $\theta$ are the $\IR$ and $S^1$ coordinates,
and have finite energy,
$$\int_{-\infty}^\infty\int_{S^1}
  \left(
  \left| \partial u\over\partial t\right|^2
  + \left| \bfJ{\partial u\over\partial \theta} + \nabla H \right|^2
  \right)
  \,dt\,d\theta < \infty.
$$
Floer proved that if the Hamiltonian $H$ is regular,
then for any finite energy solution $u$ to \qby,
there exist periodic orbits $\alpha^\minus$ and $\alpha^\plus$
such that
$$\eqalign{
    u(t,\theta)&\to\alpha^\minus(\theta) &\qquad&\hbox{as $t\to-\infty$} \cr
    u(t,\theta)&\to\alpha^\plus(\theta)  &&\hbox{as $t\to\infty$} . \cr
          }
  \eqno\qba
$$
Here $u$ and all its derivatives converge exponentially fast.
For more details, see [F4] and [R] \S3.
He then defined the moduli space
$$\eusM(\alpha^\minus,\alpha^\plus)$$
as the set of finite energy solutions to \qby{} that satisfy \qba,
and he set up a deformation theory for these moduli spaces.
He showed that
$$\eusM(\alpha^\minus,\alpha^\plus)=
  \bigcup_{d\in\IZ} \eusM^d(\alpha^\minus,\alpha^\plus)
$$
where $\eusM^d(\alpha^\minus,\alpha^\plus)$
is the zero set of a Fredholm section
$\Psi^d(\alpha^\minus,\alpha^\plus)$ of index $d$
of a Hilbert space bundle $\rmF^d(\alpha^\minus,\alpha^\plus)$
over a Hilbert manifold
$\eusP^d(\alpha^\minus,\alpha^\plus)$.
For more details, see [F4] and [R] \S2.

\proclaim Definition 2.1.5.

{
The moduli space $\eusM^d(\alpha^\minus,\alpha^\plus)$
is regular
if and only if $H$ is regular
and the section $\Psi^d(\alpha^\minus,\alpha^\plus)$ is transverse
to the zero section of $\rmF^d(\alpha^\minus,\alpha^\plus)$.

The pair $A=(H,\bfJ)$ is regular
if and only if $H$ is regular,
the moduli spaces
$\eusM^d(\alpha^\minus,\alpha^\plus)$
are regular for all $d\in\IZ$ and all $\alpha^\minus,\alpha^\plus\in\eusC_H$,
and there are no pseudoholomorphic spheres in $M$ with $\c_1=1$
that intersect periodic orbits of $H$.
}

A regular moduli space $\eusM^d(\alpha^\minus,\alpha^\plus)$
is an embedded submanifold of $\eusP^d(\alpha^\minus,\alpha^\plus)$
of dimension $d$.

\proclaim Proposition 2.1.6.

The regular pairs $A=(H,\bfJ)$
form a dense subset of $\eusH_M\times\eusJ_M$.

In fact, it is a subset of the second category in the sense of Baire.
An outline of a proof was given in [F4] Prop.~2c.2.
More details can be found in [SZ] Thm.~8.4, [HS] and [FHS].

The moduli spaces $\eusM^d(\alpha^\minus,\alpha^\plus)$ are
translation invariant.
The only translation invariant finite energy solutions to \qby{}
are the trivial ones,
$$\alphabar_0(t,\theta)=\alpha_0(\theta),$$
in $\eusM^0(\alpha_0,\alpha_0)$.
This has the following consequence.

\proclaim Proposition 2.1.7.

If $A=(H,\bfJ)$ is regular,
then the moduli space $\eusM^0(\alpha_0,\alpha_0)$
contains a single point,
the trivial solution $\alphabar_0$,
and if $\alpha^\minus\ne\alpha^\plus$,
then the moduli space $\eusM^0(\alpha^\minus,\alpha^\plus)$ is empty.

\subheading 2.2. The parameter space $\eufJ^*_{g,k^\minus,k^\plus}(A)$.

Floer considered maps from the cylinder $\IR\times S^1$ to $M$.
The cylinder $\IR\times S^1$
can also be viewed as a twice punctured sphere.
We will consider maps from a Riemann surfaces of any genus
with any number of punctures.
We need to take into account that such
a surface in general has a nontrivial moduli spaces of conformal structures.

Fix a regular pair $A=(H,\bfJ)$.
For each triple $(g,k^\minus,k^\plus)$ of non-negative integers,
we fix a compact oriented surface $\Sigmabar_{g,k^\minus,k^\plus}$
of genus $g$,
$k=k^\minus+k^\plus$ distinct points
$p^\minus_1,\ldots,p^\minus_{k^\minus},p^\plus_1,\ldots,p^\plus_{k^\plus}$
on $\Sigmabar_{g,k^\minus,k^\plus}$,
and rays $X^\pm_i \in (T_{p^\plusminus_i}M\setminus\{0\})/\IR$.
Let
$$\Sigma_{g,k^\minus,k^\plus}
  = \Sigmabar_{g,k^\minus,k^\plus}
    \setminus\{p^\minus_1,\ldots,p^\minus_{k^\minus}
               ,p^\plus_1,\ldots,p^\plus_{k^\plus}
             \} .
$$

Let $\eusJ_g$
denote the space of smooth conformal structures on
$\Sigmabar_{g,k^\minus,k^\plus}$.
Let $\eusD_{g,k}$ denote the space of
$k$-tuples
$$\Delta=(\Delta^\minus_1,\ldots,\Delta^\minus_{k^\minus},
  \Delta^\plus_1,\ldots,\Delta^\plus_{k^\plus})$$
of smooth closed pairwise disjoint disks on $\Sigmabar_{g,k^\minus,k^\plus}$
such that $p^\plusminus_i$ lies in the interior of $\Delta^\pm_i$.

It follows from Riemann's mapping theorem,
that for any $\bfj\in\eusJ_g$
and $\Delta\in\eusD_{g,k}$,
there exist unique $\bfj$-holomorphic coordinate functions
$z^\plusminus_i = x^\plusminus_i+i y^\plusminus_i
 : \Delta^\pm_i \to \{z\in\IC:|z|\le1\}$
such that $z^\plusminus_i(p^\plusminus_i)=0$
and $X^\pm_i$ is the ray spanned by the tangent vector
$\partial/\partial x^\plusminus_i$.
On $\Delta^\pm_i\setminus\{p^\plusminus_i\}$
we use cylindrical coordinates
$(t^\plusminus_i,\theta^\plusminus_i)$, with
$t^\plusminus_i+i\theta^\plusminus_i = - \log z^\plus_i$
and  $t^\minus_i+i\theta^\minus_ii = \log z^\minus_i$.
Then $(t^\plus_i,\theta^\plus_i)$ maps $\Delta^\plus_i\setminus\{p^\plus_i\}$
to
$[0,\infty)\times S^1$,
and $(t^\minus_i,\theta^\minus_i)$ maps
$\Delta^\minus_i\setminus\{p^\minus_i\}$ to
$(-\infty,0]\times S^1$.

Let $\eusJhat_{g,k^\minus,k^\plus}(A)$
be the smooth Fr\'echet space bundle over
$\eusJ_g\times\eusD_{g,k}$
whose fiber over the pair $(\bfj,\Delta)$
consists of all triples $(\bfj,\Delta,R)$ where
$R$ is a section of the bundle
$(T^{0,1}\Sigma_{g,k^\minus,k^\plus})^*\boxtimes TM$
over $\Sigma_{g,k^\minus,k^\plus}\times M$
such that
$$R = (dt^\plusminus_i - i\,d\theta^\plusminus_i)\otimes\nabla H$$
on each $\Delta^\pm_i\times M$.

We define the framed diffeomorphism group $\Diff_{g,k}$ as
the group of diffeomorphisms of $\Sigmabar_{g,k^\minus,k^\plus}$
that fix the points $p^\plusminus_i$ and the
rays $X^\pm_i$.
The group $\Diff_{g,k}$
acts on the spaces $\eusJ_g$,
$\eusD_{g,k}$
and $\eusJhat_{g,k^\minus,k^\plus}(A)$
by push-forward.
We let $\eusJhat_{g,k^\minus,k^\plus}^*(A)$
denote the largest subset of $\eusJhat_{g,k^\minus,k^\plus}(A)$
on which $\Diff_{g,k}$ acts freely.
We let
$$\eufJ_{g,k^\minus,k^\plus}(A)
=\eusJhat_{g,k^\minus,k^\plus}(A)/\Diff_{g,k}$$
and
$$\eufJ_{g,k^\minus,k^\plus}^*(A)
=\eusJhat^*_{g,k^\minus,k^\plus}(A)/\Diff_{g,k}.$$

There is an alternative description of $\eufJ^*_{g,k^\minus,k^\plus}(A)$.
Instead of fixing the rays $X^\pm_i$, we can make them part of the data.
Let $\eusJhat^+_{g,k^\minus,k^\plus}(A)$ denote the set of quadruples
$(\bfj,\Delta,X,R)$ where
$X=((X^-_1,\ldots,X^-_{k^\minus}),(X^+_1,\ldots,X^+_{k^\minus}))$
and $X^\pm_i\in (T_{p^\plusminus_i}M\setminus\{0\})/\IR$.
Let $\Diff^+_{g,k}$ be the group of
diffeomorphisms of $\Sigmabar_{g,k^\minus,k^\plus}$
that fix the points $p^\pm_i$.
Then
$$\eufJ_{g,k^\minus,k^\plus}(A)
  = \eusJhat^+_{g,k^\minus,k^\plus}(A)/\Diff^+_{g,k}.$$
Let $\eusJhat^{+*}_{g,k^\minus,k^\plus}(A)$ be the largest subset of
$\eusJhat^+_{g,k^\minus,k^\plus}(A)$ on which $\Diff^+_{g,k}$ acts freely.
Then
$$\eufJ^*_{g,k^\minus,k^\plus}(A)
  = \eusJhat^{+*}_{g,k^\minus,k^\plus}(A)/\Diff^+_{g,k}.$$
If the Hamiltonian $H$ is time independent,
then the section $R$ can be chosen independently of the rays
$X^\pm_i$.
Then the group $\Diff^+_{g,k}$ acts naturally on
$\eusJhat_{g,k^\minus,k^\plus}$, and we define
$$\eufJ^0_{g,k^\minus,k^\plus}(A)
  = \eusJhat_{g,k^\minus,k^\plus}(A)/\Diff^+_{g,k}.$$
Let $\eusJhat^{**}_{g,k^\minus,k^\plus}(A)$
denote the largest subset of $\eusJhat_{g,k^\minus,k^\plus}(A)$
on which $\Diff^+_{g,k}$ acts freely.
Then we define
$$\eufJ^{0*}_{g,k^\minus,k^\plus}(A)
  = \eusJhat^{**}_{g,k^\minus,k^\plus}(A)/\Diff^+_{g,k}.$$

\proclaim Proposition 2.2.1.

{
If $k=k^\minus+k^\plus>0$, then
$\eusJhat_{g,k^\minus,k^\plus}^*(A)
 = \eusJhat_{g,k^\minus,k^\plus}(A)$
and
$\eufJ_{g,k^\minus,k^\plus}^*(A)
 = \eufJ_{g,k^\minus,k^\plus}(A)$.
The space $\eusJhat_{g,k^\minus,k^\plus}^*(A)$
is a weakly contractible smooth Fr\'echet manifold.
The space $\eufJ_{g,k^\minus,k^\plus}^*(A)$
is a smooth Fr\'echet manifold.
The manifold $\eusJhat_{g,k^\minus,k^\plus}^*(A)$
is a smooth principal $\Diff_{g,k}$
bundle over the manifold $\eufJ^*_{g,k^\minus,k^\plus}(A)$.

If $H$ is time-independent,
then similar statements hold for
$\eusJhat^{**}_{g,k^\minus,k^\plus}(A)$,
$\Diff^+_{g,k}$,
and $\eufJ^{0*}_{g,k^\minus,k^\plus}(A)$.
}

It is tempting to say that $\eufJ^*_{g,k^\minus,k^\plus}(A)$
and $\eufJ^{0*}_{g,k^\minus,k^\plus}(A)$
are classifying spaces for the groups $\Diff_{g,k}$
and $\Diff^+_{g,k}$.
However, our arguments fall slightly short of this statement,
as we do not show that the spaces and group actions
are simplicial.

\demo Proof.

The space $\eusJhat_{g,k^\minus,k^\plus}(A)$
is clearly a smooth Fr\'echet manifold by definition.
The space $\eusJhat^*_{g,k^\minus,k^\plus}(A)$
is an open subset of $\eusJhat_{g,k^\minus,k^\plus}(A)$
and is hence a smooth Fr\'echet manifold.

The space $\eusJ_g$ is contractible by Prop.~2.1.3.
It is not hard to show,
using Riemann's mapping theorem,
that the space $\eusD_{g,k}$ is contractible.
The space of smooth sections $R$ is a vector space,
and is hence contractible.
It follows that $\eusJhat_{g,k^\minus,k^\plus}(A)$
is contractible.
The complement of $\eusJhat^*_{g,k^\minus,k^\plus}(A)$
essentially has infinite codimension.
A simple perturbation argument shows that
$\eusJhat^*_{g,k^\minus,k^\plus}(A)$
is weakly contractible.

Assume that $k>0$.
If a diffeomorphism $\phi\in\Diff_{g,k}$
fixes a conformal structure $\bfj$,
a disk $\Delta^\pm_i$ and the ray $X^\pm_i$,
it then follows from Riemann's mapping theorem that
the restriction of $\phi$ to $\Delta^\pm_i$ is the identity map.
By analytic continuation $\phi$ is then identity map.
Hence $\eusJhat_{g,k^\minus,k^\plus}(A) = \eusJhat_{g,k^\minus,k^\plus}^*(A)$.

We now turn to the main part of the theorem.
The action of $\Diff_{g,k}$ on $\eusJhat^*_{g,k^\minus,k^\plus}(A)$
is free by definition.
To show that the action defines a principal bundle over a smooth
Fr\'echet manifold,
we need to construct smooth local sections.
The tools for doing this are provided by Teichm\"uller theory.
For a global analysis approach to Teichm\"uller theory,
see [EE], [J] and [T].
Let $\Diff_{g,0}^0$ denote the identity component of $\Diff_{g,0}$.
The Teichm\"uller space $T_{g,0}$ is defined as the quotient
$\eusJ_g/\Diff^0_{g,0}$.
By Teichm\"uller's theorem,
$T_{g,0}$ is a finite dimensional complex manifold,
diffeomorphic to an open ball.
The action of $\Diff_{g,0}$ on $\eusJ_g$ is free,
except for $g=1$,
in which case the stabilizer of each point is $S^1\times S^1$,
and $g=0$,
in which case the stabilizer is the M\"obius group.
By the Earle-Eells theorem
the action of $\Diff^0_{g,0}$ on $\eusJ_g$ defines a smooth fiber bundle.

The mapping class group $\Gamma_g=\Diff_{g,0}/\Diff^0_{g,0}$
is a discrete group that acts on $T_{g,0}$.
The quotient $T_{g,0}/\Diff_{g,0} = T_{g,0}/\Gamma_g$ is
the moduli space of conformal structures on $\Sigmabar_{g,0,0}$.
This action of $\Gamma_g$ on $T_{g,0}$ is free
except at a discrete set of points
that have finite stabilizers,
and by Kravetz' theorem,
the action is properly discontinuous.

\subdemo The case $g\ge1$:

Let $\Diff^0_{g,k}=\Diff_{g,k}\cap\Diff^0_{g,0}$.
Then $\Diff_{g,k}/\Diff^0_{g,k}=\Gamma_g$.
Let $T_{g,k}=\eusJ_g/\Diff^0_{g,k}$.
It follows from Kravetz' theorem that the action of $\Gamma_k$
on $T_{g,k}$ is properly discontinuous.
In particular it admits local slices.
It follows from the Earl-Eells theorem that the action of $\Diff^0_{g,k}$
on $\eusJ_g$ defines a trivial smooth fiber bundle.
In particular the action admits local slices.
It follows that the action of $\Diff_{g,k}$ on $\eusJ_g$ admits
local slices.
Each point has a compact stabilizer.
It follows that the action of $\Diff_{g,k}$ on
$\eusJhat^*_{g,k^\minus,k^\plus}(A)$
admits local slices.
As this action is free, it gives a principal fiber bundle.

\subdemo The case $g=0$, $k\ge3$:

Select three of the points $p^\plusminus_i$.
Let $\Diff_{0,3}^*$ be the group of diffeomorphisms
that fix these three points,
but not necessarily the corresponding rays $X^\pm_i$.
The action of $\Diff_{0,3}^*$
on $\eusJ_0$ is free and transitive.
By the Earle-Eells theorem it gives a diffeomorphism between
$\Diff_{0,3}^*$ and $\eusJ_0$.
We can then argue as in the case $g\ge1$,
using the group $\Diff^*_{0,3}$ in place of $\Diff_{g,0}$.

\subdemo The case $g=0$, $k=1,2$:

If $k=2$,
then the quotient of $\eusJ_0$ by $\Diff_{0,2}$ is $S^1$,
and the stabilizer of any $\bfj\in\eusJ_0$ is $\IR$.
If $k=1$,
then the group $\Diff_{0,1}$ acts transitively on $\eusJ_0$,
and the stabilizer of any $\bfj\in\eusJ_0$
is the group of affine transformations of the form
$z\mapsto az+b$ with $a\in\IR$ and $b\in\IC$.
In either case we need to verify that
the action of the stabilizer of any $\bfj\in\eusJ_0$
on $\eusD_{0,k}$
defines a principal fiber bundle.
To show that a free action by a noncompact
finite dimensional group defines a principal bundle
one has to show that the orbits are closed.
It is not hard to see that the orbit of any $\Delta\in\eusD_{0,k}$
under these two groups is closed.

\subdemo The case $g=k=0$:

This case is handled the same way as the case $g=0$, $k=1,2$.
The theorem follows from the observation that the orbit of any nonzero
section $R$ of $(T^{0,1}\Sigmabar_{0,0,0})^*\boxtimes_\IC TM$
under the action of the M\"obius group is closed.
\enddemo

\subheading 2.3. The parametrized moduli spaces
  $\eufM^d_\sigma  ((\alpha^\minus_1,\ldots,\alpha^\minus_{k^\minus}),
                    (\alpha^\plus_1,\ldots,\alpha^\plus_{k^\plus})
                   )
  $.

Let $\eta$ be a smooth cut-off function with
$\eta=0$ on $(-\infty,{\ts{1\over3}}]$
and $\eta=1$ on $[{\ts{2\over3}},\infty)$.
We then define the energy of a map $u:\Sigma_{g,k^\minus,k^\plus}\to M$ as
$$E[u]
  = \iint_{\Sigma_{g,k^\minus,k^\plus}}
    \Bigl| du - \sum \eta(\pm t^\plusminus_i) \, d\theta^\plusminus_i
                     \otimes \bfJ\nabla H(u,\theta^\plusminus_i)
    \Bigr|^2 \, dA.
  \eqno\qbz
$$

\proclaim Definition 2.3.1.

A smooth map
$u:\Sigma_{g,k^\minus,k^\plus}\to M$ is a perturbed pseudoholomorphic curve
with data $(\bfj,\Delta,R)\in\eusJhat^*_{g,k^\minus,k^\plus}(A)$,
if and only if $E[u]<\infty$,
and,
in terms of local complex coordinates $z=x+iy$ compatible with $\bfj$,
$$(dx-i\,dy)\otimes\left({\partial u\over\partial x}
  + \bfJ\bigl(u(x,y)\bigr) {\partial u\over\partial y}\right)
  + R((x,y),u(x,y))= 0.
  \eqno\qbb
$$

We often write this equation $\dbar u + R(\cdot,u)=0$.
The left  hand side is a section of
$(T^{0,1}\Sigma_{g,k^\minus,k^\plus})^*\otimes_\IC TM$.
On the coordinate charts $(t^\plusminus_i,\theta^\plusminus_i)$
the equation takes the simpler form
$${\partial u\over\partial t^\plusminus_i}
  + \bfJ\bigl(u(t^\plusminus_i,\theta^\plusminus_i)\bigr)
    {\partial u\over\partial\theta^\plusminus_i}
  + \nabla H(u(t^\plusminus_i,\theta^\plusminus_i),\theta^\plusminus_i)
  = 0 .
$$

We define the moduli space $\eufM_{g,k^\minus,k^\plus}$
as the set of
equivalence classes $[u,c]$,
under the action of the group $\Diff_{g,k}$,
of pairs $(u,c)$,
where $u$ is a perturbed pseudoholomorphic curve
with data $c=(\bfj,\Delta,R)$.
It follows from the results of [F4],
see also [R] Prop.~3.1,
that if $[u,c]\in\eufM_{g,k^\minus,k^\plus}$,
then $u$ converges exponentially fast to a periodic orbit
$\alpha^\pm_i$ at each puncture $p^\plusminus_i$.
This gives a decomposition
$$\eufM_{g,k^\minus,k^\plus}
   = \bigcup_{\ss \alpha^\minus_1,\ldots,\alpha^\minus_{k^\minus}
              \in\eusC_H
              \atop
              \strut
              \ss \alpha^\plus_1,\ldots,\alpha^\plus_{k^\plus}\in\eusC_H
             } \mskip-20mu
     \eufM_{g,k^\minus,k^\plus}
          ((\alpha^\minus_1,\ldots,\alpha^\minus_{k^\minus}),
           (\alpha^\plus_1,\ldots,\alpha^\plus_{k^\plus})).
$$

Let
$$\eusP_{g,k^\minus,k^\plus}
  ((\alpha^\minus_1,\ldots,\alpha^\minus_{k^\minus_1}),
   (\alpha^\plus_1,\ldots,\alpha^\plus_{k^\plus_1})
  )
$$
denote the space of maps $u:\Sigma_{g,k^\minus,k^\plus}\to M$
such that $u$ is of Sobolev class $H^3$ locally,
and on each $\Delta^\pm_i$,
the map $u$ satisfies the same asymptotic conditions
as in the definition of $\eusP(\alpha^\minus,\alpha^\plus)$,
see [F4] and [R] \S2,
but with limit $\alpha^\pm_i$.
Thus defined,
$\eusP_{g,k^\minus,k^\plus}
   ((\alpha^\minus_1,\ldots,\alpha^\minus_{k^\minus_1}),
    (\alpha^\plus_1,\ldots,\alpha^\plus_{k^\plus_1})
  )
$
is a Hilbert manifold.

The tangent bundle
$T\eusP_{g,k^\minus,k^\plus}
  ((\alpha^\minus_1,\ldots,\alpha^\minus_{k^\minus_1}),
   (\alpha^\plus_1,\ldots,\alpha^\plus_{k^\plus_1})
  )
$
is the Hilbert space bundle with fibers
$$T_u\eusP_{g,k^\minus,k^\plus}
  ((\alpha^\minus_1,\ldots,\alpha^\minus_{k^\minus_1}),
   (\alpha^\plus_1,\ldots,\alpha^\plus_{k^\plus_1})
  )
  = H^3(\Sigma_{g,k^\minus,k^\plus},u^*TM).
$$
By Prop.~2.2.1,
$$\eqalign{
  & \eufP_{g,k^\minus,k^\plus}
    ((\alpha^\minus_1,\ldots,\alpha^\minus_{k^\minus_1}),
     (\alpha^\plus_1,\ldots,\alpha^\plus_{k^\plus_1})
    ) \cr
  &  \qquad = \eusP_{g,k^\minus,k^\plus}
  ((\alpha^\minus_1,\ldots,\alpha^\minus_{k^\minus_1}),
   (\alpha^\plus_1,\ldots,\alpha^\plus_{k^\plus_1})
  )
  \times_{\Diff_{g,k}}
  \eusJhat^*_{g,k^\minus,k^\plus}(A) \cr
          }
$$
is the total space of a smooth fiber bundle
$$\pi:\eufP_{g,k^\minus,k^\plus}
  ((\alpha^\minus_1,\ldots,\alpha^\minus_{k^\minus_1}),
   (\alpha^\plus_1,\ldots,\alpha^\plus_{k^\plus_1})
  )
  \to\eufJ^*_{g,k^\minus,k^\plus}(A)
  \eqno\qbt
$$
with fiber
$\eusP_{g,k^\minus,k^\plus}
  ((\alpha^\minus_1,\ldots,\alpha^\minus_{k^\minus_1}),
   (\alpha^\plus_1,\ldots,\alpha^\plus_{k^\plus_1})
  ).
$
Furthermore,
by Prop.~2.2.1,
$$\eqalign{
  & T^\fib \eufP_{g,k^\minus,k^\plus}
    ((\alpha^\minus_1,\ldots,\alpha^\minus_{k^\minus_1}),
     (\alpha^\plus_1,\ldots,\alpha^\plus_{k^\plus_1})
    ) \cr
  &\qquad  = T\eusP_{g,k^\minus,k^\plus}
    ((\alpha^\minus_1,\ldots,\alpha^\minus_{k^\minus_1}),
     (\alpha^\plus_1,\ldots,\alpha^\plus_{k^\plus_1})
    )
    \times_{\Diff_{g,k}}
    \eusJhat^*_{g,k^\minus,k^\plus}(A) . \cr
         }
$$
is a smooth Hilbert space bundle over
$\eufP_{g,k^\minus,k^\plus}
  ((\alpha^\minus_1,\ldots,\alpha^\minus_{k^\minus_1}),
   (\alpha^\plus_1,\ldots,\alpha^\plus_{k^\plus_1})
  )
$.
It is the fiber-wise tangent bundle of \qbt.
The fiber $T^\fib_{[u,c]}\eufP_{g,k^\minus,k^\plus}$
of $T^\fib\eufP_{g,k^\minus,k^\plus}$ at $[u,c]$ is
$$H^3(\Sigma_{g,k^\minus,k^\plus}, u^*TM)
.$$
Finally,
let
$\rmF_{g,k^\minus,k^\plus}$
be the smooth Hilbert space bundle over
$\eusP_{g,k^\minus,k^\plus}$
with fibers
$H^2((T^{0,1}\Sigma_{g,k^\minus,k^\plus})^*\otimes_\IC u^* TM)$.
By Prop.~2.2.1,
$$\IF_{g,k^\minus,k^\plus}
  = \rmF_{g,k^\minus,k^\plus}
    \times_{\Diff_{g,k}}
    \eusJhat^*_{g,k^\minus,k^\plus}(A)
$$
is a smooth Hilbert space bundle over
$\eufP_{g,k^\minus,k^\plus}
  ((\alpha^\minus_1,\ldots,\alpha^\minus_{k^\minus_1}),
   (\alpha^\plus_1,\ldots,\alpha^\plus_{k^\plus_1})
  ).
$
The fiber $(\IF_{g,k^\minus,k^\plus})_{[u,c]}$
of $\IF_{g,k^\minus,k^\plus}$ at $[u,c]$ is
$$H^2((T^{0,1}\Sigma_{g,k^\minus,k^\plus})^*\otimes_\IC u^* TM).$$

The left hand side of \qbb{} defines a section
$$\Psi_{g,k^\minus,k^\plus} : \eufP_{g,k^\minus,k^\plus}
  ((\alpha^\minus_1,\ldots,\alpha^\minus_{k^\minus_1}),
   (\alpha^\plus_1,\ldots,\alpha^\plus_{k^\plus_1})
  )
 \to \IF_{g,k^\minus,k^\plus} .
$$
Arguing as in [F4], see also [R] Prop.~3.1,
we see that the moduli space
$\eufM_{g,k^\minus,k^\plus}$
is contained in
$\eufP_{g,k^\minus,k^\plus}$
and is the zero locus of $\Psi_{g,k^\minus,k^\plus}$.

For
$[u,c]\in\eufM_{g,k^\minus,k^\plus}
  ((\alpha^\minus_1,\ldots,\alpha^\minus_{k^\minus_1}),
   (\alpha^\plus_1,\ldots,\alpha^\plus_{k^\plus_1})
  )
$,
we let
$$D_{[u,c]}\Psi_{g,k^\minus,k^\plus}:T_{[u,c]}\eufP_{g,k^\minus,k^\plus}
  ((\alpha^\minus_1,\ldots,\alpha^\minus_{k^\minus_1}),
   (\alpha^\plus_1,\ldots,\alpha^\plus_{k^\plus_1})
  )
 \to (\IF_{g,k^\minus,k^\plus})_{[u,c]}
$$
denote the intrinsic derivative of $\Psi_{g,k^\minus,k^\plus}$ at $[u,c]$.
We let
$$D^\fib_{[u,c]}\Psi_{g,k^\minus,k^\plus}:
T^\fib_{[u,c]}\eufP_{g,k^\minus,k^\plus}
  ((\alpha^\minus_1,\ldots,\alpha^\minus_{k^\minus_1}),
   (\alpha^\plus_1,\ldots,\alpha^\plus_{k^\plus_1})
  )
 \to (\IF_{g,k^\minus,k^\plus})_{[u,c]}$$
denote the restriction of $D_{[u,c]}\Psi_{g,k^\minus,k^\plus}$ to
$T^\fib_{[u,c]}\eufP_{g,k^\minus,k^\plus}
  ((\alpha^\minus_1,\ldots,\alpha^\minus_{k^\minus_1}),
   (\alpha^\plus_1,\ldots,\alpha^\plus_{k^\plus_1})
  )
$.
As in [R] Prop.~2.1,
$$\eqalign{
  & D^\fib_{[u,c]}\Psi_{g,k^\minus,k^\plus} \xi \cr
  & = (dx-i\,dy)\otimes \bigg(
      \nabla_x\xi
      +\bfJ\nabla_y\xi
      + \bfT\Big(\xi,{\partial u\over\partial x}\Big)
      + \bfJ\,\bfT\Big(\xi,{\partial u\over\partial y}\Big)
                     \bigg)
      + \nabla_\xi \nabla R . \cr
          }
  \eqno\qbc
$$
Arguing as in [F4], see also [R] Thm.~2.2,
we see that $D^\fib_{[u,c]}\Psi_{g,k^\minus,k^\plus}$ is Fredholm.
This gives a decomposition
$$\eufM_{g,k^\minus,k^\plus}
  ((\alpha^\minus_1,\ldots,\alpha^\minus_{k^\minus_1}),
   (\alpha^\plus_1,\ldots,\alpha^\plus_{k^\plus_1})
  )
  = \bigcup_{d\in\IZ} \eufM^d_{g,k^\minus,k^\plus}
  ((\alpha^\minus_1,\ldots,\alpha^\minus_{k^\minus_1}),
   (\alpha^\plus_1,\ldots,\alpha^\plus_{k^\plus_1})
  )
$$
according to index of $D^\fib_{[u,c]}\Psi_{g,k^\minus,k^\plus}$.
By taking \qbc{} as a definition,
we get a family $D^\fib\Psi_{g,k^\minus,k^\plus}$
of Fredholm operators $D^\fib_{[u,c]}\Psi_{g,k^\minus,k^\plus}$
parametrized by $\eufP_{g,k^\minus,k^\plus}$.
This gives a decomposition
$$\eufP_{g,k^\minus,k^\plus}
  ((\alpha^\minus_1,\ldots,\alpha^\minus_{k^\minus_1}),
   (\alpha^\plus_1,\ldots,\alpha^\plus_{k^\plus_1})
  )
 = \bigcup_{d\in\IZ} \eufP_{g,k^\minus,k^\plus}^d
  ((\alpha^\minus_1,\ldots,\alpha^\minus_{k^\minus_1}),
   (\alpha^\plus_1,\ldots,\alpha^\plus_{k^\plus_1})
  )
$$
according to the index of $D^\fib_{[u,c]}\Psi_{g,k^\minus,k^\plus}$.
We write $\IF^d_{g,k^\minus,k^\plus}$ and
$\Psi^d_{g,k^\minus,k^\plus}$ for the restrictions of
$\IF_{g,k^\minus,k^\plus}$
and $\Psi_{g,k^\minus,k^\plus}$ to $\eufP_{g,k^\minus,k^\plus}^d$.

\proclaim Proposition 2.3.2.

The moduli space
$\eufM^d_{g,k^\minus,k^\plus}
  ((\alpha^\minus_1,\ldots,\alpha^\minus_{k^\minus}),
   (\alpha^\plus_1,\ldots,\alpha^\plus_{k^\plus}))$
is an embedded smooth Fr\'echet submanifold of\/
$\eufP^d_{g,k^\minus,k^\plus}
  ((\alpha^\minus_1,\ldots,\alpha^\minus_{k^\minus}),
       (\alpha^\plus_1,\ldots,\alpha^\plus_{k^\plus}))$.
The projection
$$\pi:\eufM^d_{g,k^\minus,k^\plus}
      ((\alpha^\minus_1,\ldots,\alpha^\minus_{k^\minus}),
       (\alpha^\plus_1,\ldots,\alpha^\plus_{k^\plus}))
   \to\eufJ^*_{g,k^\minus,k^\plus}(A)
$$
is a $\sigma$-proper Fredholm map of index $d$.

\demo Proof.

Let $[u,c]\in\eufM^d_{g,k^\minus,k^\plus}$.
The annihilator of the range of $D^\fib_{[u,c]}\Psi^d_{g,k^\minus,k^\plus}$
is the kernel of the adjoint operator.
This is an elliptic differential operator.
By the Aronzajn unique continuation theorem,
a section in the annihilator of the range
can not vanish on an open subset of $\Sigma_{g,k^\minus,k^\plus}$.
On the other hand, by varying $R$ we see that
the range of $D_{[u,c]}\Psi^d_{g,k^\minus,k^\plus}$
contains all elements of $(\IF^d_{g,k^\minus,k^\plus})_{[u,c]}$
that are smooth and supported on
the complement of the disks $\Delta^\pm_i$.
Hence the annihilator of the range consists of sections supported on
the disks $\Delta^\pm_i$.
It follows that the annihilator
of the range of $D_{[u,c]}\Psi^d_{g,k^\minus,k^\plus}$ is trivial,
and that the range is dense.
As $D^\fib_{[u,c]}\Psi^d_{g,k^\minus,k^\plus}$ is Fredholm,
it is then surjective.

Let $U\subset T^\fib_{[u,c]}\eufP^d_{g,k^\minus,k^\plus}$
be the kernel of $D^\fib_{[u,c]}\Psi^d_{g,k^\minus,k^\plus}$.
Let $U^\perp$ be a linear complement to $U$ in
$T^\fib_{[u,c]}\eufP^d_{g,k^\minus,k^\plus}$.
Let $W$ be a linear complement to
$T^\fib_{[u,c]}\eufP^d_{g,k^\minus,k^\plus}$ in
$T_{[u,c]}\eufP^d_{g,k^\minus,k^\plus}$.
As $D_{[u,c]}\Psi^d_{g,k^\minus,k^\plus}$ is surjective,
there exists a finite dimensional subspace $V$ of $W$
such that
$D_{[u,c]}\Psi^d_{g,k^\minus,k^\plus}:U^\perp\oplus V\to
(\IF^d_{g,k^\minus,k^\plus})_{[u,c]}$
is invertible.
Let $V^\perp$ be a linear complement to $V$ in $W$.
We can then parametrize a neighborhood of $[u,c]$
in $\eufP^d_{g,k^\minus,k^\plus}$ by a neighborhood
of $0$ in
$T_{[u,c]}\eufP^d_{g,k^\minus,k^\plus}
 =U\oplus U^\perp\oplus V\oplus V^\perp
$.
By the implicit function theorem,
applied to the Hilbert space $U^\perp\oplus V$
with the Fr\'echet space $U\oplus V^\perp$
as parameter space,
there exist smooth maps
$S_1: U\oplus V^\perp \to U^\perp$
and
$S_2: U\oplus V^\perp \to V$
such that the map
$u\oplus v^\perp \mapsto u+v^\perp
  +S_1(u\oplus v^\perp)+S_2(u\oplus V^\perp)$
gives a local parametrization of a neighborhood of $[u,c]$
in $\eufM^d_{g,k^\minus,k^\plus}$.
Hence $\eufM^d_{g,k^\minus,k^\plus}$ is an embedded Fr\'echet submanifold
of $\eufP^d_{g,k^\minus,k^\plus}$.

We can parametrize a neighborhood of $[c]$ in $\eufJ^*_{g,k^\minus,k^\plus}(A)$
by a neighborhood of 0 in $W=V\oplus V^\perp$.
We may assume that the projection
$\eufP^d_{g,k^\minus,k^\plus}\to\eufJ^*_{g,k^\minus,k^\plus}(A)$
is given,
in terms of the parametrizations,
by the projection $U\oplus U^\perp\oplus V\oplus V^\perp\to V\oplus V^\perp$.
In terms of these parametrizations $\pi$ is the map
$$U\oplus V^\perp \to V\oplus V^\perp$$
given by
$$u\oplus v^\perp \mapsto S_2(u\oplus v^\perp) + v^\perp.$$
Hence
$$\index D\pi
  = \dim U - \dim V
  = \index D^\fib_{[u,c]}\Psi^d_{g,k^\minus,k^\plus}
  = d .
$$

Finally, that the map $\pi$ is $\sigma$-proper follows by Uhlenbeck-Gromov
compactness.
\enddemo

\proclaim Definition 2.3.3.

For any smooth $\sigma:[0,1]^q\to\eufJ^*_{g,k^\minus,k^\plus}(A)$
the parametrized moduli space
$$\eufM_\sigma^d
  ((\alpha^\minus_1,\ldots,\alpha^\minus_{k^\minus}),
   (\alpha^\plus_1,\ldots,\alpha^\plus_{k^\plus})
  )
$$
is defined as the set of pairs
$$(u,x)\in
 \eufM^{d-q}_{g,k^\minus,k^\plus}
  ((\alpha^\minus_1,\ldots,\alpha^\minus_{k^\minus}),
   (\alpha^\plus_1,\ldots,\alpha^\plus_{k^\plus})
  )
  \times [0,1]^q
$$
such that $\pi(u)=\sigma(x)$.

\proclaim Definition 2.3.4.

{
A map $\sigma:[0,1]^q\to\eufJ^*_{g,k^\minus,k^\plus}(A)$
is regular if and only if $\sigma$,
and the restrictions of $\sigma$ to the boundary faces of $[0,1]^q$
of any dimension,
are transverse to
$\pi:\eufM^d_{g,k^\minus,k^\plus}
  ((\alpha^\minus_1,\ldots,\alpha^\minus_{k^\minus}),
   (\alpha^\plus_1,\ldots,\alpha^\plus_{k^\plus})
  )
  \to \eufJ^*_{g,k^\minus,k^\plus}(A)
$
for all integers $d$ and all $\alpha^\pm_*\in\eusC_H$.

It is semi-regular if and only if $\sigma$ is transverse
to $\pi$ for $d\le-q$,
and its restriction to any boundary face of $[0,1]^q$
is transverse to $\pi$ for $d\le-q'$
where $q'$ is the dimension of the boundary face.
}

If $\sigma:[0,1]^q\to\eufJ^*_{g,k^\minus,k^\plus}(A)$
is regular,
then
$\eufM^d_\sigma$
is empty for $d<0$,
and for $d\ge0$ it is an embedded submanifold of
$\eufP^d_{g,k^\minus,k^\plus}$
of dimension $d$ with corners.
If $\sigma:[0,1]^q\to\eufJ^*_{g,k^\minus,k^\plus}(A)$
is semi-regular,
then
$\eufM^d_\sigma$
is empty for $d<0$,
and it is discrete set of points for $d=0$.

\proclaim Remark 2.3.5.

\rm
There is an alternative description of $\eufM^d_\sigma$.
Let
$\eufP^d_\sigma=\sigma^* \eufP^{d-q}_{g,k^\minus,k^\plus}$.
Then $\eufP^d_\sigma$ is a fiber bundle over $[0,1]^q$ with fiber
$\eusP^{d-q}_{g,k^\minus,k^\plus}$.
Let $\IF^d_\sigma=\sigma^*\IF^{d-q}$.
Then $\IF^d_\sigma$ is a Hilbert space bundle over $\eufP^d_\sigma$.
Let $\Psi^d_\sigma = \sigma^* \Psi^{d-q}_{g,k^\minus,k^\plus}$.
Then $\eufM^d_\sigma$ is the zero locus of $\Psi^d_\sigma$.

\proclaim Proposition 2.3.6.

A map $\sigma:[0,1]^q\to\eufJ^*_{g,k^\minus,k^\plus}(A)$
is transverse to
$\pi:\eufM^{d-q}_{g,k^\minus,k^\plus}
  ((\alpha^\minus_1,\ldots,\alpha^\minus_{k^\minus}),
   (\alpha^\plus_1,\ldots,\alpha^\plus_{k^\plus})
  )
  \to \eufJ^*_{g,k^\minus,k^\plus}(A)
$
if and only if the section $\Psi^d_\sigma$ is transverse to the
zero section of $\IF^d_\sigma$.

\demo Proof.

The map
$\sigma:[0,1]^q\to\eufJ^*_{g,k^\minus,k^\plus}(A)$
is transverse to the map
$\pi:\eufM^{d-q}_{g,k^\minus,k^\plus}\to\eufJ^*_{g,k^\minus,k^\plus}(A)$
if and only if the lifted map
$\sigma:\eufP^d_\sigma\to\eufP^{d-q}_{g,k^\minus,k^\plus}$
is transverse to
$\eufM^{d-q}_{g,k^\minus,k^\plus}$.
This is the case if and only if the lifted map
$\sigma:\IF^d_\sigma\to\IF^{d-q}_{g,k^\minus,k^\plus}$ is transverse
to the intersection of the section $\Psi^{d-q}_{g,k^\minus,k^\plus}$
and the zero section of $\IF^{d-q}_{g,k^\minus,k^\plus}$.
This map and these two manifolds are pairwise transverse.
Hence the above condition says that the triple is transverse.
But the triple is transverse if and only if the pull-backs of the section
$\Psi^{d-q}_{g,k^\minus,k^\plus}$
and the zero section of $\IF^{d-q}_{g,k^\minus,k^\plus}$ under $\sigma$
are transverse sections of $\IF^d_\sigma$.
The pull-backs are $\Psi^d_\sigma$ and the zero section of $\IF^d_\sigma$.
\enddemo

\proclaim Remark 2.3.7.

\rm
Note that
$$\eusM^d(\alpha^\minus,\alpha^\plus)
  = \eufM^d_{\sigma_0} ((\alpha^\minus), (\alpha^\plus))
$$
and
$$\eusP^d(\alpha^\minus,\alpha^\plus)
  = \eufP^d_{\sigma_0} ((\alpha^\minus), (\alpha^\plus))
$$
where $\sigma_0:[0,1]^0\to\eufJ^*_{0,1,1}$ maps the standard
0-simplex to the
standard conformal structure on $\Sigma_{0,1,1}=\IR\times S^1$
and the standard perturbation $R=(dt-i\,d\theta)\otimes\nabla H$.
For the sake of definiteness,
we can take $\Delta^\minus_1$ and $\Delta^\plus_1$ to be the
closures of $(-\infty,-1]\times S^1$ and $[1,\infty)\times S^1$ in
$\Sigmabar_{0,1,1}$.
Thus the results about
$\eufM^d_\sigma$ and $\eufP^d_\sigma$
that we will show in \S3
apply, mutatis mutandis,
to
$\eusM^d(\alpha^\minus,\alpha^\plus)$
and
$\eusP^d(\alpha^\minus,\alpha^\plus)$.

%%%%%%%%%%%%%%%%%%%%%%%%%%%%%%%%%%%%%%%%%%%%%%%%%%%%%%%%%%%%%%

\heading \S3. Excision

\subheading 3.1. The gluing maps $\gop^\ell_{ij}$ and $\cop^\ell_{ij}$.

For any $\ell\ge0$ we define a smooth map
$$\gop_{ij}^\ell :
  \eufJ^*_{g_1,k^\minus_1,k^\plus_1+1}(A)
  \times
  \eufJ^*_{g_2,k^\minus_2+1,k^\plus_2}(A)
  \to
  \eufJ_{g_1+g_2,k^\minus_1+k^\minus_2,k^\plus_1+k^\plus_2}(A)
  \eqno\qca
$$
as follows.
The definition is illustrated in fig.~1.
We denote the punctures, disks, and coordinates
on $\Sigmabar_{g_1,k^\minus_1,k^\plus_1+1}$
by $p^\plusminus_{1,i}$, $\Delta^\pm_{1,i}$,
and $(t^\plusminus_{1,i},\theta^\plusminus_{1,i})$.
We denote the punctures, disks, and coordinates
on $\Sigmabar_{g_2,k^\minus_2+1,k^\plus_2}$
by $p^\plusminus_{2,i}$, $\Delta^\pm_{2,i}$,
and $(t^\plusminus_{2,i},\theta^\plusminus_{2,i})$.
We can then form a surface
$\Sigmabar^\ell_{g_1+g_2,k^\minus_1+k^\minus_2,k^\plus_1+k^\plus_2}$
by gluing $p^\plus_{1,i}$ to $p^\minus_{2,j}$ as follows.
We remove the disk $t^\plus_{1,i}>\ell$ from
$\Sigmabar_{g_1,k^\minus_1,k^\plus_1+1}$
and the disk $t^\minus_{2,j}<-\ell$ from
$\Sigmabar_{g_2,k^\minus_2+1,k^\plus_2}$.
We then identify the annulus $0\le t^\plus_{1,i}\le\ell$ on
$\Sigmabar_{g_1,k^\minus_1,k^\plus_1+1}$
with the annulus $-\ell\le t^\minus_{2,j}\le0$ on
$\Sigmabar_{g_2,k^\minus_2+1,k^\plus_2}$
by identifying
$(t^\plus_{1,i},\theta_{1,i})$
with
$(t^\minus_{2,j},\theta^\minus_{2,j})$
if
$t^\minus_{2,j}=t^\plus_{1,i}-\ell$
and
$\theta^\minus_{2,j}=\theta^\plus_{1,i}$.
On the corresponding annulus on
$\Sigmabar^\ell_{g_1+g_2,k^\minus_1+k^\minus_2,k^\plus_1+k^\plus_2}$
we use coordinates
$(t_0,\theta_0)$
defined by
$t_0=t^\plus_{1,i}=t^\minus_{2,j}+\ell$
and $\theta_0=\theta^\plus_{1,i}=\theta^\minus_{2,j}$.
We refer to these annuli on
$\Sigmabar_{g_1,k^\minus_1,k^\plus_1+1}$,
$\Sigmabar_{g_2,k^\minus_2+1,k^\plus_2}$ and
$\Sigmabar^\ell_{g_1+g_2,k^\minus_1+k^\minus_2,k^\plus_1+k^\plus_2}$
as the {\it necks}.
We order the punctures on
$\Sigmabar^\ell_{g_1+g_2,k^\minus_1+k^\minus_2,k^\plus_1+k^\plus_2}$
 as follows:
$$\eqalign{
  & \big( (p^\minus_{2,1},\ldots,p^\minus_{2,j-1},
           p^\minus_{1,1},\ldots,p^\minus_{1,k^\minus_1},
           p^\minus_{2,j+1},\ldots,p^\minus_{2,k^\minus_2+1}), \cr
  & \qquad
          (p^\plus_{1,1},\ldots,p^\plus_{1,i-1},
           p^\plus_{2,1},\ldots,p^\plus_{2,k^\plus_2},
           p^\plus_{i,j+1},\ldots,p^\plus_{1,k^\plus_1+1})
    \big). \cr
          }
$$
Then
$\Sigmabar^\ell_{g_1+g_2,k^\minus_1+k^\minus_2,k^\plus_1+k^\plus_2}$
is diffeomorphic to the fixed surface
$\Sigmabar_{g_1+g_2,k^\minus_1+k^\minus_2,k^\plus_1+k^\plus_2}$.
Choose a diffeomorphism
$\phi^\ell:
 \Sigmabar^\ell_{g_1+g_2,k^\minus_1+k^\minus_2,k^\plus_1+k^\plus_2}
 \to \Sigmabar_{g_1+g_2,k^\minus_1+k^\minus_2,k^\plus_1+k^\plus_2}$.

Let
$(\bfj_1,\Delta_1,R_1)\in\eusJhat^*_{g_1,k^\minus_1,k^\plus_1+1}(A)$
and
$(\bfj_2,\Delta_2,R_2)\in\eusJhat^*_{g_2,k^\minus_2+1,k^\plus_2}(A).$
Identify the complement of the neck in
$\Sigmabar^\ell_{g_1+g_2,k^\minus_1+k^\minus_2,k^\plus_1+k^\plus_2}$
with appropriate parts of the complement of the necks in
$\Sigmabar_{g_1,k^\minus_1,k^\plus_1+1}
 \cup
 \Sigmabar_{g_2,k^\minus_2+1,k^\plus_2}
$.
Let $\bfj_3$ be the conformal structure on
$\Sigmabar^\ell_{g_1+g_2,k^\minus_1+k^\minus_2,k^\plus_1+k^\plus_2}$
that equals $\bfj_1$ and $\bfj_2$ off the neck,
and equals the conformal structure
given by the coordinates $(t_0,\theta_0)$ on the neck.
Let
$$\eqalign{
  \Delta_3 = & \big( (\Delta^\minus_{2,1},\ldots,\Delta^\minus_{2,j-1},
     \Delta^\minus_{1,1},\ldots,\Delta^\minus_{1,k^\minus_1},
     \Delta^\minus_{2,j+1},\ldots,\Delta^\minus_{2,k^\minus_2+1}), \cr
  & \qquad
    (\Delta^\plus_{1,1},\ldots,\Delta^\plus_{1,i-1},
     \Delta^\plus_{2,1},\ldots,\Delta^\plus_{2,k^\plus_2},
     \Delta^\plus_{i,j+1},\ldots,\Delta^\plus_{1,k^\plus_1+1})
               \big) \cr
          }
$$
Let $R_3$ be equal to $R_1$ and $R_2$ off the neck,
and equal to $(dt_0-i\,d\theta_0)\otimes\nabla H$ on the neck.
Then
$(\phi^\ell_*\,\bfj_3,\phi^\ell_*\Delta_3,\phi^\ell_*R_3)
  \in\eusJhat_{g_1+g_2,k^\minus_1+k^\minus_2,k^\plus_1+k^\plus_2}(A)
$.
The corresponding equivalence class
$[\phi^\ell_*\,\bfj_3,\phi^\ell_*\Delta_3,\phi^\ell_*R_3]
  \in \eufJ_{g_1+g_2,k^\minus_1+k^\minus_2,k^\plus_1+k^\plus_2}(A)
$
only depends on
$[\bfj_1,\Delta_1,R_1] \in \eufJ^*_{g_1,k^\minus_1,k^\plus_1+1}(A)$
and
$[\bfj_2,\Delta_2,R_2] \in \eufJ^*_{g_2,k^\minus_2+1,k^\plus_2}(A)$.
We define \qca{} by
$$[\bfj_1,\Delta_1,R_1]
  \gop^\ell_{ij}
  [\bfj_2,\Delta_2,R_2]
  =
  [\phi^\ell_*\,\bfj_3,\phi^\ell_*\Delta_3,\phi^\ell_*R_3].
$$
In particular, given smooth maps
$$\sigma_1 : [0,1]^{q_1} \to \eufJ^*_{g_1,k^\minus_1,k^\plus_1+1}(A)$$
and
$$\sigma_2 : [0,1]^{q_2}\to\eufJ^*_{g_2,k^\minus_2+1,k^\plus_2}(A),$$
we can form the map
$$\sigma_1\gop^\ell_{ij} \sigma_2 :
    [0,1]^{q_1+q_2}
    \to \eufJ_{g_1+g_2,k^\minus_1+k^\minus_2,k^\plus_1+k^\plus_2}(A) .
$$
defined as
$$(\sigma_1\gop^\ell_{ij} \sigma_2) (x_1,\ldots,x_{q_1+q_2})
  = \sigma_1(x_1,\ldots,x_{q_1})
    \gop^\ell_{ij} \sigma_2 (x_{q_1+1},\ldots,x_{q_1+q_2}) .
$$

We obtain a similar map
$$\cop^\ell_{ij} :
  \eufJ^*_{g,k^\minus+1,k^\plus+1}(A)
  \to \eufJ_{g+1,k^\minus,k^\plus}(A)
  \eqno\qcaa
$$
by gluing $p^\plus_i$ to $p^\minus_j$
on $\Sigmabar_{g,k^\minus+1,k^\plus+1}$
to form a surface $\Sigmabar^\ell_{g+1,k^\minus,k^\plus}$.
The details are the same as in the definition of $\gop^\ell_{ij}$.
We order the punctures on $\Sigmabar^\ell_{g+1,k^\minus,k^\plus}$ as follows:
$$\bigl((p^\minus_1,\ldots,
         p^\minus_{j-1},p^\minus_{j+1},\ldots,
         p^\minus_{k^\minus+1}
        )
        ,(p^\plus_1,\ldots,
          p^\plus_{i-1},p^\plus_{i+1},\ldots,
          p^\plus_{k^\plus+1}
         )
  \bigr) .
$$
In particular, given a smooth map
$$\sigma : [0,1]^q
  \to \eufJ^*_{g,k^\minus+1,k^\plus+1}(A),
$$
we can form the map
$$\cop^\ell_{ij}\sigma : [0,1]^q
  \to \eufJ_{g+1,k^\minus,k^\plus}(A).
$$

\subheading 3.2. The excision maps $\bfg^\ell_{ij}$ and $\bfc^\ell_{ij}$.

In this section we will construct a lift
$$\bfg^\ell_{ij} :
  \eufP^\plus_{g_1,k^\minus_1,k^\plus_1+1;i}[\alpha_0]
  \times
  \eufP^\minus_{g_2,k^\minus_2+1,k^\plus_2;j}[\alpha_0]
  \to
  \eufP_{g_1+g_2,k^\minus_1+k^\minus_2,k^\plus_1+k^\plus_2}
  \times
  \eusP(\alpha_0,\alpha_0)
$$
of the map $\qca$,
where
$$\eqalign{
  & \eufP^{\plus}_{g_1,k^\minus_1,k^\plus_1+1;i}[\alpha_0] \cr
  & = \eufP_{g_1,k^\minus_1,k^\plus_1+1}
      ((\alpha^\minus_{1,1},\ldots,\alpha^\minus_{1,k^\minus_1}),
       (\alpha^\plus_{1,1},\ldots,
        \alpha^\plus_{1,i-1},\alpha_0,\alpha^\plus_{1,i+1},\ldots
        \alpha^\plus_{1,k^\plus_1+1}
       )
      ) , \cr
  \noalign{\smallskip}
  & \eufP^{\minus}_{g_2,k^\minus_2+1,k^\plus_2;j}[\alpha_0] \cr
  & = \eufP_{g_2,k^\minus_2+1,k^\plus_2}
      ((\alpha^\minus_{2,1},\ldots,
        \alpha^\minus_{2,j-1},\alpha_0,\alpha^\minus_{2,j+1},\ldots
        \alpha^\minus_{2,k^\minus_2+1}),
       (\alpha^\plus_{2,1},\ldots,
        \alpha^\plus_{2,k^\plus_2}
       )
      ) \cr
          }
$$
and
$$\eqalign{
  & \eufP_{g_1+g_2,k^\minus_1+k^\minus_2,k^\plus_1+k^\plus_2} \cr
  & = \eufP_{g_1+g_2,k^\minus_1+k^\minus_2,k^\plus_1+k^\plus_2}
      ((\alpha^\minus_{2,1},\ldots,\alpha^\minus_{2,j-1},
        \alpha^\minus_{1,1},\ldots,\alpha^\minus_{1,k^\minus_1},
        \alpha^\minus_{2,j+1},\ldots,\alpha^\minus_{2,k^\minus_2+1}
       ), \cr
  & \phantom{{} = \eufP_{g_1+g_2,k^\minus_1+k^\minus_2,k^\plus_1+k^\plus_2} (
            }
       (\alpha^\plus_{1,1},\ldots,\alpha^\plus_{1,i-1},
        \alpha^\plus_{2,1},\ldots,\alpha^\plus_{2,k^\plus_2},
        \alpha^\plus_{1,i+1},\ldots,\alpha^\plus_{1,k^\plus_1+1}
       ) . \cr
         }
$$
We use similar notation for the bundles $T^\fib\eufP$ and $\IF$  over
these spaces.
The map $\bfg^\ell_{ij}$ will only be defined on an open subset
of
$\eufP^\plus_{g_1,k^\minus_1,k^\plus_1+1;i}[\alpha_0]
 \times
 \eufP^\minus_{g_2,k^\minus_2+1,k^\plus_2;j}[\alpha_0]
$.

Let $\Sigmabar^\ell_{g_1+g_2,k^\minus_1+k^\minus_2,k^\plus_1+k^\plus_2}$
be the surface obtained by gluing
$p^\plus_i$ on $\Sigmabar_{g_1,k^\minus_1,k^\plus_1+1}$
to $p^\minus_j$ on $\Sigmabar_{g_2,k^\minus_2+1,k^\plus_2}$
as in \S3.1.
We also get a copy of $\IR\times S^1$
by forming the union of the punctured disk $t^\plus_{1,i}\ge0$ in
$\Sigma_{g_1,k^\minus_1,k^\plus_1+1}$
and the punctured disk $t^\minus_{2,j}\le0$ in
$\Sigma_{g_2,k^\minus_2+1,k^\plus_2}$
and identifying
$(t^\plus_{1,i},\theta_{1,i})$
with $(t^\minus_{2,j},\theta^\minus_{2,j})$
if $t^\minus_{2,j}=t^\plus_{1,i}-\ell$
and $\theta^\minus_{2,j}=\theta^\plus_{1,i}$.
On $\IR\times S^1$ we have global coordinates $(t,\theta)$
defined by $t=t^\plus_{1,i}$ and $\theta=\theta^\plus_{1,i}$,
or $t=t^\minus_{2,j}+\ell$ and $\theta=\theta^\minus_{2,j}$.
Let $\eta$ be a smooth function $\IR\to[0,1]$
with $\eta=0$ on $(-\infty,{\ts{1\over3}}]$ and
$\eta=1$ on $[{\ts{2\over3}},\infty)$.

The domain of $\bfg^\ell_{ij}$ consists of all
$$([u_1,c_1],[u_2,c_2]) \in
  \eufP^\plus_{g_1,k^\minus_1,k^\plus_1+1;i}[\alpha_0]
  \times
  \eufP^\minus_{g_2,k^\minus_2+1,k^\plus_2;j}[\alpha_0] ,
$$
such that
$$[c_1]\gop^\ell_{ij}[c_2]
 \in \eufJ^*_{g_1+g_2,k^\minus_1+k^\minus_2,k^\plus_1+k^\plus_2}(A) ,
$$
$\dist\bigl(
    u_1(t^\plus_{1,i},\theta^\plus_{1,i}),
    \alpha_0(\theta^\plus_{1,i})
       \bigr) \le r_0/2
$
for $t^\plus_{1,i}\ge\ell/3$,
and
$\dist\bigl(
    u_2(t^\minus_{2,j},\theta^\minus_{2,j}),
    \alpha_0(\theta^\minus_{2,j})
       \bigr) \le r_0/2
$
for $t^\minus_{2,j}\le-\ell/3$.
Then
$$u_1(t^\plus_{1,i},\theta^\plus_{1,i})
  = \exp \xi_1(t^\plus_{1,i},\theta^\plus_{1,i})
$$
for $t^\plus_{1,i}\ge\ell/3$
and
$$u_2(t^\minus_{2,j},\theta^\minus_{2,j})
  = \exp \xi_2(t^\minus_{2,j},\theta^\minus_{2,j})
$$
for $t^\minus_j\le-\ell/3$
as in [R] \S2.
We define
$$\bfg^\ell_{ij}([u_1,c_1],[u_2,c_2])
  = ([u_3,c_1\gop^\ell_{ij}c_2],u_4)$$
where $u_3$ and $u_4$ are as follows.
Identify the complement of the necks in
$\Sigma_{g_1,k^\minus_1,k^\plus_1+1}\cup\Sigma_{g_2,k^\minus_2+1,k^\plus_2}$
with the complement of the necks in
$\Sigma^\ell_{g_1+g_2,k^\minus_1+k^\minus_2,k^\plus_1+k^\plus_2}
 \cup(\IR\times S^1)$.
Then we let $u_3$ and $u_4$ be equal to
to $u_1$ and $u_2$ off the necks.
Identify the necks of
$\Sigma_{g_1,k^\minus_1,k^\plus_1+1}$, $\Sigma_{g_2,k^\minus_2+1,k^\plus_2}$,
$\Sigma^\ell_{g_1+g_2,k^\minus_1+k^\minus_2,k^\plus_1+k^\plus_2}$,
and $(\IR\times S^1)$.
Then we let
$u_3=\exp_3$ and $u_4=\exp\xi_4$ on the necks,
where
$$\xi_3  = \Bigl(\cos\bigl(\pi\eta(t_0/\ell)\bigr)
                 +\bfJ\sin\bigl(\pi\eta(t_0/\ell)\bigr)
          \Bigr)
         \cos\bigl(\pi\eta(t_0/\ell)/2\bigr) \xi_1
         + \sin\bigl(\pi\eta(t_0/\ell)/2\bigr) \xi_2
$$
and
$$\xi_4  = - \Bigl(\cos\bigl(\pi\eta(t_0/\ell)\bigr)
                 +\bfJ\sin\bigl(\pi\eta(t_0/\ell)\bigr)
           \Bigr)
           \sin\bigl(\pi\eta(t_0/\ell)/2\bigr) \xi_1
         + \cos\bigl(\pi\eta(t_0/\ell)/2\bigr) \xi_2 .
$$
The point of this definition is that
$$\pmatrix{ \bigl(\cos(\pi s)
                 +i\sin(\pi s)
            \bigr)
            \cos(\pi s/2)
          & \sin(\pi s/2) \cr
            - \bigl(\cos(\pi s)
                    +i\sin(\pi s)
              \bigr)
              \sin(\pi s/2)
          & \cos(\pi s/2) \cr
          }
    \qquad 0\le s\le 1
  \eqno\qcb
$$
is a path from $\bigl({1\atop0}{0\atop1}\bigr)$
and $\bigl({0\atop1}{1\atop0}\bigr)$ in the Lie group
$\U(2)$.

The map $\bfg^\ell_{ij}$ lifts to bundle maps
$$\eqalign{
  (\bfg^\ell_{ij})_*: {}
  &  T^\fib\eufP^\plus_{g_1,k^\minus_1,k^\plus_1+1;i}[\alpha_0]
    \boxplus
    T^\fib\eufP^\minus_{g_2,k^\minus_2+1,k^\plus_2;j}[\alpha_0] \cr
  &  \to
    T^\fib\eufP_{g_1+g_2,k^\minus_1+k^\minus_2,k^\plus_1+k^\plus_2}
    \boxplus
    T\eusP(\alpha_0,\alpha_0) \cr
          }
$$
and
$$(\bfg^\ell_{ij})_*:
    \IF^\plus_{g_1,k^\minus_1,k^\plus_1+1;i}[\alpha_0]
    \boxplus
    \IF^\minus_{g_2,k^\minus_2+1,k^\plus_2;j}[\alpha_0]
    \to
    \IF_{g_1+g_2,k^\minus_1+k^\minus_2,k^\plus_1+k^\plus_2}
    \boxplus
    \rmF(\alpha_0,\alpha_0)$$
defined as follows.
Let
$([u_1,c_1],[u_2,c_2])$ be in the domain of $\bfg^\ell_{ij}$.
Let
$s_1\oplus s_2$
be an element of
$$
 T^\fib_{[u_1,c-1]}\eufP^\plus_{g_1,k^\minus_1,k^\plus_1+1;i}[\alpha_0]
 \boxplus
 T^\fib_{[u_2,c_2]}\eufP^\minus_{g_2,k^\minus_2+1,k^\plus_2;j}[\alpha_0]
$$
or
$$(\IF^\plus_{g_1,k^\minus_1,k^\plus_1+1;i}[\alpha_0])_{[u_1,c_1]}
 \boxplus
 (\IF^\minus_{g_2,k^\minus_2+1,k^\plus_2;j}[\alpha_0])_{[u_2,c_2]}
.$$
Then
$s_1 = (\exp \xi_1)_*\zeta_1$
and
$s_2 = (\exp \xi_2)_*\zeta_2$,
as in [R] \S2.
We define $(\bfg^\ell_{ij})_*(s_1\oplus s_2)=s_3\oplus s_4$
where $s_3$ and $s_4$ are as follows.
Identify the complement of the necks in
$\Sigma_{g_1,k^\minus_1,k^\plus_1+1}\cup\Sigma_{g_2,k^\minus_2+1,k^\plus_2}$
with the complement of the necks in
$\Sigma_{g_1+g_2,k^\minus_1+k^\minus_2,k^\plus_1+k^\plus_2}
 \cup(\IR\times S^1)$
We then let $s_3\oplus s_4$ be equal to $s_1\oplus s_2$ off the necks.
Identify the necks of
$\Sigma_{g_1,k^\minus_1,k^\plus_1+1}$,
$\Sigma_{g_2,k^\minus_2+1,k^\plus_2}$,
$\Sigma_{g_1+g_2,k^\minus_1+k^\minus_2,k^\plus_1+k^\plus_2}$
and
$\IR\times S^1$.
We then let
$s_3=(\exp\xi_3)_*\zeta_3$
and
$s_4=(\exp\xi_4)_*\zeta_4$
on the necks,
where
$$\zeta_3  = \Bigl(\cos\bigl(\pi\eta(t_0/\ell)\bigr)
                 +J\sin\bigl(\pi\eta(t_0/\ell)\bigr)
          \Bigr)
         \cos\bigl(\pi\eta(t_0/\ell)/2\bigr) \zeta_1
         + \sin\bigl(\pi\eta(t_0/\ell)/2\bigr) \zeta_2
$$
and
$$\zeta_4  = - \Bigl(\cos\bigl(\pi\eta(t_0/\ell)\bigr)
                 +J\sin\bigl(\pi\eta(t_0/\ell)\bigr)
           \Bigr)
           \sin\bigl(\pi\eta(t_0/\ell)/2\bigr) \zeta_1
         + \cos\bigl(\pi\eta(t_0/\ell)/2\bigr) \zeta_2.
$$

The map $\bfg^\ell_{ij}$ induces a map
$$\bfg^\ell_{ij} :
  \eufP^\plus_{\sigma_1;i}[\alpha_0]
  \times
  \eufP^\minus_{\sigma_2;j}[\alpha_0]
  \to
  \eufP_{\sigma_1\gop^\ell_{ij}\sigma_2}
  \times
  \eusP(\alpha_0,\alpha_0)
$$
the obvious way,
where
$$\eqalign{
  & \eufP^{\plus}_{\sigma_1;i}[\alpha_0]
    = \eufP_{\sigma_1}
      ((\alpha^\minus_{1,1},\ldots,\alpha^\minus_{1,k^\minus_1}),
       (\alpha^\plus_{1,1},\ldots,
        \alpha^\plus_{1,i-1},\alpha_0,\alpha^\plus_{1,i+1},\ldots
        \alpha^\plus_{1,k^\plus_1+}
       )
      ) , \cr
  \noalign{\smallskip}
  & \eufP^{\minus}_{\sigma_2;j}[\alpha_0]
    = \eufP_{\sigma_2}
      ((\alpha^\minus_{2,1},\ldots,
        \alpha^\minus_{2,j-1},\alpha_0,\alpha^\minus_{2,j+1},\ldots
        \alpha^\minus_{2,k^\minus_2+1}),
       (\alpha^\plus_{2,1},\ldots,
        \alpha^\plus_{2,k^\plus_2}
       )
      ) \cr
          }
$$
and
$$\eqalign{
  & \eufP_{\sigma_1\gop^\ell_{ij}\sigma_2} \cr
  & = \eufP_{\sigma_1\gop^\ell_{ij}\sigma_2}
      ((\alpha^\minus_{2,1},\ldots,\alpha^\minus_{2,j-1},
        \alpha^\minus_{1,1},\ldots,\alpha^\minus_{1,k^\minus_1},
        \alpha^\minus_{2,j+1},\ldots,\alpha^\minus_{2,k^\minus_2+1}
       ), \cr
  & \phantom{{} = \eufP_{\sigma_1\gop^\ell_{ij}\sigma_2}(
            }
       (\alpha^\plus_{1,1},\ldots,\alpha^\plus_{1,i-1},
        \alpha^\plus_{2,1},\ldots,\alpha^\plus_{2,k^\plus_2},
        \alpha^\plus_{1,i+1},\ldots,\alpha^\plus_{1,k^\plus_1+1}
       ) . \cr
         }
$$
This map is invertible.
The inverse of the matrices \qcb{} form a path
$$\pmatrix{ \bigl(\cos(\pi s)
                 -i\sin(\pi s)
            \bigr)
            \cos(\pi s/2)
          & - \bigl(\cos(\pi s)
                    -i\sin(\pi s)
              \bigr)
              \sin(\pi s/2) \cr
            \sin(\pi s/2)
          & \cos(\pi s/2) \cr
          }
    \qquad 0\le s\le 1 .
$$
Hence the inverse of $\bfg^\ell_{ij}$ is given by
$(\bfg^\ell_{ij})^{-1}(u_3,u_4)
  = ( u_1,u_2)
$
where $u_1$ and $u_2$ are equal $u_3$ or $u_4$ off the necks,
and
$u_1=\exp\xi_1$ with
$$\eqalign{
  \xi_1 = & \Bigl(\cos\bigl(\pi \eta(t_0/\ell)\bigr)
                 -\bfJ\sin\bigl(\pi \eta(t_0/\ell)\bigr)
            \Bigr)
            \cos\bigl(\pi \eta(t_0/\ell)/2\bigr) \xi_3 \cr
        &  - \Bigl(\cos\bigl(\pi \eta(t_0/\ell)\bigr)
                    -\bfJ\sin\bigl(\pi \eta(t_0/\ell)\bigr)
              \Bigr)
              \sin\bigl(\pi \eta(t_0/\ell)/2\bigr) \xi_4 \cr
          }
$$
and $u_2=\exp\xi_2$ with
$$\xi_2 = \sin\bigl(\pi \eta(t_0/\ell)/2\bigr)\xi_3 + \cos\bigl(\pi
\eta(t_0/\ell)/2\bigr)\xi_4
$$
on the necks.

This map lifts to bundle maps
$$(\bfg^\ell_{ij})_*:
    T^\fib\eufP^\plus_{\sigma_1;i}[\alpha_0]
    \boxplus
    T^\fib\eufP^\minus_{\sigma_2;j}[\alpha_0]
    \to
    T^\fib\eufP_{\sigma_1\gop^\ell_{ij}\sigma_2}
    \boxplus
    T\eusP(\alpha_0,\alpha_0)$$
and
$$(\bfg^\ell_{ij})_*:
    \IF^\plus_{\sigma_1;i}[\alpha_0]
    \boxplus
    \IF^\minus_{\sigma_2;j}[\alpha_0]
    \to
    \IF_{\sigma_1\gop^\ell_{ij}\sigma_2}
    \boxplus
    \rmF(\alpha_0,\alpha_0)$$
defined the same way as before.
These bundle maps are invertible.
We denote the inverse $(\bfg^\ell_{ij})^*$.

Let
$$\eqalign{
  \eufP^{\plus,\minus}_{g_1,k^\minus+1,k^\plus+1;i,j}[\alpha_0,\alpha_0]
  = \eufP_{g,k^\minus+1,k^\plus+1}
      (&(\alpha^\minus_1,\ldots,
        \alpha^\minus_{j-1},\alpha_0,\alpha^\minus_{j+1},\ldots
        \alpha^\minus_{k^\minus+1}), \cr
 &     (\alpha^\plus_1,\ldots,
        \alpha^\plus_{i-1},\alpha_0,\alpha^\plus_{i+1},\ldots
        \alpha^\plus_{k^\plus+1}
       )
      ) \cr
        }
$$
and
$$\eqalign{
  \eufP_{g+1,k^\minus,k^\plus}
  = \eufP_{g+1,k^\minus,k^\plus}
    (
  &  (\alpha^\minus_1,\ldots,\alpha^\minus_{j-1},
         \alpha^\minus_{j+1},\ldots,\alpha^\minus_{k^\minus+1}
     ), \cr
  &  (\alpha^\plus_1,\ldots,\alpha^\plus_{i-1},
      \alpha^\plus_{i+1},\ldots,\alpha^\plus_{k^\plus+1}
     )
    ) . \cr
         }
$$
Then the map $\qcaa$ lifts to a map
$$\bfc^\ell_{ij} :
  \eufP^{\plus,\minus}_{g,k^\minus+1,k^\plus+1;i,j}[\alpha_0,\alpha_0]
  \to
  \eufP_{g+1,k^\minus,k^\plus}
  \times
  \eusP(\alpha_0,\alpha_0),
$$
defined the same way as $\bfg^\ell_{ij}$.
The map $\bfc^\ell_{ij}$ lifts to bundle maps
$$(\bfc^\ell_{ij})_* :
  T^\fib\eufP^{\plus,\minus}_{g,k^\minus+1,k^\plus+1;i,j}[\alpha_0,\alpha_0]
  \to
  T^\fib\eufP_{g+1,k^\minus,k^\plus}
  \boxplus
  T\eusP(\alpha_0,\alpha_0)
$$
and
$$(\bfc^\ell_{ij})_* :
  \IF^{\plus,\minus}_{g,k^\minus+1,k^\plus+1;i,j}[\alpha_0,\alpha_0]
  \to
  \IF_{g+1,k^\minus,k^\plus}
  \boxplus
  \rmF(\alpha_0,\alpha_0) .
$$
defined the same way as $(\bfg^\ell_{ij})_*$.
Similarly,
for any smooth $\sigma:[0,1]^q\to\eufJ^*_{g,k^\minus+1,k^\plus+1}(A)$
there are local diffeomorphisms
$$\bfc^\ell_{ij} :
  \eufP^{\plus,\minus}_{\sigma;i,j}[\alpha_0,\alpha_0]
  \to
  \eufP_{\cop^\ell_{ij}\sigma}
  \times
  \eusP(\alpha_0,\alpha_0)
$$
with lifts
$$(\bfc^\ell_{ij})_* :
  T^\fib\eufP^{\plus,\minus}_{\sigma;i,j}[\alpha_0,\alpha_0]
  \to
  T^\fib\eufP_{\cop^\ell_{ij}\sigma}
  \boxplus
  T\eusP(\alpha_0,\alpha_0)
$$
and
$$(\bfc^\ell_{ij})_* :
  \IF^{\plus,\minus}_{\sigma;i,j}[\alpha_0,\alpha_0]
  \to
  \IF_{\cop^\ell_{ij}\sigma}
  \boxplus
  \rmF(\alpha_0,\alpha_0) .
$$
We denote the inverse of $(\bfc^\ell_{ij})_*$ by
$(\bfc^\ell_{ij})^*$.

\subheading 3.3. Index calculations.

The following is an excision principle for the index of
$D^\fib\Psi_{g,k^\minus,k^\plus}$.
Recall that $\eufP^d_{g,k^\minus,k^\plus}$ denotes the subset
of $\eufP_{g,\minus,k^\plus}$ where index $D^\fib\Psi_{g,k^\minus,k^\plus}$
is $d$.
Similarly we add a superscript $d$ to the notation introduced in \S3.2.

\proclaim Lemma 3.3.1.

$$\eqalign{
    \bfg^\ell_{ij}
    \bigl( \eufP^{d_1;\plus}_{g_1,k^\minus_1,k^\plus_1+1;i}[\alpha_0]
           \times
           \eufP^{d_2;\minus}_ {g_2,k^\minus_2+1,k^\plus_2;j}[\alpha_0]
    \bigr)
  & \subset
    \eufP^{d_1+d_2}_{g_1+g_2,k^\minus_1+k^\minus_2,k^\plus_1+k^\plus_2}
    \times
    \eusP^0(\alpha_0,\alpha_0) \cr
  \bfc^\ell_{ij}
    \bigl( \eufP^{d;\plus,\minus}_{g,k^\minus+1,k^\plus+1;i,j}
           [\alpha_0,\alpha_0]
    \bigr)
  & \subset
    \eufP^d_{g+1,k^\minus,k^\plus}
    \times
    \eusP^0(\alpha_0,\alpha_0) \cr
  \bfg^\ell_{ij}
    \bigl( \eufP^{d_1;\plus}_{\sigma_1;i}[\alpha_0]
           \times
           \eufP^{d_2;\minus}_ {\sigma_2;j}[\alpha_0]
    \bigr)
  & \subset
    \eufP^{d_1+d_2}_{\sigma_1\gop^\ell_{ij}\sigma_2}
    \times
    \eusP^0(\alpha_0,\alpha_0) \cr
  \bfc^\ell_{ij}
    \bigl( \eufP^{d;\plus,\minus}_{\sigma;i,j}
           [\alpha_0,\alpha_0]
    \bigr)
  & \subset
    \eufP^d_{\cop^\ell_{ij}\sigma}
    \times
    \eusP^0(\alpha_0,\alpha_0) \cr
         }
$$

\demo Proof.

By Remark 2.3.5, it suffices to prove the first two statements.
We will prove the first one.
The second is proven the same way.
Let
$([u_3,c_1\gop^\ell_{ij}c_2],u_4) = \bfg^\ell_{ij}([u_1,c_1],[u_2,c_2])$.
To keep the notation simple,
we write $u_1$, $u_2$, and $u_3$ for
$[u_1,c_1]$, $[u_2,c_2]$, and $[u_3,c_1\gop^\ell_{ij}c_2]$.
We need to show that
$$\index D^\fib_{u_3}\Psi_{g_1+g_2,k^\minus_1+k^\minus_2,k^\plus_1+k^\plus_2}
  = \index D^\fib_{u_1}\Psi_{g_1,k^\minus_1,k^\plus_1+1}
    + \index D^\fib_{u_2}\Psi_{g_2,k^\minus_2+1,k^\plus_2}
$$
and
$$\index D_{u_4}\Psi(\alpha_0,\alpha_0)=0.$$
The map $u_4$ is a perturbation of $\alphabar_0$.
Hence $\index D_{u_4}\Psi(\alpha_0,\alpha_0)=0$.
Thus it suffices to show that
$$\eqalign{
  & \index D^\fib_{u_1}\Psi_{g_1,k^\minus_1,k^\plus_1+1}
    + \index D^\fib_{u_2}\Psi_{g_2,k^\minus_2+1,k^\plus_2} \cr
  & = \index D^\fib_{u_3}\Psi_{g_1+g_2,k^\minus_1+k^\minus_2,
                               k^\plus_1+k^\plus_2}
      + \index D_{u_4}\Psi(\alpha_0,\alpha_0). \cr
          }
$$
This is simply invariance of the index under excision.
However, as we are working on non-compact surfaces,
we can not appeal to the standard excision principle of
Atiyah and Singer.
Instead we argue as follows.
It follows from [R] Prop.~2.1 that
$$\eqalign{
  & (D^\fib_{u_1}\Psi_{g_1,k^\minus_1,k^\plus_1+1})(\exp\xi_1)_*\zeta_1
    \oplus
     (D^\fib_{u_2}\Psi_{g_2,k^\minus_2+1,k^\plus_2})(\exp\xi_2)_*\zeta_2 \cr
  \noalign{\smallskip}
  & = (\exp\xi_1)_*\left(
        {\partial\zeta_1\over\partial t}
        +  \bfJ(D_{\alpha_0}\Phi_H)\zeta_1
        + (\xi_1\otimes\zeta_1)\cdot f
                   \right) \cr
  & \phantom{={}}
    \oplus (\exp\xi_2)_*\left(
        {\partial\zeta_2\over\partial t}
        +  \bfJ(D_{\alpha_0}\Phi_H)\zeta_2
        + (\xi_2\otimes\zeta_2)\cdot f
                   \right) \cr
          }
$$
on the necks.
Similarly,
$$\eqalign{
  & (D^\fib_{u_3}\Psi_{g_1+g_2,k^\minus_1+k^\minus_2,k^\plus_1+k^\plus_2})
    (\exp\xi_3)_*\zeta_3
    \oplus
     (D_{u_4}\Psi(\alpha_0,\alpha_0))(\exp\xi_4)_*\zeta_4 \cr
  \noalign{\smallskip}
  & = (\exp\xi_3)_*\left(
        {\partial\zeta_3\over\partial t}
        +  \bfJ(D_{\alpha_0}\Phi_H)\zeta_3
        + (\xi_3\otimes\zeta_3)\cdot f
                   \right) \cr
  & \phantom{={}}
    \oplus (\exp\xi_4)_*\left(
        {\partial\zeta_4\over\partial t}
        +  \bfJ(D_{\alpha_0}\Phi_H)\zeta_4
        + (\xi_4\otimes\zeta_4)\cdot f
                   \right)  \cr
          }
$$
on the necks.
It follows that the excision map $\bfg^\ell_{ij}$ intertwines
$D^\fib_{u_1}\Psi_{g_1,k^\minus_1,k^\plus_1+1}
 \oplus D^\fib_{u_2}\Psi_{g_2,k^\minus_2+1,k^\plus_2}$
and
$D^\fib_{u_3}
   \Psi_{g_1+g_2,k^\minus_1+k^\minus_2,k^\plus_1+k^\plus_2}
 \oplus
 D_{u_4}\Psi(\alpha_0,\alpha_0)
$
modulo lower order terms that are supported on the necks.
These terms define a compact operator.
Hence
$$\eqalign{
  & \index \bigl( D^\fib_{u_1}\Psi_{g_1,k^\minus_1,k^\plus_1+1}
                \oplus D^\fib_{u_2}\Psi_{g_2,k^\minus_2+1,k^\plus_2}
           \bigr) \cr
  & = \index ( \bfg^\ell_{ij})^*
               \bigl( D^\fib_{u_3}\Psi_{g_1+g_2,k^\minus_1+k^\minus_2,
                                        k^\plus_1+k^\plus_2}
                      \oplus
                      D_{u_4}\Psi(\alpha_0,\alpha_0)
               \bigr) \cr
  & = \index \bigl( D^\fib_{u_3}\Psi_{g_1+g_2,k^\minus_1+k^\minus_2,
                                      k^\plus_1+k^\plus_2}
                    \oplus
                    D_{u_4}\Psi(\alpha_0,\alpha_0)
             \bigr) . \cr
         }
$$
\enddemo

In our setup two different minimal Chern numbers arise.
The minimal Chern number of Floer is $N_0$.

\proclaim Definition 3.3.2.

We define the homotopical minimal Chern number $N_0$ by
$$\langle \c_1(TM) | \pi_2(M) \rangle = N_0\IZ$$
and the homological minimal Chern number $N_1$ by
$$\langle \c_1(TM) | H_2(M,\IZ) \rangle = N_1\IZ.$$

The following is the main index formula.

\proclaim Proposition 3.3.3.

{
There exist a unique function
$$\muH : \eusC^0_H\to \IZ/2N_0\IZ$$
such that if
$\alpha^\pm_i\in\eusC^0_H$
and
$\eufP^d_{g,k^\plus,k^\minus}
  ((\alpha^\minus_1,\ldots,\alpha^\minus_{k^\minus}),
   (\alpha^\plus_1,\ldots,\alpha^\plus_{k^\plus})
   )
$
is nonempty, then
$$d
  \equiv 2n(1-g-k^\minus)
    + \sum_{i=1}^{k^\minus} \muH(\alpha^\minus_i)
    - \sum_{i=1}^{k^\plus} \muH(\alpha^\plus_i)
  \quad \bmod 2N_1,
$$
and, for $g=0$,
$$d
  \equiv 2n(1-k^\minus)
    + \sum_{i=1}^{k^\minus} \muH(\alpha^\minus_i)
    - \sum_{i=1}^{k^\plus} \muH(\alpha^\plus_i)
  \quad \bmod 2N_0 .
$$
The function $\muH$ is independent of $\bfJ$.

If $\sigma$ is a smooth map
$[0,1]^q\to\eufJ^*_{g,k^\minus,k^\plus}(A)$,
$\alpha^\pm_i\in\eusC^0_H$, and
$\eufP^d_\sigma
  ((\alpha^\minus_1,\ldots,\alpha^\minus_{k^\minus}),\break
   (\alpha^\plus_1,\ldots,\alpha^\plus_{k^\plus})
   )
$
is nonempty, then
$$d
  \equiv q+2n(1-g-k^\minus)
    + \sum_{i=1}^{k^\minus} \muH(\alpha^\minus_i)
    - \sum_{i=1}^{k^\plus} \muH(\alpha^\plus_i)
  \quad \bmod 2N_1 .
$$
and, for $g=0$,
$$d
  \equiv q+2n(1-k^\minus)
    + \sum_{i=1}^{k^\minus} \muH(\alpha^\minus_i)
    - \sum_{i=1}^{k^\plus} \muH(\alpha^\plus_i)
  \quad \bmod 2N_0 .
$$

If $\alpha^\pm\in\eusC^0_H$
and $\eusP^d(\alpha^\minus,\alpha^\plus)$ is non-empty, then
$$d\equiv\mu(\alpha^\minus)-\mu(\alpha^\plus) \quad \bmod 2N_0.$$
}

\demo Proof.

For each $\alpha\in\eusC^0_H$,
we choose
$[u^\plus_\alpha,c^\plus_\alpha]\in\eufP_{0,0,1}(\emptyset,(\alpha))$.
To keep the notation simple we write $u$ for $[u,c]$.
In particular we write $u^\plus_\alpha$ for $[u^\plus_\alpha,c^\plus_\alpha]$.
If there exists a function $\muH$ as in the Proposition,
then
$$\index D^\fib_{u^\plus_\alpha}\Psi_{0,0,1}
 \equiv 2n - \muH(\alpha)
 \quad \bmod 2N_0 .
 \eqno\qcc
$$
This shows uniqueness.

To show existence,
we take \qcc{} as the definition of $\muH(\alpha)$,
and show that the first index formula in the Proposition holds.
We first note that if $u^0\in\eufP_{g,0,0}$,
then $D^\fib_{u^0}\Psi_{g,0,0}$
is a perturbed $\dbar$-complex on the compact Riemann surface
$\Sigma_{g,0,0}$ coupled to $(u^0)^*TM$.
By the Riemann-Roch theorem,
$$\index D^\fib_{u^0}\Psi_{g,0,0}
  = 2n(1-g) + 2\big\langle \c_1(TM) \big| (u^0)_*[\Sigma]\rangle.
$$
In particular,
$$\index D^\fib_{u^0}\Psi_{g,0,0}
  \equiv 2n(1-g) \quad \bmod  2N_1.
  \eqno\qcd
$$
If $g=0$, then this and the following congruences hold modulo $2N_0$.
We also note that by Prop.~2.1.7,
$$\index D_\alphabar\Psi(\alpha,\alpha) = 0 .
  \eqno \qce
$$

For each $\alpha\in\eusC^0_H$,
we also choose
$u^\minus_\alpha\in\eufP_{0,1,0}((\alpha),\emptyset)$.
Let
$(u^0,\alphabar) = \bfg^\ell_{11}(u^\plus_\alpha,u^\minus_\alpha).$
With a slight abuse of notation we here write ``$\,\alphabar\,$'' for a
a map that is only approximately equal to $\alphabar$.
By Lemma 3.3.1,
$$\index D^\fib_{u^\plus_\alpha}\Psi_{0,0,1}
 + \index D^\fib_{u^\minus_\alpha}\Psi_{0,1,0}
 = \index D^\fib_{u^0}\Psi_{0,0,0}
   + \index D_\alphabar\Psi(\alpha,\alpha).$$
It then follows from \qcc{}, \qcd{} and \qce{} that
$$\index D^\fib_{u^\minus_\alpha}\Psi_{0,1,0}
  \equiv \muH(\alpha) \quad \bmod 2N_1 .
  \eqno\qcf
$$
Let $u\in\eufP_{g,k^\minus,k^\plus}
  ((\alpha^\minus_1,\ldots,\alpha^\minus_{k^\minus}),
   (\alpha^\plus_1,\ldots,\alpha^\plus_{k^\plus})
   )
$.
By repeated use of Lemma 3.3.1,
$$\eqalign{
  & \index D^\fib_u\Psi_{g,k^\minus,k^\plus}
    + \index D^\fib_{u^\plus_{\alpha^\minus_1}} \Psi_{0,0,1}
    + \cdots
    + \index D^\fib_{u^\plus_{\alpha^\minus_{k^\minus}}} \Psi_{0,0,1} \cr
  & \qquad
    + \index D^\fib_{u^\minus_{\alpha^\plus_1}} \Psi_{0,1,0}
    + \cdots
    + \index D^\fib_{u^\minus_{\alpha^\plus_{k^\plus}}} \Psi_{0,1,0} \cr
  & = \index D^\fib_{u^0} \Psi_{g,0,0}
    + \index D_{\alphabar^\minus_1}\Psi(\alpha^\minus_1,\alpha^\minus_1)
    + \cdots
    + \index D_{\alphabar^\minus_{k^\minus}}
             \Psi(\alpha^\minus_{k^\minus},\alpha^\minus_{k^\minus}) \cr
  & \qquad
    + \index D_{\alphabar^\plus_1}\Psi(\alpha^\plus_1,\alpha^\plus_1)
    + \cdots
    + \index D_{\alphabar^\plus_{k^\plus}}
             \Psi(\alpha^\plus_{k^\plus},\alpha^\plus_{k^\plus})
 . \cr
          }
$$
The first index formula now follows from \qcc, \qcd, \qce,
and \qcf.
The other index formulas are easy consequences of the first index formula;
see Remark 2.3.5 and 2.3.7.
It follows from Prop.~2.1.3 by continuity
that $\muH$ does not depend on $\bfJ$.
\enddemo

A different approach to index formulas for perturbed pseudoholomorphic
curves uses the Maslov index, see [F1], [SZ] and [MS2] Sect.~9.2.

\subheading 3.4. Coherent orientations.

We will orient the moduli space by performing excisions
on the orientation sheaves.
Our construction is inspired by [D] Sect.~3, [DK] Sect.~7.1.6, and [FH].
However, the details of our construction are different,
as we do not localize the kernels and cokernels.

Let $T$ be a Fredholm operator.
We then let $\eusO_T$ denote the set of orientations of
the vector space
$\ker T \oplus (\coker T)^*$,
or, equivalently, the set of orientations of the line
$\Lambda^\max\ker T \otimes (\Lambda^\max\coker T)^*$.
Then $\eusO_T$ is a 1-dimensional vector space over $\IZ/2\IZ$.
The following properties of $\eusO_T$ will be used extensively.

\item $\bullet$

If $T$ is invertible, then $\eusO_T=\IZ/2\IZ$
canonically.

\item $\bullet$

If $T-T'$ is a compact operator,
then $\eusO_T=\eusO_{T'}$ canonically.

\item $\bullet$

$\eusO_{T_1\circ T_2}=\eusO_{T_1}\otimes\eusO_{T_2}$
canonically.

\item $\bullet$

Hence,
if $S_1T-US_2$ is compact,
and $S_1$ and $S_2$ are invertible,
then we can identify
$$\eqalign{
  & \eusO_T
    = (\IZ/2\IZ)\otimes\eusO_T
    = \eusO_{S_1}\otimes\eusO_T
    = \eusO_{S_1\circ T} \cr
  & = \eusO_{U\circ S_2}
    = \eusO_U\otimes\eusO_{S_2}
    = \eusO_U\otimes(\IZ/2\IZ)
    = \eusO_U. \cr
          }
$$
We denote the resulting isomorphism
$$S_* : \eusO_T \to \eusO_U.$$

\item $\bullet$

We can identify $\eusO_{T_1\oplus T_2}$ with
$\eusO_{T_1}\otimes\eusO_{T_2}$.
In other words,
elements of $\eusO_{T_1}$ and $\eusO_{T_2}$
determine an element of $\eusO_{T_1\oplus T_2}$.
We use the following convention.
Choose orientations of $\ker T_1$, $\coker T_1$, $\ker T_2$ and $\coker T_2$
compatible with the orientations of
$\ker T_1\oplus (\coker T_1)^*$
and $\ker T_1 \oplus (\coker T_2)^*$.
We then use the induced orientation
of
$$(\ker T_1) \oplus (\coker T_1)^* \oplus (\ker T_2) \oplus (\coker T_2)^*.$$

\item $\bullet$

The spaces $\eusO_{T_1}\otimes\eusO_{T_2}$ and
$\eusO_{T_2}\otimes\eusO_{T_1}$ are identical.
However,
if we choose elements of $\eusO_{T_1}$ and $\eusO_{T_2}$,
and make the identifications
$\eusO_{T_1}\otimes\eusO_{T_2}=\eusO_{T_1\oplus T_2}$
and
$\eusO_{T_2}\otimes\eusO_{T_1}=\eusO_{T_2\oplus T_1}$
as above,
then the induced elements of
$\eusO_{T_1}\otimes\eusO_{T_2}$
and
$\eusO_{T_2}\otimes\eusO_{T_1}$
differ by $(-1)^{\index T_1 \index T_2}$.

\item $\bullet$

More generally,
if $\rho$ is a permutation on $k$ letters,
and we choose elements of $\eusO_{T_1},\ldots,\eusO_{T_k}$,
then the induced elements of
$\eusO_{T_1}\otimes\cdots\eusO_{T_k}$
and
$\eusO_{T_{\rho(1)}}\otimes\cdots\eusO_{T_{\rho(k)}}$
differ by the graded sign of the permutation
$$(T_1,\ldots,T_k) \to (T_{\rho(1)},\ldots,T_{\rho(k)})$$
where
$T_i$ is assinged the degree $\index T_i$.

\noindent
The graded sign of a permutation is defined as the sign of the permutation
obtained by removing all elements of even degree.

\item $\bullet$

If $T=\{T_x\}_{x\in X}$ is a continuous family of Fredholm operators
parametrized by a topological space $X$,
then $\eusO_T=\bigcup_{x\in X}\eusO_{T_x}$ forms a continuous double
cover of $X$. We call $\eusO_T$ the orientation sheaf for the family $T$.

Let $\eufO_{g,k^\minus,k^\plus}$
denote the orientation sheaf for
the family $D^\fib\Psi_{g,k^\minus,k^\plus}$
parametrized by
$\eufP_{g_,k^\minus,k^\plus}
 ((\alpha^\minus_1,\ldots,\alpha^\minus_{k^\minus}),
  (\alpha^\plus_1,\ldots,\alpha^\plus_{k^\plus})
 )
$.
Let $\eufO_\sigma$
denote the orientation sheaf for the family $D\Psi_\sigma$
parametrized by
$\eufP_\sigma
 ((\alpha^\minus_1,\ldots,\alpha^\minus_{k^\minus}),
  (\alpha^\plus_1,\ldots,\alpha^\plus_{k^\plus})
 )
$.
Let $\eusO(\alpha^\minus,\alpha^\plus)$
denote the orientation sheaf for the family
$D\Psi(\alpha^\minus,\alpha^\plus)$
parametrized by $\eusP(\alpha^\minus,\alpha^\plus)$.
Let $\eusO_{[0,1]^q}$ denote the orientation  sheaf for the manifold $[0,1]^q$.

A global section of $\eufO_{g,k^\minus,k^\plus}$
determines a global section of
$$\eufO_\sigma=\eusO_{[0,1]^q}\otimes\sigma^*\eufO_{g,k^\minus,k^\plus},$$
see Remark 2.3.5,
and a global section of
$$\eusO(\alpha^\minus,\alpha^\plus)=\eufO_{\sigma_0},$$
see Remark 2.3.7.
If the moduli space $\eufM_\sigma$ is regular,
then a global section of $\eufO_\sigma$ determines an orientation
of $\eufM_\sigma$; see for instance [DK] Sect.~5.4.1.
Hence it suffices to choose global sections of the orientation sheaves
$\eufO_{g,k^\minus,k^\plus}$.

We first note that in some cases $\eufO_{g,k^\minus,k^\plus}$
has a canonical section.
If $[u,c]\in\eufP_{g,0,0}(\emptyset,\emptyset)$,
then $D_{[u,c]}\Psi_{g,0,0}$ is a perturbed $\dbar$-operator coupled to
a complex vector bundle over the compact surface $\Sigma_{g,0,0}$.
Hence the operator has a complex index line bundle.
The real index line bundle is the determinant of the complex index line bundle.
Hence the real index line bundle has a canonical orientation.
This gives a canonical section of $\eufO_{g,0,0}$.

By Prop.~2.1.7,
the operator $D_{\alphabar_0}\Psi(\alpha_0,\alpha_0)$ is invertible.
Thus
$\eusO_{\alphabar_0}(\alpha_0,\alpha_0)$
has a canonical element.
By continuity,
we also get canonical elements of $\eusO_u(\alpha_0,\alpha_0)$
for $u$ close to $\alphabar_0$.

In the proof of Lemma 3.3.1 we saw that the maps
$(\bfg^\ell_{ij})_*$
intertwine the families
$$D^\fib\Psi_{g_1,k^\minus_1,k^\plus_1+1}
  \boxplus
  D^\fib\Psi_{g_2,k^\minus_2+1,k^\plus_2}
$$
and
$$D^\fib\Psi_{g_1+g_2,k^\minus_1+k^\minus_2,k^\plus_1+k^\plus_2}
  \boxplus
  D\Psi(\alpha_0,\alpha_0)
$$
up to compact operators.
Hence there are induced excision maps
$$(\bfg^\ell_{ij})_* :
  \eufO_{g_1,k^\minus_1,k^\plus_1+1}
  \boxtimes
  \eufO_{g_2,k^\minus_2+1,k^\plus_2}
  \to
  \eufO_{g_1+g_2,k^\minus_1+k^\minus_2,k^\plus_1+k^\plus_2}
  \otimes
  \eusO_{\alphabar_0}(\alpha_0,\alpha_0).
  \eqno\qcg
$$
Here we write $\eusO_{\alphabar_0}(\alpha_0,\alpha_0)$
for $\eusO_u(\alpha_0,\alpha_0)$ with $u$ close to $\alphabar_0$.

Similarly there are excision maps
$$(\bfc^\ell_{ij})_* :
  \eufO_{g,k^\minus+1,k^\plus+1}
  \to
  \eufO_{g+1,k^\minus,k^\plus}
  \otimes
  \eusO_{\alphabar_0}(\alpha_0,\alpha_0).
  \eqno\qch
$$

\proclaim Definition 3.4.1.

{
A coherent system of orientations is a
choice of a global section of the double cover
$\eufO_{g,k^\minus,k^\plus}$
of each
$\eufP_{g,k^\minus,k^\plus}(\alpha^\minus_1,\ldots,\alpha^\minus_{k^\minus},
                            \alpha^\plus_1,\ldots,\alpha^\plus_{k^\plus}
                           )
$,
with $\alpha^\pm_i\in\eusC^0_H$,
such that the following four conditions are satisfied:

\item $\bullet$

For any $g\ge0$,
the coherent section of $\eufO_{g,0,0}$ is given by the canonical orientations.

\item $\bullet$

For each $\alpha_0\in\eusC_H$,
the coherent element of $\eusO_{\alphabar_0}(\alpha_0,\alpha_0)$ is
equal to the canonical orientation.

\item $\bullet$

The excision map \qcg{}
preserves or reverses the coherent orientations
according to the graded sign of the permutation
$$\eqalign{
  & (\alpha^\minus_{1,1},\ldots,\alpha^\minus_{1,k^\minus_1},
     \alpha^\plus_{1,1},\ldots,
     \alpha^\plus_{1,i-1},\mathop{\alpha_0}_\plus,\alpha^\plus_{1,i+1},\ldots,
     \alpha^\plus_{1,k^\plus_1+1},\cr
  & \qquad
     \alpha^\minus_{2,1},\ldots,\alpha^\minus_{2,j-1}
     ,\mathop{\alpha_0}_\minus,\alpha^\minus_{2,j+1},\ldots,
     \alpha^\minus_{2,k^\minus_2+1},
     \alpha^\plus_{2,1},\ldots,\alpha^\plus_{2,k^\plus_2})  \cr
  & \mapsto
    (\alpha^\minus_{2,1},\ldots,\alpha^\minus_{2,j-1},
     \alpha^\minus_{1,1},\ldots,\alpha^\minus_{1,k^\minus_1},
     \alpha^\minus_{2,j+1},\ldots,\alpha^\minus_{2,k^\minus_2+1},  \cr
  & \qquad\qquad
     \alpha^\plus_{1,1},\ldots,\alpha^\plus_{1,i-1},
     \alpha^\plus_{2,1},\ldots,\alpha^\plus_{2,k^\plus_2},
     \alpha^\plus_{1,i+1},\ldots,\alpha^\plus_{k^\plus_1+1},
     \mathop{\alpha_0}_\plus,\mathop{\alpha_0}_\minus), \cr
          }
  \eqno\qci
$$
where $\alpha$ has degree $\muH(\alpha)$.

\item $\bullet$

The excision map \qch{}
preserves or reverses the coherent orientations
according the graded sign of the permutation
$$\eqalign{
  & (\alpha^\minus_1,\ldots,
     \alpha^\minus_{j-1},\mathop{\alpha_0}_\minus,\alpha^\minus_{j+1},\ldots,
     \alpha^\minus_{k^\minus+1},
     \alpha^\plus_1,\ldots,
     \alpha^\plus_{i-1},\mathop{\alpha_0}_\plus,\alpha^\plus_{i+1},
     \alpha^\plus_{k^\plus+1}) \cr
  & \mapsto
    (\alpha^\minus_1,\ldots,
     \alpha^\minus_{j-1},\alpha^\minus_{j+1},\ldots,
     \alpha^\minus_{k^\minus+1}, \cr
  & \phantom{{}\mapsto(}
     \alpha^\plus_1,\ldots,
     \alpha^\plus_{i-1},\alpha^\plus_{i+1},
     \ldots,\alpha^\plus_{k^\plus+1},
     \mathop{\alpha_0}_\plus,\mathop{\alpha_0}_\minus ) , \cr
          }
  \eqno\qcj
$$
where $\alpha$ has degree $\muH(\alpha)$.

\noindent
We also refer to the induced sections of $\eufO_\sigma$
and $\eufO(\alpha^\minus,\alpha^\plus)$ as coherent orientations.
}

\noindent
Recall that graded sign means the sign of the permutation obtained by
removing all elements of even degree.

\proclaim Proposition 3.4.2.

For any regular Hamiltonian $H$,
there exist coherent systems of orientations.
Given one coherent system of orientations,
any other coherent system of orientations
can be obtained by choosing
a function $\tau:\eusC^0_H\to\{-1,1\}$
and changing the section of
$\eufO_{g,k^\minus,k^\plus}
  ((\alpha^\minus_1,\ldots,\alpha^\minus_{k^\minus}),
   (\alpha^\plus_1,\ldots,\alpha^\plus_{k^\plus})
  )$
by
$$\prod_{i=1}^{k^\minus} \tau(\alpha^\minus_i)
  \prod_{i=1}^{k^\plus} \tau(\alpha^\plus_i) .
  \eqno\qcx
$$

\demo Proof.

To keep the notation simple,
we write $u$ for $[u,c]\in\eufP_{g,k^\minus,k^\plus}$,
and $\eufO_u$ for $(\eufO_{g,k^\minus,k^\plus})_{[u,c]}$.
For each $\alpha\in\eusC^0_H$,
we choose elements $u^\plus_\alpha\in\eufP_{0,0,1}(\emptyset,(\alpha))$
and $u^\minus_\alpha\in\eufP_{0,1,0}((\alpha),\emptyset)$.
We then arbitrarily choose an element of each
$\eufO_{u^\plus_\alpha}$
that will serve as a coherent orientation.
For $\ell$ large enough,
$(u^0_\alpha,\alphabar) = \bfg^\ell_{11} (u^\plus_\alpha,u^\minus_\alpha)$
is defined.
(Recall that we write
``$\,\alphabar\,$''
for maps that are only approximately equal to $\alphabar$.)
The excision map
$$(\bfg^\ell_{11})_* :
  \eufO_{u^\plus_\alpha} \otimes \eufO_{u^\minus_\alpha}
  \to \eufO_{u^0_\alpha} \otimes \eusO_{\alphabar}
$$
has to preserve coherent orientations.
As we have chosen the coherent element of $\eufO_{u^\plus_\alpha}$,
and
$\eufO_{u^0_\alpha}$ and $\eusO_{\alphabar}$
have canonical elements,
this determines the coherent element of $\eufO_{u^\minus_\alpha}$.

Let
$$u\in\eufP_{g,k^\minus,k^\plus}
    ((\alpha^\minus_1,\ldots,\alpha^\minus_{k^\minus_1}),
   (\alpha^\plus_1,\ldots,\alpha^\plus_{k^\plus_1})
  )
.$$
By repeated use of the excision maps we get a map
$$\eqalign{
  & \eufO_{u^\plus_{\alpha^\minus_{k^\minus}}}
    \otimes\cdots\otimes
    \eufO_{u^\plus_{\alpha^\minus_1}}
    \otimes
    \eufO_u
    \otimes
    \eufO_{u^\minus_{\alpha^\plus_{k^\plus}}}
    \otimes\cdots\otimes
    \eufO_{u^\minus_{\alpha^\plus_1}} \cr
  & \to
    \eusO_{\alphabar^\minus_{k^\minus}}
    \otimes\cdots\otimes
    \eusO_{\alphabar^\minus_1}
    \otimes
    \eufO_{u^0}
    \otimes
    \eusO_{\alphabar^\plus_{k^\plus}}
    \otimes\cdots\otimes
    \eusO_{\alphabar^\plus_1} \cr
          }
$$
where $u^0\in\eufP_{g,0,0}(\emptyset,\emptyset)$.
This excision map has to preserve the coherent orientations.
This determines the coherent element of $\eufO_u$.

These orientations are unique,
once we have chosen elements of $\eufO_{u^\plus_\alpha}$
for all $\alpha\in\eusC^0_H$.
A different choice of elements of $\eufO_{u^\plus_\alpha}$
changes the coherent orientation by a factor \qcx.

We have chosen elements of
$\eufO_{u^\plus_\alpha}$, $\eufO_{u^\minus_\alpha}$
and $\eusO_{\alphabar}$ twice;
first in the initial step of the construction,
and then when we defined the coherent element
of an arbitrary $\eufO_u$.
Both elements of
$\eufO_{u^\plus_\alpha}$ and $\eufO_{u^\minus_\alpha}$
are determined by requiring the map $(\bfg^\ell_{11})_*$
to preserve the orientations.
Thus the two elements of
$\eufO_{u^\plus_\alpha}$ and $\eufO_{u^\minus_\alpha}$ are identical.

%%%%%%%%%%%%%%%%%%%%%%%%%%%%%%%%%%%%%%%%%%%%%%%%%%%%%%%%

\pageinsert

\vss

$$\comdia{
    \eqalign{
    & \eufO_{u_1}\otimes\eufO_{u_2} \cr
    & \otimes\eufO_{u^\minus_{\alpha^\plus_{2,k^\plus_2}}}
      \otimes\cdots\otimes
      \eufO_{u^\minus_{\alpha^\plus_{2,1}}} \cr
    & \otimes\eufOtilde_{u^\plus_{\alpha^\minus_{2,k^\minus_1}}}
      \otimes\cdots\otimes
      \eufOtilde_{u^\plus_{\alpha^\minus_{2,j+1}}} \cr
    & \otimes\eufOtilde_{u^\plus_{\alpha_0}}
      \otimes\eufOtilde_{u^\plus_{\alpha^\minus_{2,j-1}}}
      \otimes\cdots\otimes
      \eufOtilde_{u^\plus_{\alpha^\minus_{2,1}}} \cr
    & \otimes\eufO_{u^\minus_{\alpha^\plus_{1,k^\plus_1}}}
      \otimes\cdots\otimes
      \eufO_{u^\minus_{\alpha^\plus_{1,i+1}}} \cr
    & \otimes
      \eufO_{u^\minus_{\alpha_0}}
      \otimes
      \eufO_{u^\minus_{\alpha^\plus_{1,i-1}}}
      \otimes\cdots\otimes
      \eufO_{u^\minus_{\alpha^\plus_{1,1}}} \cr
    & \otimes\eufOtilde_{u^\plus_{\alpha^\minus_{1,k^\minus_1}}}
      \otimes\cdots\otimes
      \eufOtilde_{u^\plus_{\alpha^\minus_{1,1}}} \cr
            }
  & \quad\mapright{}\quad
  & \eqalign{
    & \eufO_{u_1}
      \otimes
      \eufO_{u_2} \cr
    & \otimes
      \eufO_{u^\minus_{\alpha^\plus_{1,k^\plus_1}}}
      \otimes\cdots\otimes
      \eufO_{u^\minus_{\alpha^\plus_{1,i+1}}} \cr
    & \otimes
      \eufO_{u^\minus_{\alpha^\plus_{2,k^\plus_2}}}
      \otimes\cdots\otimes
      \eufO_{u^\minus_{\alpha^\plus_{2,1}}} \cr
    & \otimes\eufO_{u^\minus_{\alpha^\plus_{1,i-1}}}
      \otimes\cdots\otimes
      \eufO_{u^\minus_{\alpha^\plus_{1,1}}} \cr
    & \otimes
      \eufOtilde_{u^\plus_{\alpha^\minus_{2,k^\minus_1}}}
      \otimes\cdots\otimes
      \eufOtilde_{u^\plus_{\alpha^\minus_{2,j+1}}} \cr
    & \otimes
      \eufOtilde_{u^\plus_{\alpha^\minus_{1,k^\minus_1}}}
      \otimes\cdots\otimes
      \eufOtilde_{u^\plus_{\alpha^\minus_{1,1}}} \cr
    & \otimes
      \eufOtilde_{u^\plus_{\alpha^\minus_{2,j-1}}}
      \otimes\cdots\otimes
      \eufOtilde_{u^\plus_{\alpha^\minus_{2,1}}} \cr
    & \otimes
      \eufO_{u^\minus_{\alpha_0}}
      \otimes
      \eufOtilde_{u^\plus_{\alpha_0}} \cr
            } \cr
    \mapdown{}
  &
  & \mapdown{} \cr
    \eqalign{
    & \eufO_{u_1}\otimes\eufO_{u^0_2} \cr
    & \otimes
      \eusO_{\alphabar^\plus_{2,k^\plus_2}}
      \otimes\cdots\otimes
      \eusO_{\alphabar^\plus_{2,1}} \cr
    & \otimes
      \eusO_{\alphabar^\minus_{2,k^\minus_1}}
      \otimes\cdots\otimes
      \eusO_{\alphabar^\minus_{2,j+1}} \cr
    & \otimes\eusO_{\alphabar_0}
      \otimes\eusO_{\alphabar^\minus_{2,j-1}}
      \otimes\cdots\otimes
      \eusO_{\alphabar^\minus_{2,1}} \cr
    & \otimes\eufO_{u^\minus_{\alpha^\plus_{1,k^\plus_1}}}
      \otimes\cdots\otimes
      \eufO_{u^\minus_{\alpha^\plus_{1,i+1}}} \cr
    & \otimes
      \eufO_{u^\minus_{\alpha_0}}
      \otimes
      \eufO_{u^\minus_{\alpha^\plus_{1,i-1}}}
      \otimes\cdots\otimes
      \eufO_{u^\minus_{\alpha^\plus_{1,1}}} \cr
    & \otimes\eufOtilde_{u^\plus_{\alpha^\minus_{1,k^\minus_1}}}
      \otimes\cdots\otimes
      \eufOtilde_{u^\plus_{\alpha^\minus_{1,1}}} \cr
            }
  &
  & \eqalign{
    & \eufO_{u_3}
      \otimes
      \eusO_{\alphabar_0} \cr
    & \otimes
      \eufO_{u^\minus_{\alpha^\plus_{1,k^\plus_1}}}
      \otimes\cdots\otimes
      \eufO_{u^\minus_{\alpha^\plus_{1,i+1}}} \cr
    & \otimes
      \eufO_{u^\minus_{\alpha^\plus_{2,k^\plus_2}}}
      \otimes\cdots\otimes
      \eufO_{u^\minus_{\alpha^\plus_{2,1}}} \cr
    & \otimes\eufO_{u^\minus_{\alpha^\plus_{1,i-1}}}
      \otimes\cdots\otimes
      \eufO_{u^\minus_{\alpha^\plus_{1,1}}} \cr
    & \otimes
      \eufOtilde_{u^\plus_{\alpha^\minus_{2,k^\minus_1}}}
      \otimes\cdots\otimes
      \eufOtilde_{u^\plus_{\alpha^\minus_{2,j+1}}} \cr
    & \otimes
      \eufOtilde_{u^\plus_{\alpha^\minus_{1,k^\minus_1}}}
      \otimes\cdots\otimes
      \eufOtilde_{u^\plus_{\alpha^\minus_{1,1}}} \cr
    & \otimes
      \eufOtilde_{u^\plus_{\alpha^\minus_{2,j-1}}}
      \otimes\cdots\otimes
      \eufOtilde_{u^\plus_{\alpha^\minus_{2,1}}} \cr
    & \otimes
      \eufO_{u^\minus_{\alpha_0}}
      \otimes
      \eufOtilde_{u^\plus_{\alpha_0}} \cr
            } \cr
    \mapdown{}
  &
  & \mapdown{} \cr
    \eqalign{
    & \eufO_{u^0_2}
      \otimes
      \eusO_{\alphabar^\plus_{2,k^\plus_2}}
      \otimes\cdots\otimes
      \eusO_{\alphabar^\plus_{2,1}} \cr
    & \otimes
      \eusO_{\alphabar^\minus_{2,k^\minus_1}}
      \otimes\cdots\otimes
      \eusO_{\alphabar^\minus_{2,j+1}} \cr
    & \otimes\eusO_{\alphabar_0}
      \otimes\eusO_{\alphabar^\minus_{2,j-1}}
      \otimes\cdots\otimes
      \eusO_{\alphabar^\minus_{2,1}} \cr
    & \otimes\eufO_{u^0_1}
      \otimes\eusO_{\alpha^\plus_{1,k^\plus_1}}
      \otimes\cdots\otimes
      \eusO_{\alpha^\plus_{1,i+1}} \cr
    & \otimes
      \eusO_{\alpha_0}
      \otimes
      \eusO_{\alpha^\plus_{1,i-1}}
      \otimes\cdots\otimes
      \eusO_{\alpha^\plus_{1,1}} \cr
    & \otimes\eusO_{\alpha^\minus_{1,k^\minus_1}}
      \otimes\cdots\otimes
      \eusO_{\alpha^\minus_{1,1}} \cr
            }
  & \mapright{}
  & \eqalign{
    & \eufO_{u^0_3}
      \otimes
      \eusO_{\alphabar_0} \cr
    & \otimes
      \eusO_{\alphabar^\plus_{1,k^\plus_1}}
      \otimes\cdots\otimes
      \eusO_{\alphabar^\plus_{1,i+1}} \cr
    & \otimes
      \eusO_{\alphabar^\plus_{2,k^\plus_2}}
      \otimes\cdots\otimes
      \eusO_{\alphabar^\plus_{2,1}} \cr
    & \otimes
      \eusO_{\alphabar^\plus_{1,i-1}}
      \otimes\cdots\otimes
      \eusO_{\alphabar^\plus_{1,1}} \cr
    & \otimes
      \eusO_{\alphabar^\minus_{2,k^\minus_1}}
      \otimes\cdots\otimes
      \eusO_{\alphabar^\minus_{2,j+1}} \cr
    & \otimes
      \eusO_{\alphabar^\minus_{1,k^\minus_1}}
      \otimes\cdots\otimes
      \eusO_{\alphabar^\minus_{1,1}} \cr
    & \otimes
      \eusO_{\alphabar^\minus_{2,j-1}}
      \otimes\cdots\otimes
      \eusO_{\alphabar^\minus_{2,1}} \cr
    & \otimes
      \eufO_{u^0_{\alpha_0}}
      \otimes
      \eusO_{\alphabar_0} \cr
            } \cr
       }
$$

\vss

\endinsert

%%%%%%%%%%%%%%%%%%%%%%%%%%%%%%%%%%%%%%%%%%%%%%%%%%%%%%%

We need to verify that these orientations satisfy the four
conditions of Def.~3.4.1.
The first condition is satisfied trivially.
The third condition is verified as follows.
We write $\eufOtilde_{u^\plus_\alpha}$
for $\eufO_{u^\plus_\alpha}$,
but with a designated element that is $(-1)^{\mu(\alpha)}$
times the coherent element.
It follows from Prop.~3.3.3,
or rather the simplified formula
$$d \equiv \sum_{i=1}^{k^\minus} \mu(\alpha^\minus_i)
      + \sum_{i=1}^{k^\plus} \mu(\alpha^\plus_i)
  \quad \bmod 2,$$
and the discussion of permutations earlier in this section
that the excision maps
$$\eufO_{u^\minus_\alpha} \otimes \eufOtilde_{u^\plus_\alpha}
  \to \eufO_{u^0_\alpha} \otimes \eusO_{\alphabar}
$$
and
$$\eqalign{
  & \eufO_u
    \otimes
    \eufO_{u^\minus_{\alpha^\plus_{k^\plus}}}
    \otimes\cdots\otimes
    \eufO_{u^\minus_{\alpha^\plus_1}}
    \otimes
    \eufOtilde_{u^\plus_{\alpha^\minus_{k^\minus}}}
    \otimes\cdots\otimes
    \eufOtilde_{u^\plus_{\alpha^\minus_1}} \cr
  & \to
    \eufO_{u^0}
    \otimes
    \eusO_{\alphabar^\minus_{k^\plus}}
    \otimes\cdots\otimes
    \eusO_{\alphabar^\plus_1}
    \otimes
    \eusO_{\alphabar^\minus_{k^\minus}}
    \otimes\cdots\otimes
    \eusO_{\alphabar^\minus_1}  \cr
          }
$$
preserve the orientations induced by the coherent orientations
and the designated elements of the spaces
$\eufOtilde_{u^\plus_{\alpha^\minus_i}}$.
That the third condition is satisfied now
follows from the commutative diagram on the next page.
The two maps on the left and the map on the lower right
preserve the orientations.
The bottom map is defined
to make the diagram commute.
It is induced by
a fairly complicated composition of  excision operators
on the bundles $T^\fib\eufP$ and $\IE$.
However,
a moment's thought shows that if we choose the parameter $\ell$
for the gluing of the $u^\pm_\alpha$'s at least three times
larger than the parameter $\ell$ for the gluing of $u_1$ and $u_2$,
then this operator splits as a direct sum of
operators
$T^\fib_{u^0_1}\eufP_{g_1,0,0}
  \oplus
  T^\fib_{u^0_2}\eufP_{g_2,0,0}
  \to
  T^\fib_{u^0_3}\eufP_{g_1+g_2,0,0}
  \oplus
  T^\fib_{u^0_{\alpha_0}}\eufP_{0,0,0}
$
and
$(\IE_{g_1,0,0})_{u^0_1}
  \oplus
  (\IE_{g_2,0,0})_{u^0_2}
  \to
  (\IE_{g_1+g_2,0,0})_{u^0_3}
  \oplus
  (\IE_{0,0,0})_{u^0_{\alpha_0}},
$
and the identity operators on
$T_{\alphabar^\minus_{2,1}}\eusP
 \cdots,
 T_{\alphabar^\plus_{1,k^\plus_1}}\eusP,
 T_{\alphabar_0}\eusP
$
and
$\rmE_{\alphabar^\minus_{2,1}},
 \cdots,
 \rmE_{\alphabar^\plus_{1,k^\plus_1}},
 \rmE_{\alphabar_0} .
$
The operators
$T^\fib_{u^0_1}\eufP_{g_1,0,0}
  \oplus
  T^\fib_{u^0_2}\eufP_{g_2,0,0}
  \to
  T^\fib_{u^0_3}\eufP_{g_1+g_2,0,0}
  \oplus
  T^\fib_{u^0_{\alpha_0}}\eufP_{0,0,0}
$
and
$(\IE_{g_1,0,0})_{u^0_1}
  \oplus
  (\IE_{g_2,0,0})_{u^0_2}
  \to
  (\IE_{g_1+g_2,0,0})_{u^0_3}
  \oplus
  (\IE_{0,0,0})_{u^0_{\alpha_0}}
$
commute with $\bfJ$,
up to lower order terms.
In this case the coherent orientations are induced by the complex structures.
It follows that the induced map
$\eufO_{u^0_1}\otimes \eufO_{u^0_2}\to \eufO_{u^0_3}
\otimes \eufO_{u^0_{\alpha_0}}$
preserves the coherent orientations.
Hence the bottom map is orientation preserving.

The top map is a permutation,
and changes the orientation according to the
graded sign of the permutation \qci.
As the diagram commutes,
it follows that the upper right map also
changes the orientation according to the graded sign of the permutation \qci.
But this map is given by the excision map
$\eufO_{u_1}\otimes\eufO_{u_2}\to\eufO_{u_3}\otimes\eusO_{\alphabar_0}$.
It follows that the third condition is satisfied.
The fourth condition is verified the same way.

Finally we verify the second condition.
By definition,
the excision map $\bfg^\ell_{11}$ above preserves
the orientations given by the canonical element of $\eufO_{\alphabar}$.
By the third condition $\bfg^\ell_{11}$ preserves the orientations  given by
the coherent element of $\eufO_{\alphabar}$.
Hence the two elements of  $\eusO_{\alphabar}$ are identical.
\enddemo

Given $\sigma_1 : [0,1]^{q_1}\to \eufJ^*_{g_1,k^\minus_1,k^\plus_1+1}(A)$
and $\sigma_2 : [0,1]^{q_2}\to \eufJ^*_{g_2,k^\minus_2+1,k^\plus_2}(A)$
we get two families
$D\Psi_{\sigma_1}
  \boxplus
  D\Psi_{\sigma_2}
$
and
$D\Psi_{\sigma_1\gop^\ell_{ij}\sigma_2}
 \oplus
 D\Psi_{\sigma_{\alphabar_0}}
$
of Fredholm operators parametrized by $[0,1]^{q_1+q_2}$.
These families are intertwined by the excision operators
$(\bfg^\ell_{ij})_*$,
up to compact operators.
Thus there is an induced excision map
$$(\bfg^\ell_{ij})_* :
  \eufO_{\sigma_1} \boxtimes \eufO_{\sigma_2}
  \to
  \eufO_{\sigma_1\gop^\ell_{ij}\sigma_2} \otimes \eusO_{\alphabar_0}.
  \eqno\qck
$$
Similarly there is an excision map
$$(\bfc^\ell_{ij})_* :
  \eufO_\sigma
  \to
  \eufO_{\cop^\ell_{ij}\sigma} \otimes \eusO_{\alphabar_0}.
  \eqno\qcl
$$

\proclaim Proposition 3.4.3.

{
Any coherent system of orientations of the moduli spaces
$\eufM^d_\sigma$ has the following three properties.

\item $\bullet$

The moduli space $\eusM^0(\alpha_0,\alpha_0)$
consists of a single point $\alphabar_0$,
and this point is positively oriented.

\item $\bullet$

The map \qck{}
preserves or reverses the orientation according to the graded sign
of the permutation
$$\eqalign{
  & (\sigma_1,
     \alpha^\minus_{1,1},\ldots,\alpha^\minus_{1,k^\minus_1},
     \alpha^\plus_{1,1},\ldots,
     \alpha^\plus_{1,i-1},\mathop{\alpha_0}_\plus,\alpha^\plus_{1,i+1},\ldots,
     \alpha^\plus_{1,k^\plus_1+1},\cr
  & \qquad
     \sigma_2,
     \alpha^\minus_{2,1},\ldots,
     \alpha^\minus_{2,j-1},\mathop{\alpha_0}_\minus,\alpha^\minus_{2,j+1},
     \ldots,
     \alpha^\minus_{2,k^\minus_2+1},
     \alpha^\plus_{2,1},\ldots,\alpha^\plus_{2,k^\plus_2})  \cr
  & \mapsto
    (\sigma_1,\sigma_2,\alpha^\minus_{2,1},\ldots,\alpha^\minus_{2,j-1},
     \alpha^\minus_{1,1},\ldots,\alpha^\minus_{1,k^\minus_1},
     \alpha^\minus_{2,j+1},\ldots,\alpha^\minus_{2,k^\minus_2+1},  \cr
  & \qquad\qquad
     \alpha^\plus_{1,1},\ldots,\alpha^\plus_{1,i-1},
     \alpha^\plus_{2,1},\ldots,\alpha^\plus_{2,k^\plus_2},
     \alpha^\plus_{1,i+1},\ldots,\alpha^\plus_{k^\plus_1+1},
     \mathop{\alpha_0}_\plus,\mathop{\alpha_0}_\minus). \cr
          }
  \eqno\qcm
$$
where $\alpha$ has degree $\muH(\alpha)$,
$\sigma_1$ has degree $q_1$,
and $\sigma_2$ has degree $q_2$.

\item $\bullet$

The map \qcl{}
preserves or reverses the orientations according to the graded sign
of the permutation \qcj.

}

To show this one argues the same way as when we verified
the third condition of Def.~3.4.1 in the proof of Prop.~3.4.2,
keeping track of the additional factors $\eusO_{[0,1]^{q_1}}$,
$\eusO_{[0,1]^{q_2}}$, and $\eusO_{[0,1]^q}$.

Let $S_k$ denote the symmetric group on $k$ letters.
Let $\rho=(\rho^\minus,\rho^\plus)$
with
$\rho^\minus : \{1,\ldots,k^\minus\} \to \{1,\ldots,k^\minus\}$
and $\rho^\plus : \{1,\ldots,k^\plus\} \to \{1,\ldots,k^\plus\}$
be an element of $S_{k^\minus}\times S_{k^\plus}$.
Let $\phi$ be any diffeomorphism
with $\phi(p^\plusminus_i) = p^\plusminus_{\rho^\plusminus(i)}$
and $\phi_* X^\pm_i=X^\pm_{\rho^\plusminus(i)}$.
Then for any $c\in\eusJhat^*_{g,k^\minus,k^\plus}(A)$,
$[\phi_* c] \in\eufJ^*_{g,k^\minus,k^\plus}(A)$
is uniquely determined by $[c]\in\eufJ^*_{g,k^\minus,k^\plus}(A)$
and $(\rho^\minus,\rho^\plus)\in S_{k^\minus,k^\plus}$.
This defines a natural group action
$$(S_{k^\minus}\times S_{k^\plus})
  \times \eufJ^*_{g,k^\minus,k^\plus}(A)
  \to \eufJ^*_{g,k^\minus,k^\plus}(A) .
$$
For each $\rho\in S_{k^\minus}\times S_{k^\plus}$
a similar construction gives a map
$$\eqalign{
  & \eufP_\sigma
    ((\alpha^\minus_1,\ldots,\alpha^\minus_{k^\minus}),
      (\alpha^\plus_1,\ldots,\alpha^\plus_{k^\plus})
     ) \cr
 & \quad\to
   \eufP_{\rho.\sigma}
    ((\alpha^\minus_{\rho^\minus(1)},\ldots,
      \alpha^\minus_{\rho^\minus(k^\minus)}),
     (\alpha^\plus_{\rho^\plus(1)},\ldots,
      \alpha^\plus_{\rho^\plus(k^\plus)})
    ). \cr
          }
$$
This map has natural lifts to bundle maps
$T^\fib\eufP_\sigma\to T^\fib\eufP_{\rho.\sigma}$
and $\IF_\sigma\to\IF_{\rho.\sigma}$.
These bundle maps intertwine $D\Psi_\sigma$
and $D\Psi_{\rho.\sigma}$.
Hence there is an induced map
$$\eufO_\sigma \to \eufO_{\rho.\sigma}. \eqno\qco$$

\proclaim Proposition 3.4.4.

The map \qco{} preserves or reverses the coherent orientation
according to the graded sign of the permutation
$$\eqalign{
  & (\alpha^\minus_1,\ldots,\alpha^\minus_{k^\minus},
     \alpha^\plus_1,\ldots,\alpha^\plus_{k^\plus}
    ) \cr
  & \quad\to
    (\alpha^\minus_{\rho^\minus(1)},\ldots,
     \alpha^\minus_{\rho^\minus(k^\minus)},
     \alpha^\plus_{\rho^\plus(1)},\ldots,
     \alpha^\plus_{\rho^\plus(k^\plus)}
    ). \cr
          }
  \eqno\qnp
$$

This is shown the same way as Prop.~3.4.3.

\subheading 3.5. Uhlenbeck-Gromov compactness.

The usual integration by parts argument,
shows that if
$$u\in\eufP^d((\alpha^\minus_1,\ldots,\alpha^\minus_{k^\minus}),
              (\alpha^\minus_1,\ldots,\alpha^\minus_{k^\minus})
             )
$$
then
$$E[u]\le C +\int_{\Sigma_{g,k^\minus,k^\plus}}u^*\omega,$$
where $C$ only depends on $H$ and $R$.
The integral of $u^*\omega$ is uniquely determined by
the periodic orbits $\alpha^\pm_i$ and $d$.
In fact, if $u_1\in\eufP^{d_1}$
and $u_2\in\eufP^{d_2}$,
then $\int_\Sigma (u_1^*\omega-u_2^*\omega)$
equals $k$ times the relative first Chern class of $u_1^*TM$
and $u_2^*TM$,
which,
by excision and the Atiyah-Singer theorem,
equals $k(d_1-d_2)/2$.
With a uniform estimate of $E[u]$ at hand,
the following compactness result follows by the standard
Uhlenbeck-Gromov argument;
see for instance [PW].

\proclaim Proposition 3.5.1.

If $\sigma:[0,1]^q\to\eufJ^*_{g,k^\minus,k^\plus}(A)$ is semi-regular,
then the moduli space
$\eufM^d_\sigma((\alpha^\minus_1,\ldots,\alpha^\minus_{k^\minus}),
                (\alpha^\plus_1,\ldots,\alpha^\plus_{k^\plus}))
$
is empty for all $d<0$,
and any sequence in
$\eufM^0_\sigma((\alpha^\minus_1,\ldots,\alpha^\minus_{k^\minus}),
                (\alpha^\plus_1,\ldots,\alpha^\plus_{k^\plus}))$
has a strongly convergent subsequence.

There is a similar compactness result for 1-dimensional moduli spaces.

\proclaim Proposition 3.5.2.

{
If $\sigma:[0,1]^q\to\eufJ^*_{g,k^\minus,k^\plus}(A)$ is regular,
then any sequence $[u_n,c_n]$ in
$\eufM^1_\sigma
   ((\alpha^\minus_1,\ldots,\alpha^\minus_{k^\minus}),
    (\alpha^\plus_1,\ldots,\alpha^\plus_{k^\plus}))
$
has a subsequence,
that we also denote $[u_n,c_n]$,
such that one of the following conditions holds:

\item $\bullet$

The subsequence converges strongly to an element of
$$\qquad\qquad\eufM^1_\sigma
   ((\alpha^\minus_1,\ldots,\alpha^\minus_{k^\minus}),
    (\alpha^\plus_1,\ldots,\alpha^\plus_{k^\plus}))
.$$

\item $\bullet$

There exists $i\in\{1,\ldots,k^\minus\}$, $\alpha_0\in\eusC_H$,
and a sequence $\ell_n$, such that
$\ell_n\to\infty$ as $n\to\infty$,
$(\alphabar_0,[u_n,c_n])$ is in the range of $\bfg^{\ell_n}_{1i}$,
and
$(\bfg^{\ell_n}_{1i})^{-1} (\alphabar_0,[u_n,c_n])$
converges strongly to a pair in
$$\qquad\qquad\eusM^1(\alpha^\minus_i,\alpha_0)
  \times\eufM^0_\sigma
    ((\alpha^\minus_1
      ,\ldots,\alpha^\minus_{i-1},\alpha_0,\alpha^\minus_{i+1}
      ,\ldots,\alpha^\minus_{k^\minus}
     )
     ,(\alpha^\plus_1,\ldots,\alpha^\plus_{k^\plus})
    )
.$$
If $\alpha^\minus_i\in\eusC^0_H$,
then $\alpha_0\in\eusC^0_H$.

\item $\bullet$

There exists $i\in\{1,\ldots,k^\plus\}$, $\alpha_0\in\eusC_H$,
and a sequence $\ell_n$, such that
$\ell_n\to\infty$ as $n\to\infty$,
and
$([u_n,c_n],\alphabar_0)$ is in the range of $(\bfg^{\ell_n}_{i1})^{-1}$
and
$(\bfg^{\ell_n}_{i1})^{-1}([u_n,c_n],\alphabar_0)$
converges strongly to a pair in
$$\qquad\qquad\eufM^0_\sigma
     ((\alpha^\minus_1,\ldots,
       \alpha^\minus_{k^\minus}
      ),
      (\alpha^\plus_1,\ldots,
       \alpha^\plus_{i-1},\alpha_0,\alpha^\plus_{i+1},\ldots,
       \alpha^\plus_{k^\plus})
     )
\times \eusM^1(\alpha_0,\alpha^\plus_i).
$$
If $\alpha^\plus_i\in\eusC^0_H$,
then $\alpha_0\in\eusC^0_H$.

}

There are similar estimates for $E[u]$ with
$u\in\eufM^d_{\sigma_1\gop^\ell_{ij}\sigma_2}$
that are uniform in $\ell$,
if the energy integrand \qbz{} is redefined as
$|\nabla u -\eta(t_0)\,\eta(\ell-t_0)\,d\theta_0\otimes\nabla H|^2$
on the neck.
Then the standard Uhlenbeck-Gromov argument
gives the following compactness result.

\proclaim Proposition 3.5.3.

If $\sigma_1:[0,1]^{q_1}\to \eufJ^*_{g_1,k^\minus_1,k^\plus_1+1}(A)$
and $\sigma_2:[0,1]^{q_2}\to \eufJ^*_{g_2,k^\minus_2+1,k^\plus_2}(A)$
are semi-regular,
and $\sigma_1\gop^\ell_{ij}\sigma_2$ takes values in
$\eufJ^*_{g_1+g_2,k^\minus_1+k^\minus_2,k^\plus_1+k^\plus_2}(A)$
for all $\ell\ge0$,
then,
for $\ell$ large enough,
$$\eqalign{
   \eufM^d_{\sigma_1\gop^\ell_{ij}\sigma_2}
      (
   &   (\alpha^\minus_{2,1},\ldots,\alpha^\minus_{2,j-1},
        \alpha^\minus_{1,1},\ldots,\alpha^\minus_{1,k^\minus_1},
        \alpha^\minus_{2,j+1},\ldots,\alpha^\minus_{2,k^\minus_2+1}
       ), \cr
   &   (\alpha^\plus_{1,1},\ldots,\alpha^\plus_{1,i-1},
        \alpha^\plus_{2,1},\ldots,\alpha^\plus_{2,k^\plus_2},
        \alpha^\plus_{1,i+1},\ldots,\alpha^\plus_{1,k^\plus_1+1}
       )
      ) \cr
          }
$$
is empty for all $d<0$,
and for any sequence $\ell_n$ of positive numbers
such that $\ell_n\to\infty$ as $n\to\infty$
and for any sequence
$$\eqalign{
   [u_n,c_n]\in\eufM^0_{\sigma_1\gop^{\ell_n}_{ij}\sigma_2}
      (
   &   (\alpha^\minus_{2,1},\ldots,\alpha^\minus_{2,j-1},
        \alpha^\minus_{1,1},\ldots,\alpha^\minus_{1,k^\minus_1},
        \alpha^\minus_{2,j+1},\ldots,\alpha^\minus_{2,k^\minus_2+1}
       ), \cr
   &   (\alpha^\plus_{1,1},\ldots,\alpha^\plus_{1,i-1},
        \alpha^\plus_{2,1},\ldots,\alpha^\plus_{2,k^\plus_2},
        \alpha^\plus_{1,i+1},\ldots,\alpha^\plus_{1,k^\plus_1+1}
       )
      ) \cr
          }
$$
there exists $\alpha_0\in\eusC_H$ and a subsequence,
that we also denote $\ell_n$ and $[u_n,c_n]$,
such that $([u_n,c_n],\alphabar_0)$
lies in the range of $\bfg^{\ell_n}_{ij}$
and $(\bfg^{\ell_n}_{ij})^{-1}([u_n,c_n],\alphabar_0)$ converges strongly
to a pair in
$$\eqalign{
& \eufM^0_{\sigma_1}
  ((\alpha^\minus_{1,1},\ldots,\alpha^\minus_{1,k^\minus_1}),
   (\alpha^\plus_{1,1},\ldots,
    \alpha^\plus_{1,i-1},\alpha_0,\alpha^\plus_{1,i+1},\ldots
    \alpha^\plus_{1,k^\plus_1+1}
   )
  ) \cr
& \times\eufM^0_{\sigma_2}
  ((\alpha^\minus_{2,1},\ldots,
    \alpha^\minus_{2,j-1},\alpha_0,\alpha^\minus_{2,j+1},\ldots
    \alpha^\minus_{2,k^\minus_2+1}),
   (\alpha^\plus_{2,1},\ldots,
    \alpha^\plus_{2,k^\plus_2}
   )
  ) . \cr
        }
$$
If $g_1=0$ and $\alpha^\pm_{1,\nu}\in\eusC^0_H$ for
$\nu=1,\ldots,k^\plusminus_1$,
or if $g_2=0$ and $\alpha^\pm_{2,\nu}\in\eusC^0_H$
for $\nu=1,\ldots,k^\plusminus_2$,
then $\alpha_0\in\eusC^0_H$.

There is a similar compactness theorem for $\eufM^d_{\cop^\ell_{ij}\sigma}$
with $d\le0$.

\subheading 3.6. Gluing of moduli spaces.

We now change our notation slightly,
and let $A$ denote a triple $(H,\bfJ,\euf0)$
where $(H,\bfJ)$ is a regular pair of a time-dependent Hamiltonian
and almost complex structure on $M$
and $\eufo$ is a coherent system of orientations.
It follows from Prop.~3.5.1
that if $H$ is regular and
$\sigma:[0,1]^q\to\eufJ^*_{g,k^\minus,k^\plus}(A)$
is transverse to
$$\pi : \eufM^{-q}
        ((\alpha^\minus_1,\ldots,\alpha^\minus_{k^\minus_1}),
         (\alpha^\plus_1,\ldots,\alpha^\plus_{k^\plus_1})
        )
    \to \eufJ^*_{g,k^\minus,k^\plus}(A)
,$$
then the moduli space
$\eufM^0_\sigma
  ((\alpha^\minus_1,\ldots,\alpha^\minus_{k^\minus_1}),
   (\alpha^\plus_1,\ldots,\alpha^\plus_{k^\plus_1})
  )
$
consists of a finite number of points.
If all $\alpha^\pm_i\in\eusC^0_H$,
then the moduli space becomes a finite set of points,
with signs
given by $\eufo$.
We denote the number of points,
counted with signs,
by
$$\#\eufM^0_\sigma
  ((\alpha^\minus_1,\ldots,\alpha^\minus_{k^\minus_1}),
   (\alpha^\plus_1,\ldots,\alpha^\plus_{k^\plus_1})
  )
.$$

Similarly,
$\eusM^1(\alpha^\minus,\alpha^\plus)$ is a finite union of affine
lines.
The affine structures give the lines
canonical orientations.
The coherent system $\eufo$ also induces orientations.
We count a line as positive if the coherent and
the canonical orientations agree,
and negative otherwise.
We denote the number of lines, counted with signs,
by
$$\#\eusM^1(\alpha^\minus,\alpha^\plus).$$

\proclaim Proposition 3.6.1.

For any contractible periodic orbit $\alpha$,
$$\#\eusM^0(\alpha,\alpha)=1.$$

This is an immediate consequence of Prop.~3.4.3.

\proclaim Proposition 3.6.2.

If $\rho=(\rho^\minus,\rho^\plus)\in S_{k^-}\times S_{k^+}$,
then
$$\eqalign{
  & \#\eufM^0_{\rho.\sigma}
    ((\alpha^\minus_{\rho^\minus(1)},\ldots,
      \alpha^\minus_{\rho^\minus(k^\minus)}
     ) ,
     (\alpha^\plus_{\rho^\plus(1)},\ldots,
      \alpha^\plus_{\rho^\plus(k^\plus)}
     )
    ) \cr
  & = \pm \#\eufM^0_{\rho.\sigma}
          ((\alpha^\minus_1,\ldots,\alpha^\minus_{k^\minus}),
           (\alpha^\plus_1,\ldots,\alpha^\plus_{k^\plus})
          ) \cr
          }
$$
where the sign is the graded sign of the permutation \qnp.

This is an immediate consequence of Prop.~3.4.4.

\proclaim Proposition 3.6.3.

{
If
$\sigma_1 : [0,1]^{q_1} \to \eufJ^*_{g_1,k^\minus_1,k^\plus_1+1}(A)$
and
$\sigma_2 : [0,1]^{q_2} \to \eufJ^*_{g_2,k^\minus_2+1,k^\plus_2}(A)$
are semi-regular
and
$$\sigma_1\gop^\ell_{ij}\sigma_2
  : [0,1]^{q_1+q_2}
    \to \eufJ_{g_1+g_2,k^\minus_1+k^\minus_2,k^\plus_1+k^\plus_2}(A)
$$
takes values in
$\eufJ^*_{g_1+g_2,k^\minus_1+k^\minus_2,k^\plus_1+k^\plus_2}(A)$
for all $\ell\ge0$,
then,
for sufficiently large $\ell$,
$\sigma_1\gop^\ell_{ij}\sigma_2$
is semi-regular,
and if all periodic orbits of $H$ are contractible,
then
$$\openup1\jot\vbox{
    \+    $\ds \#\eufM^0_{\sigma_1\gop^\ell_{ij}\sigma_2}
               (
          $%
    &     $\ds  (\alpha^\minus_{2,1},\ldots,\alpha^\minus_{2,j-1},
                 \alpha^\minus_{1,1},\ldots,\alpha^\minus_{1,k^\minus_1},
                 \alpha^\minus_{2,j+1},\ldots,\alpha^\minus_{2,k^\minus_2+1}
                 ) ,
          $\cr
    \+&   $\ds   (\alpha^\plus_{1,1},\ldots,\alpha^\plus_{1,i-1},
                  \alpha^\plus_{2,1},\ldots,\alpha^\plus_{2,k^\plus_2}
                  \alpha^\plus_{1,i+1},\ldots,\alpha^\plus_{1,k^\plus_1+1},
                 )
                )
          $\cr \cleartabs
    \smallskip
    \+    $\ds  =  (-1)^{q_1q_2}
                \smash{\sum_{\alpha_0\in\eusC^0_H}}
                \pm\,
          $
    &     $\ds  \#\eufM^0_{\sigma_1}
                (
          $%
    &     $\ds   (\alpha^\minus_{1,1},\ldots,\alpha^\minus_{1,k^\minus_1}),
          $\cr
    \+&&  $\ds   (\alpha^\plus_{1,1},\ldots,
                  \alpha^\plus_{1,i-1},\alpha_0,\alpha^\plus_{1,i+1},\ldots,
                  \alpha^\plus_{1,k^\plus_1+1}
                 )
                )
          $\cr
    \+&   $\ds \#\eufM^0_{\sigma_2}
               ((\alpha^\minus_{2,1},\ldots,
                 \alpha^\minus_{2,j-1},\alpha_0,\alpha^\minus_{2,j+1},\ldots,
                 \alpha^\minus_{2,k^\minus_2+1}
                ),
          $\cr
    \+&&  $\ds  (\alpha^\plus_{2,1},\ldots,\alpha^\plus_{2,k^\plus_2})
               )
          $ ,\cr
    }
$$
where the signs equal the graded signs of the permutations \qci.

The assumption that all periodic orbits of $H$ are contractible
can be replaced by $g_1=0$ and $\alpha^\pm_{1,\nu}\in\eusC^0_H$
for $\nu=1,\ldots,k^\plusminus_1$,
or by $g_2=0$ and $\alpha^\pm_{2,\nu}\in\eusC^0_H$
for $\nu=1,\ldots,k^\plusminus_2$,
}

\proclaim Proposition 3.6.4.

If
$\sigma:[0,1]^q\to\eufJ^*_{g,k^\minus+1,k^\plus+1}(A)$
is semi-regular
and
$$\cop^\ell_{ij}\sigma:[0,1]^q\to\eufJ_{g+1,k^\minus,k^\plus}(A)$$
takes values in $\eufJ^*_{g+1,k^\minus,k^\plus}(A)$
for all $\ell\ge0$, then,
for sufficiently large $\ell$,
$\cop^\ell_{ji}\sigma$
is semi-regular,
and if all periodic orbits of $H$ are contractible,
then
$$\openup1\jot\vbox{
  \+   $\ds \#\eufM^0_{\cop^\ell_{ij}\sigma}
            ((\alpha^\minus_1,\ldots,
              \alpha^\minus_{j-1},\alpha^\minus_{j+1},\ldots,
              \alpha^\minus_{k^\minus+1}
             ),
             (\alpha^\plus_1,\ldots,
              \alpha^\plus_{i-1},\alpha^\plus_{i+1},\ldots,
              \alpha^\plus_{k^\plus+1}
             )
            )
       $\cr
  \smallskip
  \+   $\ds = \smash{\sum_{\hidewidth\alpha_0\in\eusC^0_H\hidewidth}}
            \pm\,\#\eufM^0_\sigma
            (
       $%
  &    $\ds  (\alpha^\minus_1,\ldots,
              \alpha^\minus_{j-1},\alpha_0,\alpha^\minus_{j+1},\ldots,
              \alpha^\minus_{k^\minus+1}
             ),
       $\cr
  \+&  $\ds  (\alpha^\plus_1,\ldots,
              \alpha^\plus_{i-1},\alpha_0,\alpha^\plus_{i+1},\ldots,
              \alpha^\plus_{k^\plus+1}
             )
            )
       $ , \cr
      }
$$
where the signs equal the graded signs of the permutation \qcj.

Given a smooth map
$$\sigma : [0,1]^q \to \eufJ^*_{g,k^\minus,k^\plus}(A),$$
we can form the boundary maps
$$\partial^0_i\sigma,\partial^1_i\sigma
  : [0,1]^{q-1} \to \eufJ^*_{g,k^\minus,k^\plus}(A)$$
given by
$$\partial^0_\nu\sigma(x_1,\ldots,x_{q-1})
  = \sigma(x_1,\ldots,x_{\nu-1},0,x_\nu,\ldots,x_{q-1})$$
and
$$\partial^1_\nu\sigma(x_1,\ldots,x_{q-1})
  = \sigma(x_1,\ldots,x_{\nu-1},1,x_\nu,\ldots,x_{q-1}).$$

We use the following convention for orienting boundaries of manifolds.
Let $X$ be a manifold with boundary $\partial X$.
We orient the normal bundle $N\partial X$ of $\partial X$
by taking the outward direction to be positive.
Given an orientation of $X$ this determines an orientation of $\partial X$
by the requirement that $N\partial X\oplus T\partial X = TX$ as
oriented vector spaces.

\proclaim Proposition 3.6.5.

If $\sigma:[0,1]^q\to\eufJ^*_{g,k^\minus,k^\plus}(A)$ is semi-regular
and $\alpha^\pm_i\in\eusC^0_H$ for $i=1,\ldots,k^\plusminus$,
then
$$\hss
  \eqalign{
  & \smash{\sum_{j=1}^{k^\minus} \sum_{\alpha_0\in\eusC^0_H}}
      (-1)^{q+\mu(\alpha^\minus_1)+\cdots+\mu(\alpha^\minus_{j-1})} \cr
  & \qquad\qquad\qquad\quad
      \#\eusM^1(\alpha^\minus_j,\alpha_0) \;
      \#\eufM^0_\sigma
    ((\alpha^\minus_1
      ,\ldots,\alpha^\minus_{j-1},\alpha_0,\alpha^\minus_{j+1}
      ,\ldots,\alpha^\minus_{k^\minus}
     )
     ,(\alpha^\plus_1,\ldots,\alpha^\plus_{k^\plus})
    ) \cr
  & - \smash{\sum_{i=1}^{k^\plus} \sum_{\alpha_0\in\eusC^0_H}}
      (-1)^{q+\mu(\alpha^\minus_1)+\cdots+\mu(\alpha^\minus_{k^\minus})
            +\mu(\alpha^\plus_1)+\cdots+\mu(\alpha^\plus_{i-1})+\mu(\alpha_0)
           } \cr
  & \qquad\qquad\qquad\quad
    \#\eufM^0_\sigma
      ((\alpha^\minus_1,\ldots,\alpha^\minus_{k^\minus})
       ,(\alpha^\plus_1
         ,\ldots,\alpha^\plus_{i-1},\alpha_0,\alpha^\plus_{i+1}
         ,\ldots,\alpha^\plus_{k^\plus}
        )
      ) \;
      \#\eusM^1(\alpha_0,\alpha^\plus_i) \cr
  \noalign{\smallskip}
  & + \sum_{\nu=1}^q (-1)^{\nu+1}
    \Big(
      \#\eufM^0_{\partial^1_\nu\sigma}
      ((\alpha^\minus_1
        ,\ldots
        ,\alpha^\minus_{k^\minus}
       )
       ,(\alpha^\plus_1
         ,\ldots
         ,\alpha^\plus_{k^\plus}
        )
      ) \cr
  & \qquad\qquad\qquad\qquad
    - \#\eufM^0_{\partial^0_\nu\sigma}
      ((\alpha^\minus_1
        ,\ldots
        ,\alpha^\minus_{k^\minus}
       )
       ,(\alpha^\plus_1
         ,\ldots
         ,\alpha^\plus_{k^\plus}
        )
      )
    \Big) = 0  . \cr
    }
  \hss
$$

\proclaim Proposition 3.6.6.  \rm ([F4] Thm.~4.)

For $\alpha^\minus,\alpha^\plus\in\eusC^0_H$,
$$\sum_{\alpha_0\in\eusC_H}
    \#\eusM^1(\alpha^\minus,\alpha_0) \;
    \#\eusM^1(\alpha_0,\alpha^\plus) = 0 .
$$

\demo Proof of Thm.~3.6.3.

We denote the moduli spaces
$\eufM^{0,\plus}_{\sigma_1,i}[\alpha_0]$,
$\eufM^{0,\minus}_{\sigma_2,j}[\alpha_0]$
and
$\eufM^0_{\sigma_1\gop^\ell_{ij}\sigma_2}$.
In this notation, we have to show that for $\ell$ large enough,
$\eufM^0_{\sigma_1\gop^\ell_{ij}\sigma_2}$ is semi-regular,
and
$$\#\eufM^0_{\sigma_1\gop^\ell_{ij}\sigma_2}
  = \sum_{\alpha_0\in\eusC_H}
     \pm\#\eufM^{0,\plus}_{\sigma_1,i}[\alpha_0]
        \#\eufM^{0,\minus}_{\sigma_2,j}[\alpha_0].
$$
The moduli space
$\eufM^{0,\plus}_{\sigma_1,i}[\alpha_0]
 \times
 \eufM^{0,\minus}_{\sigma_2,j}[\alpha_0]
$
is the zero locus of the section
$\Psi^{0,\plus}_{\sigma_1,i}[\alpha_0]
 \boxplus
 \Psi^{0,\minus}_{\sigma_2,j}[\alpha_0]
$
of the bundle
$\IF^{0,\plus}_{\sigma_1,i}[\alpha_0]
 \boxplus
 \IF^{0,\minus}_{\sigma_2,j}[\alpha_0]
$
over
$\eufP^{0,\plus}_{\sigma_1,i}[\alpha_0]
 \times
 \eufP^{0,\minus}_{\sigma_2,j}[\alpha_0]
$.
By Prop.~2.3.6 this section is transverse to the zero section.
Let $\eufU$ be a neighborhood of
$\eufM^{0,\plus}_{\sigma_1,i}[\alpha_0]
 \times
 \eufM^{0,\minus}_{\sigma_2,j}[\alpha_0]
$
in
$\eufP^{0,\plus}_{\sigma_1,i}[\alpha_0]
 \times
 \eufP^{0,\minus}_{\sigma_2,j}[\alpha_0]
$.
It follows from the implicit function theorem
that if we choose $\eufU$ small enough,
then any sufficiently small perturbation of
the restriction of
$\Psi^{0,\plus}_{\sigma_1,i}[\alpha_0]
 \boxplus
 \Psi^{0,\minus}_{\sigma_2,j}[\alpha_0]
$
to $\eufU$
is also transverse to the zero section
and its zero locus is a small perturbation of
$\eufM^{0,\plus}_{\sigma_1,i}[\alpha_0]
 \times\eufM^{0,\minus}_{\sigma_2,j}[\alpha_0]
$.

By Prop.~3.4.3,
we can count the number of points in
$\eufM^0_{\sigma_1\gop^\ell_{ij}\sigma_2}
 \times
 \eusM^0(\alpha_0,\alpha)$
instead of
$\eufM^0_{\sigma_1\gop^\ell_{ij}\sigma_2}$.
It follows from Prop.~3.5.1 that for $\ell$
large enough,
$\eufM^{0,\plus}_{\sigma_1,i}[\alpha_0]
 \times
 \eufM^{0,\minus}_{\sigma_2,j}[\alpha_0]
$
is contained in the domain of $\bfg^\ell_{ij}$.
It follows from Prop.~3.5.3
that for $\ell$ large enough,
$\eufM^0_{\sigma_1\gop^\ell_{ij}\sigma_2}
 \times
 \eusM^0(\alpha_0,\alpha)$
lies in the range of $\bfg^\ell_{ij}$.
It also follows from Prop.~3.5.3 that
$(\bfg^\ell_{ij})^{-1}
 (\eufM^0_{\sigma_1\gop^\ell_{ij}\sigma_2}
  \times
  \eusM^0(\alpha_0,\alpha)
 )
\subset\eufU
$
for $\ell$ large enough.

Now
$\Psi^{0,\plus}_{\sigma_1,i}[\alpha_0]
 \boxplus
 \Psi^{0,\minus}_{\sigma_2,j}[\alpha_0]
$
and
$(\bfg^\ell_{ij})^*
 \bigl( \Psi^0_{\sigma_1\gop^\ell_{ij}\sigma_2}
        \boxplus\Psi^0(\alpha_0,\alpha_0)
 \bigr)
$
are equal off the necks.
By [R] Prop.~2.1
that on the necks,
with $u_1=\exp\xi_1$ and $u_2=\exp\xi_2$,
$$\eqalign{
  \Psi^{0,\plus}_{\sigma_1,i}[\alpha_0](u_1)
  \boxplus
 \Psi^{0,\minus}_{\sigma_2,j}[\alpha_0](u_2)
  = {}
  & (\exp\xi_1)_*
    ( d\xi_1/dt
      + \bfJ (D_{\alpha_0}\Phi_H) \xi_1
      + (\xi_1 \otimes \xi_1)\cdot f_1 ) \cr
  & \oplus  (\exp\xi_2)_*
    ( d\xi_2/dt
      + \bfJ (D_{\alpha_0}\Phi_H) \xi_2
      + (\xi_2\otimes\xi_2)\cdot f_2 ) \cr
              }
$$
where $f_1$ and $f_2$ are smooth functions of
$\xi$, $\nabla\xi$, $\eta$, and $\nabla\eta$.
Similarly, on the necks,
with $u_3=\exp\xi_3$ and $u_4=\exp\xi_4$,
$$\eqalign{
  & \Psi^0_{\sigma_1\gop^\ell_{ij}\sigma_2}
    \boxplus
    \Psi^0(\alpha_0,\alpha_0) \cr
  & = (\exp\xi_3)_*
      \big( d\xi_3/dt
            + \bfJ (D_{\alpha_0}\Phi_H) \xi_3
            + (\xi_3\otimes\xi_3)\cdot f
      \big) \cr
  & \phantom{{}={}}\oplus  (\exp\xi_4)_*
    \big( d\xi_4/dt
          + \bfJ (D_{\alpha_0}\Phi_H) \xi_4
          + (\xi_4\otimes\xi_4)\cdot f
    \big) . \cr
              }
$$
It follows that
$$\bigl( \Psi^{0,\plus}_{\sigma_1,i}[\alpha_0]
         \boxplus
         \Psi^{0,\minus}_{\sigma_2,j}[\alpha_0]
  \bigr)
  - (\bfg^\ell_{ij})^*
    \bigl( \Psi^0_{\sigma_1\gop^\ell_{ij}\sigma_2}
           \boxplus
           \Psi^0(\alpha_0,\alpha_0)
    \bigr)
$$
vanishes off the necks,
and on the necks it is given by terms of the form
$(\xi_\mu\otimes\xi_\nu)\cdot f$ and $(d\eta/dt)\,\xi_\mu\cdot f$.
We can make these terms,
and their first order Fr\'echet derivatives,
arbitrarily small by choosing $\eufU$ small and $\ell$ large.
It then follows from the implicit function theorem that
$\eufM^{0,\plus}_{\sigma_1,i}[\alpha_0]
 \times
 \eufM^{0,\minus}_{\sigma_2,j}[\alpha_0]
$
is a small perturbation of
$(\bfg^\ell_{ij})^{-1}
 \bigl(
   \eufM^0_{\sigma_1\gop^\ell_{ij}\sigma_2}
   \times
   \eusM^0(\alpha_0,\alpha_0)
 \bigr)
$.
The Proposition follows,
if we disregard signs.

Finally,
the excision map $\bfg^\ell_{ij}$ induces a map \qck.
By Prop.~3.4.3
this map changes the coherent orientations by
the graded sign of the permutation \qcm.
By Prop.~3.3.3,
if the moduli spaces are non-empty,
then
$$\eqalign{
    q_1 \equiv {} &
         \muH(\alpha^\minus_{1,1})
         + \cdots
         + \muH(\alpha^\minus_{1,k^\minus_1}) \cr
    & \qquad
         + \muH(\alpha^\plus_{1,1})
         + \cdots
         + \muH(\alpha^\plus_{1,i-1})
         + \muH(\alpha^\plus_0)
         + \muH(\alpha^\plus_{1,i+1})
         + \cdots
         + \muH(\alpha^\plus_{1,k^\plus_1+1}) \cr
          }
$$
modulo 2.
It follows that the signs of the permutations \qcm{} and \qci{}
differ by $(-1)^{q_1q_2}$.
This accounts for the signs in the Proposition.
\enddemo

\demo Proof of Thm.~3.6.4.

This is proven the same way as Thm.~3.6.3.
\enddemo

\demo Proof of Thm.~3.6.5.

For simplicity we denote the moduli spaces
$\eufM^{0,\plusminus}_{\sigma,i}[\alpha_0]$,
$\eufM^0_{\partial^0_\nu\sigma}$,
$\eufM^0_{\partial^1_\nu\sigma}$.
In this notation we have to show that
$$\eqalign{
  & \sum_{j,\alpha_0}\pm\#\eusM^1(\alpha^\minus_j,\alpha_0)\,
        \#\eufM^{0,\minus}_{\sigma,j}[\alpha_0]
    - \sum_{i,\alpha_0}\pm\#\eufM^{0,\plus}_{\sigma,i}[\alpha_0]\,
          \#\eusM^1(\alpha_0,\alpha^\plus_i) \cr
  & + \sum_\nu(-1)^{\nu+1}\left(
       \#\eufM^0_{\partial^1_\nu\sigma}
       - \#\eufM^0_{\partial^0_\nu\sigma}
                    \right)
    = 0 . \cr
          }
  \eqno\qcq
$$
We also let
$$\eufM^1_\sigma
  = \eufM^1_\sigma
  ((\alpha^\minus_1,\ldots,\alpha^\minus_{k^\minus}),
   (\alpha^\plus_1,\ldots,\alpha^\plus_{k^\plus})
   ) .
$$
At first we assume that $\sigma$ is transverse to
$\pi:\eufM^{1-q}_{g,k^\minus,k^\plus}\to\eufJ^*_{g,k^\minus,k^\plus}(A)$.
At the end of the proof we will show how to eliminate this assumption.
The idea behind the proof is to identify the boundary of $\eufM^1_\sigma$,
and use the fact that the number of boundary points counted with
signs is zero.
We also have to consider that $\eufM^1_\sigma$
may be non-compact,
and include the number of ends, also counted with signs.
This accounts for the first two sums in \qcq.

The map $\sigma$,
and the boundary maps $\partial^0_\nu\sigma$ and $\partial^1_\nu\sigma$,
are transverse to
$\pi:\eufM^{1-q}_{g,k^\minus,k^\plus}\to\eufJ^*_{g,k^\minus,k^\plus}(A)$.
It follows that the moduli space $\eufM^1_\sigma$
is a smooth 1-manifold with boundary,
and its boundary points are given by
the moduli spaces $\eufM^0_{\partial^0_\nu\sigma}$
and $\eufM^0_{\partial^1_\nu\sigma}$.
This accounts for the third sum in \qcq.

The moduli spaces $\eufM^0_{\partial^0_\nu\sigma}$
and $\eufM^0_{\partial^1_\nu\sigma}$
have two natural orientations,
the coherent orientation
and the orientation as boundary of $\eufM^1_\sigma$.
We have
$\eufO_{\partial^0_\nu \sigma}
 = \eufO_{[0,1]^{q-1}} \otimes \sigma^*\eufO_{g,k^\minus,k^\plus}$.
The coherent orientation of $\eufM^0_{\partial^0_\nu\sigma}$
is given by the canonical orientation of $[0,1]^{q-1}$
and the coherent orientation of $\eufO_{g,k^\minus,k^\plus}$.
The boundary orientation of $\eufM^0_{\partial^0_\nu\sigma}$
is given by the boundary orientation of $[0,1]^{q-1}$
and the coherent orientation of $\eufO_{g,k^\minus,k^\plus}$.
These orientations differ by $(-1)^\nu$ for $\partial^0_\nu\sigma$
and by $(-1)^{\nu+1}$ for $\partial^1_\nu\sigma$.
This accounts for the signs in the third sum in \qcq.

Define
$\eufM^{1,\alpha_0,\ell,\plus}_{\sigma,i}$
by
$$\eufM^{1,\alpha_0,\ell,\plus}_{\sigma,i} \times \eusM^0(\alpha_0,\alpha_0)
  = (\eufM^1_\sigma \times \eusM^0(\alpha_0,\alpha_0) )
    \cap \range \bfg^\ell_{i1}.
$$
where $\bfg^\ell_{i1}$ is the excision map
$$(\bfg^\ell_{i1})
  : \eufP^{0,\plus}_{\sigma,i}[\alpha_0]\times\eusP^1(\alpha_0,\alpha^\plus_i)
    \to \eufP^1_\sigma\times\eusP^0(\alpha_0,\alpha_0).$$
Similarly,
define $\eusM^{1,\ell,\minus}(\alpha_0,\alpha^\plus_i)$
by
$$\eufM^{0,\plus}_{\sigma,i}[\alpha_0]
  \times\eusM^{1,\ell,\minus}(\alpha_0,\alpha^\plus_i)
  = (\eufM^{0,\plus}_{\sigma,i}[\alpha_0]
     \times\eusM^1(\alpha_0,\alpha^\plus_i)
    ) \cap \domain \bfg^\ell_{i1} .
$$
It follows from Prop.~3.5.2 that for $\ell$ large enough,
$$\eufM^1_\sigma
  \setminus
  \bigcup_{\plusminus,i,\alpha_0}
\eufM^{1,\alpha_0,\ell,\plusminus}_{\sigma,i}$$
is a compact subset of $\eufM^1_\sigma$.
Thus it suffices to count the number of ends of each
$\eufM^{1,\alpha_0,\ell,\plusminus}_{\sigma,i}$.

First we count the number of ends of
$\eufM^{1,\alpha_0,\ell,\plus}_{\sigma,i}$.
The moduli space
$\eufM^{0,\plus}_{\sigma,i}[\alpha_0]
 \times \eusM^1(\alpha_0,\alpha^\plus_i)
$
is a union of affine lines.
It is the zero locus of a translation invariant section
$\Psi^{0,\plus}_{\sigma,i}[\alpha_0]
 \boxplus
 \Psi^1(\alpha_0,\alpha^\plus_i) .
$
By Prop.~2.3.6,
this section is transverse to the zero section.
Hence there exists a translation invariant neighborhood $\eufU$ of
$\eufM^{0,\alpha_0,\plus}_{\sigma,i}
  \times
  \eusM^1(\alpha_0,\alpha^\plus_i)
$
in
$\eufP^{0,\plus}_{\sigma,i}[\alpha_0]
 \boxplus
 \eusP^1(\alpha_0,\alpha^\plus_i)
$
such that  the zero locus of any sufficiently small perturbation of
$\Psi^{0,\plus}_{\sigma,i}[\alpha_0]
 \times
 \Psi^1(\alpha_0,\alpha^\plus_i)
$
is a small perturbation of
$\eufM^{0,\alpha_0,\plus}_{\sigma,i}
  \times
  \eusM^1(\alpha_0,\alpha^\plus_i)
$.

It follows from Prop.~3.5.2 that
$(\bfg^\ell_{i1})^{-1}
 (\eufM^{1,\alpha_0,\ell,\plus}_{\sigma,i} \times \eusM^0(\alpha_0,\alpha_0))
 \subset\eufU
$
for $\ell$ large enough.
Now
$(\bfg^\ell_{i1})^{-1}
 (\eufM^{1,\alpha_0,\ell,\plus}_{\sigma,i} \times \eusM^0(\alpha_0,\alpha_0))
$
is the zero locus of
$(\bfg^\ell_{i1})^*
 \bigl( \Psi^{1,\alpha_0,\ell,\plus}_{\sigma,i}
        \boxplus
        \Psi^0(\alpha_0,\alpha_0)
 \bigr)
$.
As in the proof of Prop.~3.6.3,
these two sections of
$\IF^{0,\alpha_0,\plus}_{\sigma,i}
  \boxplus
  \IF^1(\alpha_0,\alpha^\plus_i)
$
differ by terms of the form $(d\eta/dt)\xi_\mu\cdot f$
and $(\xi_\mu\otimes\xi_\nu)\cdot f$ supported on the neck.
These terms, and their first order Fr\'echet derivatives,
can be made arbitrarily small by choosing $\eufU$ small and $\ell$ large.
It then follows from the implicit function theorem that
$(\bfg^\ell_{i1})^{-1}
 (\eufM^{1,\alpha_0,\ell,\plus}_{\sigma,i} \times \eusM^0(\alpha_0,\alpha_0))$
is a small perturbation of
$\eufM^{0,\plus}_{\sigma,i}[\alpha_0]
 \times \eusM^1(\alpha_0,\alpha^\plus_i).
$
Thus the ends of
$(\eufM^{1,\alpha_0,\ell,\plus}_{\sigma,i} \times \eusM^0(\alpha_0,\alpha_0))$
can be identified with the ends of
$\eufM^{0,\plus}_{\sigma,i}[\alpha_0]
 \times \eusM^1(\alpha_0,\alpha^\plus_i).
$
This accounts for the first sum in \qcq.

It follows from Prop.~3.4.3
that the map $\bfg^\ell_{i1}$ changes the orientations by
the graded sign of the permutation
$$\eqalign{
  & (\sigma,
     \alpha^\minus_1,\ldots,\alpha^\minus_{k^\minus},
     \alpha^\plus_1,\ldots,
     \alpha^\plus_{i-1},\mathop{\alpha_0}_\plus,\alpha^\plus_{i+1},\ldots,
     \alpha^\plus_{k^\plus},
     \mathop{\alpha_0}_\minus,\alpha^\plus_i
    ) \cr
  & \mapsto
    (\sigma,
     \alpha^\minus_1,\ldots,\alpha^\minus_{k^\minus},
     \alpha^\plus_1,\ldots,\alpha^\plus_{k^\plus},
     \mathop{\alpha_0}_\plus,\mathop{\alpha_0}_\minus
    ) , \cr
          }
$$
which is
$$(-1)^{\mu(\alpha^\plus_{i+1})+\cdots+\mu(\alpha^\plus_{k^\plus})}
  = (-1)^{q+\mu(\alpha^\minus_1)+\cdots+\mu(\alpha^\minus_{k^\minus})
          +\mu(\alpha^\plus_1)+\cdots+\mu(\alpha^\plus_{i-1})
          +\mu(\alpha_0)
         } .
$$
This accounts for the signs in the first sum in \qcq.
We get an additional minus sign as the
number of ends of $\eusM^{1,\ell,-}(\alpha_0,\alpha^\plus_i)$,
counted with signs,
is $-(\#\eusM^1(\alpha_0,\alpha^\plus_i))$.

Similarly, the second sum in \qcq{} gives the number of ends of
$\eufM^{1,\alpha_0,\ell,\minus}_{\sigma,j}$.
The sign of this term is given by the graded sign of the permutation
$$\eqalign{
  & (\alpha^\minus_j,\mathop{\alpha_0}_\plus,
     \sigma,
     \alpha^\minus_1,\ldots,
     \alpha^\minus_{j-1},\mathop{\alpha_0}_\minus,\alpha^\minus_{j+1},\ldots,
     \alpha^\minus_{k^\minus},
     \alpha^\plus_1,\ldots,\alpha^\plus_{k^\plus}
    ) \cr
  & \mapsto
    (\sigma,
     \alpha^\minus_1,\ldots,\alpha^\minus_{k^\minus},
     \alpha^\plus_1,\ldots,\alpha^\plus_{k^\plus},
     \mathop{\alpha_0}_\plus,\mathop{\alpha_0}_\minus) , \cr
          }
$$
which is
$$(-1)^{q+\mu(\alpha^\minus_1)+\cdots+\mu(\alpha^\minus_{j-1})}.$$

Finally, assume that $\sigma$ is not transverse to
$\pi:\eufM^{1-q}_{g,k^\minus,k^\plus}\to\eufJ^*_{g,k^\minus,k^\plus}(A)$.
By assumption $\partial^0_\nu\sigma$ and $\partial^1_\nu\sigma$
are semi-regular and hence transverse to
$\pi:\eufM^{1-q}_{g,k^\minus,k^\plus}\to\eufJ^*_{g,k^\minus,k^\plus}(A)$.
Then there exists a sequence of maps
$\sigma_n:[0,1]^q\to\eufJ^*_{g,k^\minus,k^\plus}(A)$ such that
$\sigma_n$ is transverse to
$\pi:\eufM^{1-q}_{g,k^\minus,k^\plus}\to\eufJ^*_{g,k^\minus,k^\plus}(A)$,
$\sigma_n\to\sigma$ as $n\to\infty$,
$\partial^0_\nu\sigma_n=\partial^0_\nu\sigma$, and
$\partial^1_\nu\sigma_n=\partial^1_\nu\sigma$.
The Proposition holds with $\sigma$ replaced by $\sigma_n$.
Thus we only have to verify that
$\#\eufM^{0,\plusminus}_{\sigma_n,i}[\alpha_0]
 = \#\eufM^{0,\plusminus}_{\sigma,i}[\alpha_0]
$
for $n$ large enough.

Let $\eufV$ be a neighborhood of
$\eufM^{0,\plusminus}_{\sigma,i}[\alpha_0]$
in
$\eufP^{0,\plusminus}_{\sigma,i}[\alpha_0]$.
The moduli spaces
$\eufM^{0,\plusminus}_{\sigma_n,i}[\alpha_0]$
and
$\eufM^{0,\plusminus}_{\sigma,i}[\alpha_0]$
are the zero loci of sections
$\Psi^{0,\plusminus}_{\sigma_n,i}[\alpha_0]$
and
$\Psi^{0,\plusminus}_{\sigma,i}[\alpha_0]$.
It follows from the implicit function theorem that if we choose
$\eufV$ small enough and $n$ large enough,
then
$\eufM^{0,\plusminus}_{\sigma_n,i}[\alpha_0]\cap\eufV$
is a small perturbation of
$\eufM^{0,\plusminus}_{\sigma,i}[\alpha_0]\cap\eufV$.
On the other hand,
in Prop.~3.5.1 we can let $u_n\in\eufM^0_{\sigma_n}$
where $\sigma_n\to\sigma$ as $n\to\infty$.
It follows that for any $\eufV$,
$\eufM^{0,\plusminus}_{\sigma_n,i}[\alpha_0]
 \subset\eufV$
for $n$ large enough.
\enddemo

\demo Remarks on Prop.~3.6.6.

A proof of this is outlined in [F4].
No complete proof has appeared in print.
A gap in [F4],
the possibility that, when $N_0=1$, a sequence in $\eusM^2(\alpha,\alpha)$
degenerates to an element of $\eusM^0(\alpha,\alpha)$
and a holomorphic sphere with $\c_1=1$,
was bridged in [HS]
by adding an additional transversality condition;
see Definition 2.1.5.

Our scheme for orienting the moduli spaces differs from the method
used by Floer,
so we need to check the signs in Prop.~3.6.6.
We get the same signs as Floer, for by Prop.~3.4.3,
the excision maps
$$(\bfg^\ell_{11})_* :
 \eusP^1(\alpha^\minus,\alpha_0)
 \times \eusP^1(\alpha_0,\alpha^\plus)
 \to
 \eusP^2(\alpha^\minus,\alpha^\plus)
  \times \eusP^0(\alpha_0,\alpha_0)
$$
preserve the coherent orientations.
\enddemo

%%%%%%%%%%%%%%%%%%%%%%%%%%%%%%%%%%%%%%%%%%%%%%%%%%%%%%%%%%%%%%%%

\heading \S4. Products and Relations

\subheading 4.1.  The homological gluing maps $\gop_{ij}$ and $\cop_{ij}$.

The singular homology groups
$H_*(\eufJ^*_{g,k^\minus,k^\plus}(A),\IZ)$
can be defined as the homology of the complex
$C_*(\eufJ^*_{g,k^\minus,k^\plus}(A),\IZ)$
of smooth cubical simplicial chains,
modulo degenerate chains; see [HW] Sect.~8.3.
These chains are linear combinations with integer
coefficients of smooth cubical simplices,
$\sigma:[0,1]^q \to \eufJ^*_{g,k^\minus,k^\plus}(A)$.
The boundary operator is defined by
$$d\sigma
  = \sum_{\nu=1}^q (-1)^{\nu+1}
    ( \partial^1_\nu\sigma - \partial^0_\nu\sigma) .
$$

If $k^\plus_1+k^\minus_1+k^\plus_2+k^\minus_2\ge1$, then
$\eufJ^*_{g_1+g_2,k_1^\minus+k_2^\minus,k_1^\plus+k^\plus_2}(A)
 =\eufJ_{g_1+g_2,k_1^\minus+k_2^\minus,k_1^\plus+k^\plus_2}(A)
$,
and the gluing map \qca{} induces a chain map
$$\gop^\ell_{ij} : C_*(\eufJ^*_{g_1,k_1^\minus,k_1^\plus1}(A))
              \otimes C_*(\eufJ^*_{g_2,k^\minus_2+1,k^\plus_2}(A))
              \to
C_*(\eufJ^*_{g_1+g_2,k_1^\minus+k_2^\minus,k_1^\plus+k^\plus_2}(A)) ,
$$
and hence a homomorphism
$$\gop_{ij} : H_*(\eufJ^*_{g_1,k_1^\minus,k_1^\plus1}(A))
              \otimes H_*(\eufJ^*_{g_2,k^\minus_2+1,k^\plus_2}(A))
              \to
H_*(\eufJ^*_{g_1+g_2,k_1^\minus+k_2^\minus,k_1^\plus+k^\plus_2}(A)).
$$
Similarly,
if $k^\minus+k^\plus\ge1$ then the gluing map \qcaa{} induces
a chain map
$$\cop^\ell_{ij} : C_*(\eufJ^*_{g,k^\minus+1,k^\plus+1}(A))
              \to C_*(\eufJ^*_{g+1,k^\minus,k^\plus}(A)) ,
$$
and hence a homomorphism
$$\cop_{ij} : H_*(\eufJ^*_{g,k^\minus+1,k^\plus+1}(A))
              \to H_*(\eufJ^*_{g+1,k^\minus,k^\plus}(A)).
$$
The homology maps are independent of $\ell$.

There are natural maps
$$\gop_{11} : H_*(\eufJ^*_{g_1,0,1}(A))
              \otimes H_*(\eufJ^*_{g_2,1,0}(A))
              \to H_*(\eufJ^*_{g_1+g_2,0,0}(A))
$$
and
$$\cop_{11} : H_*(\eufJ^*_{g,1,1}(A))
              \to H_*(\eufJ^*_{g+1,0,0}(A))
$$
as well.
These are defined using the following Lemma.

\proclaim Lemma 4.1.1.

{
For any $\sigma \in C_q(\eufJ^*_{g,1,1}(A))$
there exists
$\sigma' \in C_q(\eufJ^*_{g,1,1}(A))$
and $\tau \in C_{q+1}(\eufJ^*_{g,1,1}(A))$
such that $\sigma-\sigma' = d\tau$
and $\cop^\ell_{11}\sigma'\in C_q(\eufJ^*_{g+1,0,0}(A))$ for all $\ell\ge0$.

Similarly,
for any $\sigma_1 \in C_{q_1}(\eufJ^*_{g_1,0,1}(A))$ and
$\sigma_2 \in C_{q_2}(\eufJ^*_{g_2,1,0}(A))$
there exist
$\sigma'_1 \in C_{q_1}(\eufJ^*_{g_1,0,1}(A))$
and $\tau \in C_{q_1+1}(\eufJ^*_{g_1,k^\minus_1,k^\plus_1}(A))$
such that $\sigma_1-\sigma'_1 = d\tau$
and
$\sigma'_1\gop^\ell_{11}\sigma_2\in C_{q_1+q_2}(\eufJ^*_{g_1+g_2,0,0}(A))$
for all $l\ge0$.
}

\demo Proof of Lemma 4.1.1.

We say that $c=[\bfj,R,\Delta]\in\eufJ^*_{g,1,1}(A)$ is somewhere injective
if it has the following property.
There exists an open subset $U$ of $\Sigma_{g,1,1}$ such that
for any $p\in U$ there does not exist $q\in\Sigma_{g,1,1}$
and a complex linear map $L:(T^{0,1}_p\Sigma)^*\to (T^{0,1}_q\Sigma_{g,1,1})^*$
such that for all $m\in M$, $R(q,m)=(L\otimes{\rm id})R(p,m)$.
If $c$ is somewhere injective, $\phi_*\bfj=\bfj$ and $\phi_*R=R$,
then $\phi$ is the identity on $U$.
By analytic continuation, $\phi$ is the identity
on all of $\Sigmabar_{g,1,1}$.
In other words, if $c$ is somewhere injective,
then $c\in\eufJ^*_{g,1,1}(A)$.

The complement of the somewhere injective elements
has infinite codimension in the following sense.
Any map
from a finite dimensional manifold into $\eufJ^*_{g,1,1}(A)$,
such that the boundary values are somewhere injective,
can be perturbed in the interior
so that all its values are somewhere injective.
If $c$ is somewhere injective, then
$\cop^\ell_{ij}c$ is somewhere injective.
The first part of the lemma follows.

The second part is proven in a similar way.
\enddemo

\subheading  4.2. Floer (co)homology and the contraction maps $\gop_{ij}$ and
$\cop_{ij}$.

We first recall a few facts from homological algebra.
The tensor product $C_* \otimes C'_*$
of two chain complexes $C_*$ and $C'_*$ is the
chain complex defined by
$(C\otimes C')_k = \bigoplus_{i+j=k} C_i\otimes C'_j,$
and
$d(x \otimes x') = dx \otimes x' + (-1)^i x \otimes dx'$
for $x\otimes x' \in C_i\otimes C'_j$.
With this choice of signs,
$(C_*\otimes C'_*)\otimes C''_*
 = C_*\otimes(C'_*\otimes C''_*)$
and
$\underline\IZ_*\otimes C_*
 = C_*\otimes\underline\IZ_*
 = C_*
$,
where $\underline\IZ_*$ is the chain complex given by $\underline\IZ_0=\IZ$
and $\underline\IZ_i=(0)$ for $i\ne0$.

The dual $C^*$ of a complex $C_*$ is defined by
$C^i=\Hom(C_i,\IZ)$
and
$d^*x = (-1)^i d^\dagger x$
for $x\in C^i$
where $d^\dagger$ is the adjoint of $d$.
With this choice of signs,
the contraction
$C_*\otimes C^*\to\underline\IZ_*$
is a chain map.

If $\rho$ is a permutation on $k$ letters,
then there is a chain map
$C^{(1)}_*\otimes\cdots\otimes C^{(k)}_*
  \to C^{(\rho(1))}_* \otimes\cdots\otimes C^{(\rho(k))}_*
$
defined by
$\alpha_1\otimes\cdots\otimes\alpha_k
  \mapsto\pm\alpha_{\rho(1)}\otimes\cdots\otimes\alpha_{\rho(k)}
$
where the sign is the graded sign of the permutation
$(\alpha_1,\ldots,\alpha_k)
  \to(\alpha_{\rho(1)},\ldots,\alpha_{\rho(k)}).
$
(Recall that the graded sign of a permutation is defined
as the sign of the permutation obtained by removing
all elements of even degree.)
To show this,
it suffices to check that the homomorphism
$C_*\otimes C'_*\to C'_*\otimes C_*$ given by
$\alpha\otimes\beta\mapsto(-1)^{ij}\beta\otimes\alpha$
for $\alpha\in C_i$ and $\beta\in C'_j$ is a chain map.
In particular, the symmetric group $S_k$ acts on $C_*^{\otimes k}$
by chain maps.

Fix a triple $A=(\bfJ,H,\eufo)$ where
$\bfJ$ is an almost complex structure on $M$ compatible with $\omega$,
$H$ is a time-dependent Hamiltonian on $M$
and $\eufo$ is a coherent system of orientations.
We assume that the pair $(\bfJ,H)$ is regular.

\proclaim Definition 4.2.1.

{
The Floer chain complex $CF_*(M,A)$ and the Floer cochain complex $CF^*(M,A)$
are $\IZ/2N_0\IZ$ graded and defined as follows.
For $p\in\IZ/2N_0\IZ$,
$CF_p(M,A)$ and $CF^p(M,A)$ are both the free $\IZ$ module generated by
periodic orbits $\alpha\in\eusC^0_H$ with $\muH(\alpha)=p$.
The differential
$$d:CF_i(M,A)\to CF_{i-1}(M,A)$$
is defined by
$$d\alpha
= \sum_{\beta \in \eusC^0_H} (-1)^{\mu(\beta)}\#\eusM^1(\alpha,\beta) \,
\beta.$$
The differential
$$d^*:CF^i(M,A)\to CF^{i+1}(M,A)$$
is defined by
$$d^*\alpha
= \sum_{\beta \in \eusC^0_H} \#\eusM^1(\beta,\alpha) \, \beta.$$
The pairing
$$CF_p(M,A) \otimes CF^p(M,A) \to \IZ$$
is given by
$$\langle \alpha , \beta \rangle
  = \cases { 1 & if $\alpha=\beta$ \cr 0 & otherwise \cr }
$$
for $\alpha, \beta \in \eusC^0_H$.
The Floer homology groups $HF_*(M,A)$
are the homology groups of the chain complex $CF_*(M,A)$.
The Floer cohomology groups $HF^*(M,A)$
are the cohomology groups of $CF^*(M,A)$.
}

This definition is consistent for
it follows from Prop.~3.3.3 that if $\eusM^1(\alpha,\beta)$ is non-empty,
then $\mu(\beta)=\mu(\alpha)-1$,
and it follows from Prop.~3.6.6 that $d^2=0$.
Note that $CF^*(M,A)$ is the dual of $CF_*(M,A)$.

In the following,
we will view $CF^*(M,A)$ as a chain complex,
where $CF^q(M,A)$ has degree $-q$,
rather than as a cochain complex.

\proclaim Definition 4.2.4.

{
The homomorphism
$$\eqalign{
    \gop_{ij}:{}
  &  CF^*(M,A)^{\otimes k^\minus_1}
     \otimes CF_*(M,A)^{\otimes k^\plus_1+1}
     \otimes CF^*(M,A)^{\otimes k^\minus_2+1}
     \otimes CF_*(M,A)^{\otimes k^\plus_2} \cr
  &  \to CF^*(M,A)^{\otimes(k_1^\minus+k_2^\minus)}
         \otimes CF_*(M,A)^{\otimes(k_1^\plus+k_2^\plus)} \cr
         }
$$
is defined by
$$ \eqalign{
&    \alpha^\minus_{1,1} \otimes \ldots \otimes \alpha^\minus_{1,k^\minus_1}
     \otimes
     \alpha^\plus_{1,1} \otimes \ldots \otimes \alpha^\plus_{1,k^\plus_1+1}
     \gop_{ij}
     \alpha^\minus_{2,1} \otimes \ldots \otimes \alpha^\minus_{2,k^\minus_2+1}
     \otimes
     \alpha^\plus_{2,1} \otimes \ldots \otimes \alpha^\plus_{2,k^\plus_2} \cr
& =  \pm \langle \alpha^\plus_{1,i},\alpha^\minus_{2,j} \rangle \,
     \alpha^\minus_{2,1} \otimes \ldots \otimes \alpha^\minus_{2,j-1}
     \otimes \alpha^\minus_{1,1} \otimes \ldots \otimes
\alpha^\minus_{1,k^\minus_1}
     \otimes \alpha^\minus_{2,j+1} \otimes \ldots \otimes
\alpha^\minus_{2,k^\minus_2+1} \cr
& \phantom{{}=  \pm \langle \alpha^\plus_{1,i},\alpha^\minus_{2,j} \rangle}
     \otimes \alpha^\plus_{1,1} \otimes \ldots \otimes \alpha^\plus_{1,i-1}
     \otimes \alpha^\plus_{2,1} \otimes \ldots \otimes
\alpha^\plus_{2,k^\plus_2}
     \otimes \alpha^\plus_{1,i+1} \otimes \ldots \otimes
\alpha^\plus_{1,k^\plus_1+1} \cr
         }$$
where the sign is the graded sign of the permutation
obtained by comparing the ordering of the $\alpha$'s
on the left and the right hand sides.

The homomorphism
$$\cop_{ij}
  :   CF^*(M,A)^{\otimes(k^\minus+1)}
      \otimes CF_*(M,A)^{\otimes(k^\plus+1)}
  \to CF^*(M,A)^{\otimes k^\minus}
      \otimes CF_*(M,A)^{\otimes k^\plus}
$$
is defined by
$$\eqalign{
  & \cop_{ij}
    \left(\alpha^\minus_1 \otimes \ldots \otimes \alpha^\minus_{k-} \otimes
    \alpha^\plus_1 \otimes \ldots \otimes \alpha^\plus_{k^\plus} \right) \cr
& =  \pm \langle \alpha^\plus_i,\alpha^\minus_j \rangle\,
     \alpha^\minus_1 \otimes \ldots \otimes
     \alpha^\minus_{j-1} \otimes \alpha^\minus_{j+1} \otimes \ldots
     \otimes \alpha^\minus_{k^\minus} \cr
& \phantom{{}=  \pm \langle \alpha^\plus_i,\alpha^\minus_j \rangle }
     \otimes \alpha^\plus_1 \otimes \ldots \otimes
     \alpha^\plus_{i-1} \otimes \alpha^\plus_{i+1} \otimes\ldots
     \otimes \alpha^\plus_{k^\plus} \cr
          }$$
where the sign again is the graded sign of the permutation
obtained by comparing the ordering of the $\alpha$'s
on the left and the right hand sides.
}

\proclaim Lemma 4.2.5.

The homomorphisms $\gop_{ij}$ and $\cop_{ij}$ are chain maps.

\demo Proof.

The homomorphisms $\gop_{ij}$ and $\cop_{ij}$
are compositions of the permutation maps
and contraction maps discussed at the beginning of this section.
\enddemo

\subheading 4.3. The chain map $Q$.

Let
$C_*^\pi(\eufJ^*_{g,k^\minus,k^\plus}(A))$
be the sub-complex of
$C_*(\eufJ^*_{g,k^\minus,k^\plus}(A))$
generated by semi-regular cubical simplices $\sigma$ as in Def.~2.3.4.
If the boundary faces of a simplex are semi-regular,
then the simplex can be perturbed,
keeping the boundary fixed,
to be semi-regular.
It follows that the inclusion
$$C_*^\pi(\eufJ^*_{g,k^\minus,k^\plus}(A)) \to
C_*(\eufJ^*_{g,k^\minus,k^\plus}(A))$$
is a chain homotopy equivalence.
In the following we will work with the chain complex
$C_*^\pi(\eufJ^*_{g,k^\minus,k^\plus}(A))$
rather than $C_*(\eufJ^*_{g,k^\minus,k^\plus}(A))$.
Recall that for semi-regular simplices $\sigma$, the numbers
$\#\eufM^0_\sigma((\alpha^\minus_1,\ldots,\alpha^\minus_{k^\minus}),
              (\alpha^\plus_1,,\ldots,\alpha^\plus_{k^\plus})
             )
$
are well defined.

\proclaim Definition 4.3.1.

The homomorphism
$$Q:C_*^\pi(\eufJ^*_{g,k^\minus,k^\plus}(A))
\to CF^*(M,A)^{\otimes k^\minus} \otimes CF_*(M,A)^{\otimes k^\plus}$$
is defined by
$$\eqalign{
    Q\sigma
  & = (-1)^{{1\over2}q(q-1)}
      \smash{\sum_{\alpha^\pm_*\in\eusC^0_H}}
      \#\eufM^0_\sigma((\alpha^\minus_1,\ldots,\alpha^\minus_{k^\minus}),
                 (\alpha^\plus_1,\ldots,\alpha^\plus_{k^\plus})
                ) \cr
  & \qquad\qquad\qquad\qquad\qquad \,
           \alpha^\minus_1\otimes\cdots\otimes\alpha^\minus_{k^\minus}
           \otimes\alpha^\plus_1\otimes\cdots\otimes\alpha^\plus_{k^\plus}. \cr
          }
$$

If $\sigma$ is degenerate and semi-regular,
then the moduli spaces $\eufM^0_\sigma$ are empty.
Thus $Q$ is well defined.

\proclaim Theorem 4.3.2.

If $g=0$, then $Q$ is a chain map of degree $2n(1-k^\minus)$.
In general, $Q$ is a chain map of mixed degree $2n(1-g-k^\minus)+2\nu N_1$,
$\nu\in\IZ$.
In particular, there is an induced homomorphism,
$$Q :   H_*(\eufJ^*_{g,k^\minus,k^\plus}(A))
    \to H(CF^*(M,A)^{\otimes k^\minus} \otimes CF_*(M,A)^{\otimes k^\plus}) .
$$

\demo Proof.

It follows from Prop.~3.6.5 that if
$\sigma:[0,1]^q\to \eufJ^*_{g,k^\minus,k^\plus}(A)$ is smooth and semi-regular,
then
$$\hss\eqalign{
\noalign{\medskip}
  dQ\sigma
  & = \smash{\sum_{j=1}^{k^\minus}}\,
      \smash{\sum_{\ss\alpha^\pm_*,\alpha_0\in\eusC^0_H
                   \atop
                   \ss{\rm omit}\>\alpha^\minus_j
                  }
            }
    (-1)^{{1\over2}q(q-1)+\mu(\alpha^\minus_1)
    +\cdots+\mu(\alpha^\minus_{j-1})} \cr
  & \qquad\qquad\qquad\qquad
      \#\eufM^0_\sigma((\alpha^\minus_1,\cdots,\alpha^\minus_{j-1},\alpha_0,
                  \alpha^\minus_{j+1},\ldots,\alpha^\minus_{k^\minus}
                 ),
                 (\alpha^\plus_1,\cdots,\alpha^\plus_{k^\plus})
                ) \cr
  & \qquad\qquad\qquad\qquad
         \alpha^\minus_1\otimes\cdots\otimes\alpha^\minus_{j-1}
         \otimes d^*\alpha_0\otimes\alpha^\minus_{j+1}
         \otimes\cdots\otimes\alpha^\minus_{k^\minus}
         \otimes\alpha^\plus_1\otimes\cdots\otimes\alpha^\plus_{k^\plus} \cr
\noalign{\smallskip}
  & + \smash{\sum_{i=1}^{k^\plus}}\,
      \smash{\sum_{\ss\alpha^\pm_*,\alpha_0\in\eusC^0_H
                   \atop
                   \ss{\rm omit}\>\alpha^\plus_i
                  }
            }
    (-1)^{{1\over2}q(q-1)+\mu(\alpha^\minus_1)
            +\cdots+\mu(\alpha^\minus_{k^\minus})
            +\mu(\alpha^\plus_1)+\cdots+\mu(\alpha^\plus_{i-1})
           } \cr
  & \qquad\qquad\qquad\qquad
      \#\eufM^0_\sigma((\alpha^\minus_1,\cdots,\alpha^\minus_{k^\minus}),
                       (\alpha^\plus_1,\cdots,\alpha^\plus_{i-1},\alpha_0,
                        \alpha^\plus_{i+1},
                        \ldots,\alpha^\plus_{k^\plus})
                ) \cr
  & \qquad\qquad\qquad\qquad
         \alpha^\minus_1\otimes\cdots\otimes\alpha^\minus_{k^\minus}
         \otimes\alpha^\plus_1\otimes\cdots\otimes\alpha^\plus_{i-1}
         \otimes d\alpha_0\otimes\alpha^\plus_{i+1}\cdots
         \otimes\alpha^\plus_{k^\plus} \cr
\noalign{\smallskip}
  & = \smash{\sum_{j=1}^{k^\minus}}\,
      \smash{\sum_{\alpha^\pm_*,\alpha_0\in\eusC^0_H}}
      (-1)^{{1\over2}q(q-1)+\mu(\alpha^\minus_1)
            +\cdots+\mu(\alpha^\minus_{j-1})} \cr
  & \qquad\qquad\qquad\qquad
      \#\eusM^1(\alpha^\minus_j,\alpha_0)\,
\#\eufM^0_\sigma((\alpha^\minus_1,\ldots,\alpha^\minus_{j-1},
                  \alpha_0,\alpha^\minus_{j+1},
                  \ldots,\alpha^\minus_{k^\minus}),
                 (\alpha^\plus_1,\ldots,\alpha^\plus_{k^\plus})
                ) \cr
  & \qquad\qquad\qquad\qquad
         \alpha^\minus_1\otimes\cdots\otimes\alpha^\minus_{k^\minus}
         \otimes\alpha^\plus_1\otimes\cdots\otimes\alpha^\plus_{k^\plus} \cr
\noalign{\smallskip}
  & - \smash{\sum_{i=1}^{k^\plus}}\,
        \smash{\sum_{\alpha^\pm_*,\alpha_0\in\eusC^0_H}}
    (-1)^{{1\over2}q(q-1)+\mu(\alpha^\minus_1)
    +\cdots+\mu(\alpha^\minus_{k^\minus})
            +\mu(\alpha^\plus_1)+\cdots+\mu(\alpha^\plus_{i-1})+\mu(\alpha_0)
           } \cr
  & \qquad\qquad\qquad\qquad
      \#\eufM^0_\sigma((\alpha^\minus_1,\cdots,\alpha^\minus_{k^\minus}),
(\alpha^\plus_1,\cdots,\alpha^\plus_{i-1},\alpha_0,\alpha^\plus_{i+1},
                  \ldots,\alpha^\plus_{k^\plus})
                ) \,
      \#\eusM^1(\alpha_0,\alpha^\plus_i,) \cr
  & \qquad\qquad\qquad\qquad
         \alpha^\minus_1\otimes\cdots\otimes\alpha^\minus_{k^\minus}
         \otimes\alpha^\plus_1\otimes\alpha^\plus_{k^\plus} \cr
}\hss$$$$\eqalign{
  & = (-1)^{{1\over2}(q-1)(q-2)}\,
      \smash{\sum_{\nu=1}^q (-1)^{\nu+1}}
      \smash{\sum_{\alpha^\pm_*\in\eusC^0_H}}
    \Big(
      \#\eufM^0_{\partial^1_\nu\sigma}
      ((\alpha^\minus_1
        ,\ldots
        ,\alpha^\minus_{k^\minus}
       )
       ,(\alpha^\plus_1
         ,\ldots
         ,\alpha^\plus_{k^\plus}
        )
      ) \cr
  & \qquad\qquad\qquad\qquad\qquad\qquad\qquad\qquad\quad
    - \#\eufM^0_{\partial^0_\nu\sigma}
      ((\alpha^\minus_1
        ,\ldots
        ,\alpha^\minus_{k^\minus}
       )
       ,(\alpha^\plus_1
         ,\ldots
         ,\alpha^\plus_{k^\plus}
        )
      )
    \Big)  \cr
  & \qquad\qquad\qquad\qquad\qquad\qquad\qquad\qquad
         \alpha^\minus_1\otimes\cdots\otimes\alpha^\minus_{k^\minus}
         \otimes\alpha^\plus_1\otimes\cdots\otimes\alpha^\plus_{k^\plus} \cr
  & = \sum_{\nu=1}^q
              (-1)^{\nu+1} ( Q\partial^1_\nu\sigma - Q\partial^0_\nu\sigma)
    = Qd\sigma .\cr
          }
$$
It follows that $Q$ is a chain map.
The formula for the degree of $Q$ follows from Prop.~3.3.3.
\enddemo

Recall that $S_{k^\minus} \times S_{k^\plus}$ acts
on $\eufJ^*_{g,k^\minus,k^\plus}(A)$ and on
$CF^*(M,A)^{\otimes k^\minus} \otimes CF_*(M,A)^{\otimes k^\plus}$.

\proclaim Theorem 4.3.3.

If $\rho \in S_{k^\minus} \times S_{k^\plus}$,
and $\sigma:[0,1]^q\to\eufJ_{g,k^\minus,k^\plus}(A)$
is semi-regular, then
$Q(\rho.\sigma) = \rho.Q(\sigma)$.
In particular, the diagram
$$\hss\comdia{
    (S_{k^\minus}\times S_{k^\plus})
    \times H_*(\eufJ_{g,k^\minus,k^\plus}(A))
  & \mapright{{\rm id}\times Q}
  & (S_{k^\minus}\times S_{k^\plus})
    \times H \bigl( CF^*(M,A)^{\otimes k^\minus}
                    \otimes CF_*(M,A)^{\otimes k^\plus}
             \bigr) \cr
  \mapdown{}
  &
  & \mapdown{} \cr
    H_*(\eufJ_{g,k^\minus,k^\plus}(A))
  & \mapright{Q}
  & H \bigl( CF^*(M,A)^{\otimes k^\minus}\otimes CF_*(M,A)^{\otimes k^\plus}
      \bigr) \cr
         }
  \hss
$$
commutes.

\demo Proof.

This is an immediate consequence of Proposition 3.6.2.
\enddemo

Next we state the gluing relations.

\proclaim Theorem 4.3.4.

If the simplices $\sigma_1:[0,1]^{q_1}\to\eufJ^*_{g_1,k^\minus_1,k^\plus_1}(A)$
and $\sigma_2:[0,1]^{q_2}\to\eufJ^*_{g_2,k^\minus_2,k^\plus_2}(A)$
are semi-regular,
$\sigma_1\gop^\ell_{ij}\sigma_2$ takes values in
$\eufJ^*_{g_1+g_2,k^\minus_1+k^\minus_2,k^\plus_1+k^\plus_2}(A)$
for all $\ell\ge0$,
and either all periodic orbits of $H$
are contractible, $g_1=0$, or $g_2=0$,
then
$Q(\sigma_1\gop^\ell_{ij}\sigma_2))
  = (Q\sigma_1) \gop_{ij} (Q\sigma_2)$
for all sufficiently large $\ell$,
In particular,
the diagram
$$\comdia{
    {\ninepoint\eqalign{
        & H_*(\eufJ^*_{g_1,k_1^\minus,k_1^\plus+1}(A)) \cr
        & \otimes H_*(\eufJ^*_{g_2,k_2^\minus+1,k_2^\plus}(A)) \cr
                       }
    }
  & \mapright{Q\otimes Q}
  & {\ninepoint\eqalign{
        & H\bigl(CF^*(M,A)^{\otimes k_1^\minus}
          \otimes
          CF_*(M,A)^{\otimes(k^\plus_1+1)}\bigr) \cr
        & \otimes H\bigl(CF^*(M,A)^{\otimes(k_2^\minus+1)}
          \otimes
          CF_*(M,A)^{\otimes k_2^\plus}\bigr) \cr
                       }
    } \cr
    \mapdown{\gop^{}_{ij}}
  &
  & \mapdown{\gop^{}_{ij}} \cr
    H_*(\eufJ^*_{g_1+g_2,k_1^\minus+k^\minus_2,k_1^\plus+k_2^\plus}(A))
  & \mapright{Q}
  & H\bigl(CF^*(M,A)^{\otimes(k_1^\minus+k_2^\minus)}
    \otimes
    CF_*(M,A)^{\otimes(k^\plus_1+k_2^\plus)}\bigr) \cr
          }
$$
commutes.

\demo Proof.

By Prop.~3.6.3,
$$\eqalign{
  & (Q\sigma_1) \gop_{ij} (Q\sigma_2) \cr
\noalign{\smallskip}
  & = (-1)^{{1\over2}q_1(q_1-1)}
      \smash{
        \sum_{\alpha^\pm_{1,*}\in\eusC^0_H}
            }
      \#\eufM^0_{\sigma_1}
      \big((\alpha^\minus_{1,1},\ldots,\alpha^\minus_{1,k^\minus_1+1}) ,
           (\alpha^\plus_{1,1},\ldots,\alpha^\plus_{1,k^\plus_1+1})
      \big) \cr
  & \qquad\qquad\qquad\qquad\qquad\qquad
     \alpha^\minus_{1,1}\otimes\cdots\otimes\alpha^\minus_{1,k^\minus_1+1}
     \otimes\alpha^\plus_{1,1}\otimes\cdots\otimes\alpha^\plus_{1,k^\plus_1+1}
\cr
  & \phantom{{}={}}
    \gop_{ij} (-1)^{{1\over2}q_2(q_2-1)}
       \smash{
          \sum_{\alpha^\pm_{2,*}\in\eusC^0_H}
             }
      \#\eufM^0_{\sigma_2}
      \big((\alpha^\minus_{2,1},\ldots,\alpha^\minus_{2,k^\minus_2+1}) ,
           (\alpha^\plus_{2,1},\ldots,\alpha^\plus_{2,k^\plus_2+1})
      \big) \cr
  & \qquad\qquad\qquad\qquad\qquad\qquad\qquad
     \alpha^\minus_{2,1}\otimes\cdots\otimes\alpha^\minus_{2,k^\minus_2+1}
     \otimes\alpha^\plus_{2,1}\otimes\cdots\otimes\alpha^\plus_{2,k^\plus_2+1}
\cr
\noalign{\smallskip}
  & = \smash{
        \sum_{\alpha^\pm_{1,*},\alpha^\pm_{2,*}\in\eusC^0_H}
            }
      \pm (-1)^{{1\over2}q_1(q_1-1)+{1\over2}q_2(q_2-1)}\,
      \langle \alpha^\plus_{1,i},\alpha^\minus_{2,j} \rangle \cr
  & \qquad\qquad\qquad\qquad
      \#\eufM^0_{\sigma_1}
      \big((\alpha^\minus_{1,1},\ldots,\alpha^\minus_{1,k^\minus_1+1}) ,
           (\alpha^\plus_{1,1},\ldots,\alpha^\plus_{1,k^\plus_1+1})
      \big) \cr
  & \qquad\qquad\qquad\qquad
       \#\eufM^0_{\sigma_2}
      \big((\alpha^\minus_{2,1},\ldots,\alpha^\minus_{2,k^\minus_2+1}) ,
           (\alpha^\plus_{2,1},\ldots,\alpha^\plus_{2,k^\plus_2+1})
      \big) \cr
  & \qquad\qquad\qquad\qquad
     \alpha^\minus_{2,1} \otimes \ldots \otimes \alpha^\minus_{2,j-1}
     \otimes \alpha^\minus_{1,1} \otimes \ldots \otimes
\alpha^\minus_{1,k^\minus_1}
     \otimes \alpha^\minus_{2,j+1} \otimes \ldots \otimes
\alpha^\minus_{2,k^\minus_2+1} \cr
  & \qquad\qquad\qquad\qquad
     \otimes \alpha^\plus_{1,1} \otimes \ldots \otimes \alpha^\plus_{1,i-1}
     \otimes \alpha^\plus_{2,1} \otimes \ldots \otimes
\alpha^\plus_{2,k^\plus_2}
     \otimes \alpha^\plus_{1,i+1} \otimes \ldots \otimes
\alpha^\plus_{1,k^\plus_1+1} \cr
\noalign{\smallskip}
  & =   \smash{
        \sum_{\ss\alpha^\pm_{1,*},\alpha^\pm_{2,*}\in\eusC^0_H
              \atop
              \ss {\rm omit}\>  \alpha^\plus_{1,i}, \alpha^\minus_{2,j}
             }
            }
     (-1)^{{1\over2}q(q-1)} \cr
  & \qquad\qquad\qquad\qquad
        \#\eufM^0_{\sigma_1\gop^\ell_{ij}\sigma_2}
               ((\alpha^\minus_{2,1},\ldots,\alpha^\minus_{2,j-1},
                 \alpha^\minus_{1,1},\ldots,\alpha^\minus_{1,k^\minus_1},
                 \alpha^\minus_{2,j+1},\ldots,\alpha^\minus_{2,k^\minus_2+1}
                 ) ,  \cr
  & \qquad\qquad\qquad\qquad
        \phantom{\#\eufM^0_{\sigma_1\gop^\ell_{ij}\sigma_2})}
                 (\alpha^\plus_{1,1},\ldots,\alpha^\plus_{1,i-1},
                  \alpha^\plus_{2,1},\ldots,\alpha^\plus_{2,k^\plus_2}
                  \alpha^\plus_{1,i+1},\ldots,\alpha^\plus_{1,k^\plus_1+1}
                 )
                ) \cr
  & \qquad\qquad\qquad\qquad
     \alpha^\minus_{2,1} \otimes \ldots \otimes \alpha^\minus_{2,j-1}
     \otimes \alpha^\minus_{1,1} \otimes \ldots \otimes
\alpha^\minus_{1,k^\minus_1}
     \otimes \alpha^\minus_{2,j+1} \otimes \ldots \otimes
\alpha^\minus_{2,k^\minus_2+1} \cr
  & \qquad\qquad\qquad\qquad
     \otimes \alpha^\plus_{1,1} \otimes \ldots \otimes \alpha^\plus_{1,i-1}
     \otimes \alpha^\plus_{2,1} \otimes \ldots \otimes
\alpha^\plus_{2,k^\plus_2}
     \otimes \alpha^\plus_{1,i+1} \otimes \ldots \otimes
\alpha^\plus_{1,k^\plus_1+1} \cr
\noalign{\smallskip}
  & = Q(\sigma_1\gop^\ell_{ij}\sigma_2) \cr
           }
$$
for sufficiently large $\ell$,
where the sign is the graded sign of the permutation \qci.
The Proposition follows.
\enddemo

\proclaim Theorem 4.3.5.

If $\sigma:[0,1]^q\to\eufJ^*_{g,k^\minus+1,k^\plus+1}(A)$ is semi-regular,
$\cop^\ell_{ij}$ takes values in $\eufJ^*_{g+1,k^\minus,k^\plus}(A)$
for all $\ell\ge0$,
and all periodic orbits of $H$ are contractible,
then $Q\cop^\ell_{ij}\sigma = \cop_{ij}Q\sigma$
for all sufficiently large $\ell$.
In particular,
the diagram
$$\comdia{
    H_*(\eufJ^*_{g,k^\minus+1,k^\plus+1}(A))
  & \mapright{Q}
  & H\bigl(CF^*(M,A)^{\otimes(k^\minus+1)}
    \otimes
    CF_*(M,A)^{\otimes(k^\plus+1)}\bigr) \cr
    \mapdown{\cop_{ij}}
  &
  & \mapdown{\cop_{ij}} \cr
    H_*(\eufJ^*_{g+1,k^\minus,k^\plus}(A))
  & \mapright{Q}
  & H\bigl(CF^*(M,A)^{\otimes k^\minus}
    \otimes
    CF_*(M,A)^{\otimes k^\plus }\bigr) \cr
          }
$$
commutes.

\demo Proof.

By Prop.~3.6.4,
$$\eqalign{
    \cop_{ij} Q\sigma
  & = (-1)^{{1\over2}q(q-1)}\,
      \cop_{ij}
      \smash{\sum_{\alpha^\pm_*\in\eusC^0_H}}
      \#\eufM^0_\sigma
      \big((\alpha^\minus_1,\ldots,\alpha^\minus_{k^\minus+1}) ,
           (\alpha^\plus_1,\ldots,\alpha^\plus_{k^\plus+1})
      \big) \cr
  & \qquad\qquad\qquad\qquad\qquad\qquad
     \alpha^\minus_1\otimes\cdots\otimes\alpha^\minus_{k^\minus+1}
     \otimes\alpha^\plus_1\otimes\cdots\otimes\alpha^\plus_{k^\plus+1} \cr
  \noalign{\smallskip}
  &  = (-1)^{{1\over2}q(q-1)}
       \smash{\sum_{\alpha^\pm_*\in\eusC^0_H}}
     \pm\langle \alpha^\plus_i,\alpha^\minus_j \rangle
      \#\eufM^0_\sigma
      \big((\alpha^\minus_1,\ldots,\alpha^\minus_{k^\minus+1}) ,
           (\alpha^\plus_1,\ldots,\alpha^\plus_{k^\plus+1})
      \big) \cr
  &  \qquad\qquad\qquad\qquad\quad\qquad
     \alpha^\minus_1\otimes\cdots\alpha^\minus_{j-1}
     \otimes\alpha^\minus_{j+1}\otimes\cdots\otimes\alpha^\minus_{k^\minus+1}
\cr
  &  \qquad\qquad\qquad\qquad\quad\qquad\phantom{(}
     \otimes\alpha^\plus_1\otimes\cdots\otimes\alpha^\plus_{i-1}
     \otimes\alpha^\plus_{i-1}\otimes\cdots\otimes\alpha^\plus_{k^\plus+1} \cr
  \noalign{\smallskip}
  &  = (-1)^{{1\over2}q(q-1)}
       \smash{\sum_{
           \ss\alpha^\pm_*\in\eusC^0_H
           \atop
           \ss{\rm omit}\>\alpha^\minus_j,\alpha^\plus_i
                   }
             }
       \#\eufM^0_{\cop^\ell_{ij}\sigma}
       \big((\alpha^\minus_1,\ldots,\alpha^\minus_{j-1},
             \alpha^\minus_{j+1},\ldots, \alpha^\minus_{k^\minus+1}
            ), \cr
  & \qquad\qquad\qquad\qquad\qquad\qquad\qquad\quad
           (\alpha^\plus_1,\ldots,\alpha^\plus_{i-1},
             \alpha^\plus_{i+1},\ldots, \alpha^\plus_{k^\plus+1}
            )
       \big) \cr
  & \qquad\qquad\qquad\qquad\qquad\qquad
     \alpha^\minus_1\otimes\cdots\alpha^\minus_{j-1}
     \otimes\alpha^\minus_{j+1}\otimes\cdots\otimes\alpha^\minus_{k^\minus+1}
\cr
  & \qquad\qquad\qquad\qquad\qquad\qquad
     \phantom{(}
     \otimes\alpha^\plus_1\otimes\cdots\otimes\alpha^\plus_{i-1}
     \otimes\alpha^\plus_{i-1}\otimes\cdots\otimes\alpha^\plus_{k^\plus+1} \cr
  \noalign{\smallskip}
  & = Q \cop^\ell_{ij}\sigma \cr
       }
$$
for sufficiently large $\ell$,
where the sign is the graded sign of the permutation \qcj.
The Proposition follows.
\enddemo

Every compact symplectic manifold admits a regular Hamiltonian
for which all periodic orbits are contractible, and for which
Thm.~4.3.4 and Thm.~4.3.5 therefore apply.
In \S4.5 we
prove that the map $Q$ is essentially independent of $A$,
and thus that the diagrams in Thm.~4.3.4 and
Thm.~4.3.5 will commute for any $A$.

The rays $X^\pm_i$ are needed to define the gluing maps.
However,
if $H$ is time-independent,
then they play no role in the construction of
the moduli spaces $\eufM^d_\sigma$.
Thus there is a chain map
$Q^0: C_*(\eufJ^{0*}_{g,k^\minus,k^\plus}(A))
    \to
    CF^*(M,A)^{\otimes k^\minus} \otimes CF_*(M,A)^{\otimes k^\plus}
$
defined the same way as $Q$.

There is a natural projection map
$\eufJ_{g,k^\minus,k^\plus}(A)\to\eufJ^0_{g,k^\minus,k^\plus}(A)$.
The induced map
$\eufJ^*_{g,k^\minus,k^\plus}(A)\to\eufJ^{0*}_{g,k^\minus,k^\plus}(A)$
is in general only defined on a large subset of
$\eufJ^*_{g,k^\minus,k^\plus}(A)$.
A simple perturbation argument as in \S4.1 still gives a homomorphism
$H_*(\eufJ^*_{g,k^\minus,k^\plus}(A))\to
H_*(\eufJ^{0*}_{g,k^\minus,k^\plus}(A))$.
This proves the following Theorem.

\proclaim Theorem 4.3.6.

If the Hamiltonian $H$ is time-independent,
then $Q$ factors through a homomorphism
$$Q^0: H_*(\eufJ^{0*}_{g,k^\minus,k^\plus}(A))
     \to
     H(CF^*(M,A)^{\otimes k^\minus} \otimes CF_*(M,A)^{\otimes k^\plus}).
$$

\subheading 4.4. Products and relations.

In this section we show how the the homomorphisms $Q$ give products
and other operations in Floer (co)homology,
and how Thm.~4.3.4 and 4.3.5 give relations for these products.
The definitions and proofs are illustrated in fig.~2 and 3.

Let $\Theta_{g,k^\minus,k^\plus}$ denote the canonical generator of
$H_0(\eufJ^*_{g,k^\minus,k^\plus}(A))$.

\proclaim Lemma 4.4.1.

The elements $Q(\Theta_{g,k^\minus,k^\plus})$  have degree $2n(1-g-k^\minus)$
modulo $2N_0$.

\demo Proof.

For $g=0$, this follows directly from Prop.~4.3.2.
For $g\ge1$ Prop.~4.3.2 only gives the degree modulo $2N_1$.
If all perirodic orbits of $H$ are contractible,
then it follows from Thm.~4.3.5 that
$$Q(\Theta_{g,k^\minus,k^\plus})
  =  Q(\cop_{11}^g \Theta_{0,k^\minus+g,k^\plus+g})
  =  \cop_{11}^g Q(\Theta_{0,k^\minus+g,k^\plus+g}) .
$$
By Prop.~4.3.2 the element $Q(\Theta_{0,k^\minus+g,k^\plus+g})$
has degree $2n(1-g-k^\minus)$ modulo $2N_0$.
The chain map $\cop_{11}$ has degree $0$ modulo $2N_0$.
Hence $Q(\Theta_{g,k^\minus,k^\plus})$ has degree $2n(1-g-k^\minus)$
modulo $2N_0$.
By the results of \S4.5, this holds even if $H$ has non-contractible
periodic orbits.
\enddemo

To simplify the notation,
for
$x_1\in H(CF_*(M,A)^{\otimes k^\minus_1}\otimes CF^*(M,A)^{\otimes k^\plus_1})$
and
$x_2\in H(CF_*(M,A)^{\otimes k^\minus_2}\otimes CF^*(M,A)^{\otimes k^\plus_2})$
we write
$x_1 \gop x_2$
for
$x_1 \gop_{k^\plus_1\,1} x_2$.
In this case the sign in the definition of $\gop$ is positive.

\subsubheading 4.4.1. The identity map.

The class $\Theta_{0,1,1}(A,A)$ can be represented by
the simplex $\sigma_0$ as in Remark 2.3.7.
By Prop.~3.6.1,
$Q(\sigma_0) = \sum_{\alpha\in\eusC^0} \alpha\otimes\alpha$.
In other words,
the class $Q(\Theta_{0,1,1})$ can be represented by the Floer cycle
$$\sum_{\alpha\in\eusC^0} \alpha\otimes\alpha .$$
It follows that for $a\in HF^*(M,A)$,
$$a = Q(\Theta_{0,1,1}) \gop a.$$
Similarly, for $x\in HF^*(M,A)$,
$$x = x \gop Q(\Theta_{0,1,1}).$$
Thus the class $Q(\Theta_{0,1,1})$ induces the identity maps
$HF^*(M,A)\to HF^*(M,A)$ and $HF_*(M,A)\to HF_*(M,A)$.
Other homomorphisms are given by the classes $Q(\Theta_{g,1,1})$.
By Lemma 4.4.1, these have degree $2ng$ and $-2ng$
respectively.

\subsubheading 4.4.2. The symplectic cup product.

We define the symplectic cup product of $a,b\in HF^*(M,A)$ as
$$a \cup b = (Q(\Theta_{0,1,2}) \gop a) \gop b.$$
By Lemma 4.4.1 the cup product has degree $0$.

It follows from Thm.~4.3.4 that the symplectic cup product is associative,
$$\eqalign{
    (a \cup b) \cup c
  & = \bigl( Q(\Theta_{0,1,2})
           \gop ( ( Q(\Theta_{0,1,2}) \gop a ) \gop b )
    \bigr) \gop c \cr
  & = \bigl( ( ( Q(\Theta_{0,1,2}) \gop_{21} Q(\Theta_{0,1,2}) )
               \gop a
             ) \gop b
      \bigr) \gop c \cr
  & = \bigl( ( Q(\Theta_{0,1,3}) \gop a ) \gop b \bigr) \gop c \cr
  & = \bigl( ( ( Q(\Theta_{0,1,2}) \gop_{11} Q(\Theta_{0,1,2}) )
               \gop a
             ) \gop b
      \bigr) \gop c \cr
  & = ( Q(\Theta_{0,1,2}) \gop a )
    \gop \bigl( ( Q(\Theta_{0,1,2}) \gop b ) \gop c \bigr) \cr
  & = a \cup ( b \cup c ).\cr
          }
$$
Even though we have permuted $a$ and $Q(\Theta_{0,1,2})$,
the signs are positive,
for $Q(\Theta_{0,1,2})$ has degree 0.

It follows from Thm.~4.3.4 and 4.3.3
that the cup product is graded commutative.
In fact,
let $\rho_0$ denote the permutation of order two in the symmetric group
$S_2$ on two letters.
Then
$$\eqalign{
  a \cup b
  & = ( Q(\Theta_{0,1,2}) \gop a ) \gop b \cr
  & = (-1)^{ij} ( (1,\rho_0).Q(\Theta_{0,1,2}) \gop b ) \gop a \cr
  & = (-1)^{ij} ( Q((1,\rho_0).\Theta_{0,1,2}) \gop b ) \gop a \cr
  & = (-1)^{ij} ( Q(\Theta_{0,1,2}) \gop b ) \gop a \cr
  & = (-1)^{ij} b \cup a . \cr
          }
$$
for $a\in HF^i(M,A)$ and $b\in HF^j(M,A)$.
Finally, define the unit class 1 as
$$1 = Q(\Theta_{0,1,0}) \in HF^0(M).$$
It follows from Prop 4.3.4 and \S4.4.1 that 1 is an identity
element for the symplectic cup product,
$$1\cup a
  = ( Q(\Theta_{0,1,2}) \gop  Q(\Theta_{0,1,0}) ) \gop a
  = Q(\Theta_{0,1,1}) \gop a
  = a .
$$

The classes $Q(\Theta_{g,1,2})$ give higher genus cup products.
However, by Prop.~3.4.3,
$Q(\Theta_{g,1,2}) = Q(\Theta_{g,1,1}) \gop Q(\Theta_{0,1,2})$.
Hence the higher genus cup products are simply
the composition of the cup product and the homomorphisms
$HF^*(M,A)\to HF^*(M,A)$
given by $Q(\Theta_{g,1,1})$.

\subsubheading 4.4.3. The symplectic intersection product.

We define the symplectic intersection  product of $x,y\in HF_*(M,A)$ as
$$x \cdot y = x \gop (y \gop Q(\Theta_{0,2,1})).$$
The same argument as in \S4.4.2 shows that the intersection product
is graded commutative,
associative, and has degree $-2n$.
The symplectic top class $[M]_\omega=Q(\Theta_{0,0,1})\in HF_{2n}(M,A)$
is an indentity element for the symplectic intersection product.

\subsubheading 4.4.4. The symplectic cap product.

We define the symplectic cap product of $x \in HF_*(M)$ and $a \in HF^*(M)$ as
$$x \cap a = x \gop Q(\Theta_{0,1,2}) \gop a.$$
Then
$$\eqalign{
   (x \cap a) \cap b
& = (x \gop Q(\Theta_{0,1,2}) \gop a)
    \gop Q(\Theta_{0,1,2})) \gop b \cr
& = x \gop \bigl( (Q(\Theta_{0,1,2}) \gop_{11} Q(\Theta_{0,1,2}))
                  \gop a
           \bigr) \gop b \cr
& = x \gop (Q(\Theta_{0,1,3}) \gop a ) \gop b \cr
& = x \gop \bigl((Q(\Theta_{0,1,2}) \gop_{21} Q(\Theta_{0,1,2}))
    \gop a \bigr)\gop b \cr
& = x \gop Q(\Theta_{0,1,2})
      \gop \bigl( (Q(\Theta_{0,1,2}) \gop a) \gop b \bigr) \cr
& = x \cap (a \cup b). \cr
           }$$

\subsubheading 4.4.5.  Symplectic Poincar\'e duality.

We define the symplectic Poincar\'e duals of $x\in HF_*(M,A)$
and $a\in HF^*(M,A)$ as
$$\eqalign{
    x^\sharp & = x \gop Q(\Theta_{0,2,0}) \cr
    a^\flat  & = Q(\Theta_{0,0,2}) \gop a \cr
          }
$$
By Thm.~4.3.4 and \S4.4.1, these maps are inverses,
$$\eqalign{
  (x^\sharp)^\flat
  & = Q(\Theta_{0,0,2}) \gop ( x \gop Q(\Theta_{0,2,0}) ) \cr
  & =  x \gop ( Q( \Theta_{0,0,2}) \gop_{22} Q(\Theta_{0,2,0}) ) \cr
  & =  x \gop  Q(\Theta_{0,1,1}) = x , \cr
          }
$$
and similarly
$$(a^\flat)^\sharp = a .$$
It also follows from Thm.~4.3.4 that
$$\eqalign{
   (a\cup b)^\flat & = a^\flat \cdot b^\flat \cr
   x\cap a & = x \cdot a^\flat \cr
   1^\flat & = [M]_\omega .\cr
          }
$$

\subsubheading 4.4.6. The Euler characteristic and the number of tori.

The Floer cohomology groups are the ordinary cohomology groups
with the grading reduced modulo $2N_0$;
see [F3] and [F4] Thm.~5.
In particular, the Euler characteristic of the Floer
cohomology is the same as the ordinary Euler characteristic $\chi(M)$.
Choose $A$ such that all periodic orbits for $H$ are contractible.
It then follows from Thm.~4.3.5 and \S4.4.1 that
the number of perturbed pseudo-holomorphic tori
with a given conformal structure is
$$Q(\Theta_{1,0,0}) = \cop_{11} Q(\Theta_{0,1,1})
= \cop_{11} \sum_{\alpha \in \eusC^0_H} \alpha\otimes\alpha
= \sum_{\alpha \in \eusC^0_H} (-1)^{\mu(\alpha)}
= \chi(M).$$

\subsubheading 4.4.7. Twists.

If the Hamiltonian $H$ is time-independent,
then we get a loop in $\eusJhat^{+*}_{g,k^\minus,k^\plus}(A)$
by a full counterclockwise twist of the ray $X^\pm_i$.
The image of this loop in $\eufJ^*_{g,k^\minus,k^\plus}(A)$
defines an element $\Xi^\pm_i\in H_1(\eufJ^*_{g,k^\minus,k^\plus}(A))$.
This class exists even if $H$ is time-dependent.
If $H$ is time-independent,
then the image of $\Xi^\pm_i$ in $H_1(\eufJ^{0*}_{g,k^\minus,k^\plus}(A))$
is zero.
It follows from Thm.~4.3.6 that if $H$ is time-independent,
then $Q(\Xi^\pm_i)=0$.
The class $\cop_{ij}\Xi^+_i=-\cop_{ij}\Xi^-_j$ corresponds to a Dehn twist
of $\Sigmabar_{g+1,k^\minus-1,k^\plus-1}$.
It follows from Thm.~4.3.5 that $Q(\cop_{ij}\Xi^\plus_i)=0$ as well.
By the results of \S4.5, these identities also hold for time-dependent
Hamiltonians.

\subsubheading 4.4.8. The symplectic Massey product.

The associativity relation $a\cup(b\cup c)=(a\cup b)\cup c$
gives rise to a secondary operation, the Massey product.
The Massey product $M(a,b,c)$ is an element of
$HF^*(M,A)/(a\cup HF^*(M,A) + HF^*(M,A)\cup c)$,
which is defined for $a,b,c\in HF^*(M,A)$
with $a\cup b= b\cup c=0$.

Let $\theta_{0,1,2}\in C_0(\eufJ_{0,1,2}(A))$ be a cycle
that represents the class $\Theta_{0,1,2}$.
Now
$\Theta_{0,1,2}\gop_{11}\Theta_{0,1,2}
 = \Theta_{0,1,2}\gop_{21}\Theta_{0,1,2}.
$
Hence
$\theta_{0,1,2}\gop^\ell_{11}\theta_{0,1,2}
  - \theta_{0,1,2}\gop^\ell_{21}\theta_{0,1,2}
  =d\lambda$
for some $\lambda\in C_1(\eufJ_{0,1,3}(A))$.
Let $\alpha$, $\beta$, $\gamma$
be Floer cocycles that represent $a$, $b$ and $c$.
Then there exist cochains $\zeta$ and $\xi$ such that
$\alpha\cup\beta = d\zeta$ and $\beta\cup\gamma=d\xi$.
Then the Massey product $M(a,b,c)$ is represented by the cochain
$$\bigl((Q(\lambda) \gop \alpha) \gop \beta \bigr)\gop \gamma
  + (Q(\theta_{0,1,2}) \gop \zeta) \gop \gamma
   - (-1)^{\mu(\alpha)} (Q(\theta_{0,1,2}) \gop \alpha) \gop \xi .
$$
A similar construction is discussed in [Fu].

\subsubheading 4.4.9 Other operations.

Using analogous constructions,
the higher dimensional homology classes in
$H_*(\eufJ^*_{g,k^\minus,k^\plus}(A))$ also give products.
Thms.~4.3.4 and 4.3.5 give relations between these products,
and corresponding secondary operations can be defined.

\subheading 4.5. Independence of the choice of Hamiltonian,
almost complex structure, and coherent system of orientations.

In this section we put Floer's proof, [F4] Thm.~4,
that the Floer (co)homology groups $HF_*(M,A)$ and $HF^*(M,A)$
are independent of the choice of $A=(H,\bfJ,\eufo)$
into our framework.
We then show that the chain map $Q$ is also independent of
the choice of $A$.

That $HF_*(M,A)$ is independent of $A$
means that there exists a canonical system of homomorphisms
$L_{AA'}:HF_*(M,A)\to HF_*(M,A')$ which has the following two properties:

\vskip\parskip
\halign{\hskip\parindent\sl#&\qquad#\cr
  Reflexivity:
  & $L_{AA}$ is the identity\cr
\noalign{\smallskip}
  Transitivity:
  & $L_{A'A''}\circ L_{AA'} = L_{AA''}$ \cr
        }

To define the chain maps $L_{AA'}$ we need a generalization of the
homomorphism
$Q:H_*(\eufJ^*_{g,k^\minus,k^\plus}(A))
 \to H(CF^*(M,A)^{\otimes k^\minus}\otimes CF_*(M,A)^{\otimes k^\plus})
$
of \S4.3.
Let
$$\eufJ^*_{g,k^\minus,k^\plus}
((A^\minus_1,\ldots,A^\minus_{k^\minus}),(A^\plus_1,\ldots,A^\plus_{k^\plus})),
$$
with $A^\pm_i=(H^\pm_i,\bfJ^\pm_i,\eufo^\plusminus_i)$,
be a space defined the same way as $\eufJ^*_{g,k^\minus,k^\plus}(A)$ in \S2.2,
but we now let $\bfJ$ be a family of conformal structures on $M$
parametrized by $\Sigmabar_{g,k^\minus,k^\plus}$.
We require that
$\bfJ=\bfJ^\pm_i$
and
$R=(dt^\plusminus_i-id\theta^\plusminus_i)\otimes\nabla H^\pm_i$
on $\Delta^\pm_i$.
There then exists a gluing map $\gop^\ell_{ij}$ as in \S4.1,
which is defined if $A^\plus_{1,i}=A^\minus_{2,j}$.
Similarly there exists a gluing map $\cop^\ell_{ij}$
which is defined if $A^\minus_j=A^\plus_i$.

If $A^\plus_{1,i}=A^\minus_{2,j}$,
there is a chain map
$$\eqalign{
    \gop_{ij}:{}
  & \left(
     \bigotimes_{\nu=1}^{k^\minus_1} CF^*(M,A^\minus_{1,\nu})
     \otimes \bigotimes_{\nu=1}^{k^\plus_1} CF_*(M,A^\plus_{1,\nu})
     \right) \cr
  &  \qquad\qquad\qquad\qquad
     \otimes
     \left(
     \bigotimes_{\nu=1}^{k^\minus_2} CF^*(M,A^\minus_{2,\nu})
     \otimes \bigotimes_{\nu=1}^{k^\plus_2} CF_*(M,A^\plus_{2,\nu})
     \right) \cr
  &  \to \bigotimes_{\nu=1}^{j-1} CF^*(M,A^\minus_{2,\nu})
         \otimes \bigotimes_{\nu=1}^{k^\minus_1} CF^*(M,A^\minus_{1,\nu})
         \otimes \bigotimes_{\nu=j+1}^{k^\minus_2} CF^*(M,A^\minus_{2,\nu}) \cr
  & \qquad
         \otimes \bigotimes_{\nu=1}^{i-1} CF_*(M,A^\plus_{1,\nu})
         \otimes \bigotimes_{\nu=1}^{k^\plus_2} CF_*(M,A^\plus_{2,\nu})
         \otimes \bigotimes_{\nu=i+1}^{k^\plus_1} CF_*(M,A^\plus_{1,\nu}) \cr
          }
$$
defined the same way as the chain map $\gop_{ij}$ of Def.~4.2.4.
Similarly, if $A^\minus_j=A^\plus_i$,
there is a chain map
$$\cop_{ij}:{}
  \bigotimes_{\nu=1}^{k^\minus} CF^*(M,A^\minus_\nu)
     \otimes \bigotimes_{\nu=1}^{k^\plus} CF_*(M,A^\plus\nu)
    \to
     \bigotimes_{\ss\nu=1\atop\ss\nu\ne j}^{k^\minus} CF^*(M,A^\minus_\nu)
     \otimes \bigotimes_{\ss\nu=1\atop\ss\nu\ne i}^{k^\plus}
CF_*(M,A^\plus_\nu) .
$$

It is straightforward to generalize the construction of $Q$
to get a chain map
$$\eqalign{
  & Q:C^\pi_*(\eufJ^*_{g,k^\minus,k^\plus}
((A^\minus_1,\ldots,A^\minus_{k^\minus}),(A^\plus_1,\ldots,A^\plus_{k^\plus}))
         ) \cr
  & \qquad \to
     \bigotimes_{\nu=1}^{k^\minus} CF^*(M,A^\minus_\nu)
     \otimes
     \bigotimes_{\nu=1}^{k^\plus} CF_*(M,A^\plus_\nu) . \cr
          }
$$
The main modification is that we replace $\bfJ(u(x,y))$
in \qbb{} by $\bfJ((x,y),u(x,y))$.
This chain map satisfies the appropriate
analogues of the theorems of \S4.3.

Let $\Theta_{0,1,1}(A^\minus,A^\plus)$ be the canonical generator
of $H_0(\eufJ^*_{0,1,1}(A^\minus,A^\plus))$.
We define the homomorphism
$$L_{AA'} : HF_*(M,A) \to HF_*(M,A')$$
as
$$x \mapsto x\gop Q(\Theta_{0,1,1}(A,A')).$$
Reflexivity follows from \S4.4.1.
By Thm.~4.3.4,
$$\eqalign{
  L_{A'A''}L_{AA'}x
  & = ( x \gop Q(\Theta_{0,1,1}(A,A')) )
      \gop Q(\Theta_{0,1,1}(A',A'')) \cr
  & = x
      \gop
      ( Q(\Theta_{0,1,1}(A,A')) \gop_{11} Q(\Theta_{0,1,1}(A',A'')) ) \cr
  & = x
      \gop
       Q(\Theta_{0,1,1}(A,A')\gop_{11}\Theta_{0,1,1}(A',A'')) \cr
  & = x \gop  Q(\Theta_{0,1,1}(A,A''))
    = L_{AA''} x  \cr
          }
$$
and we have established transitivity.
Therefore we can identify the groups $HF_*(M,A)$ for different
$A$ using the maps $L_{AA'}$.
This gives the abstract Floer homology groups $HF_*(M)$.

Similarly we define a map $HF^*(M,A)\to HF^*(M,A')$
by
$$L_{A'A}^* x  = Q(\Theta_{0,1,1}(A',A)) \gop x.$$
More generally we define homomorphisms
$$\hss
  H\left(
  \bigotimes_{\nu=1}^{k^\minus} CF^*(M,A^\minus_\nu)
  \otimes
  \bigotimes_{\nu=1}^{k^\plus} CF_*(M,A^\plus_\nu)
  \right)
  \to
  H\left(
  \bigotimes_{\nu=1}^{k^\minus} CF^*(M,A^{-\prime}_\nu)
  \otimes
  \bigotimes_{\nu=1}^{k^\plus} CF_*(M,A^{+\prime}_\nu)
  \right)
  \hss
$$
by
$$\eqalign{
    x
    \mapsto {}
  & Q(\Theta_{0,1,1}(A^{-\prime}_{k^\minus},A^\minus_{k^\minus}))
    \gop_{1k^\minus}
    \cdots \gop_{12}
    Q(\Theta_{0,1,1}(A^{-\prime}_1,A^\minus_1))\gop_{11} x \cr
  & \qquad
    \gop_{11}
    Q(\Theta_{0,1,1}(A^{+\prime}_1,A^\plus_1))
    \gop_{21}
    \cdots \gop_{k^\plus1}
    Q(\Theta_{0,1,1}(A^{+\prime}_{k^\plus},A^\plus_{k^\plus})).\cr
          }
$$
These homomorphisms are also reflexive and transitive.

The spaces
$\eufJ^*_{g,k^\minus,k^\plus}
      ((A^\minus_1,\ldots,A^\minus_{k^\minus})
       ,(A^\plus_1,\ldots,A^\plus_{k^\minus})
      )
$
with different $A^\pm_i$ are canonically homotopy equivalent.
Thus we can simply write $H_*(\eufJ^*_{g,k^\minus,k^\plus})$.
The canonical homotopy equivalence
$$\eufJ^*_{g,k^\minus,k^\plus}
      ((A^\minus_1,\ldots,A^\minus_{k^\minus})
       ,(A^\plus_1,\ldots,A^\plus_{k^\minus})
      )
  \to
  \eufJ^*_{g,k^\minus,k^\plus}
      ((A^{-\prime}_1,\ldots,A^{-\prime}_{k^\minus})
       ,(A^{+\prime}_1,\ldots,A^{+\prime}_{k^\minus})
      )
$$
can be realized explicitly by the map
$$\eqalign{
  c
    \mapsto {}
  & c_0(A^{-\prime}_{k^\minus},A^\minus_{k^\minus})
          \gop^\ell_{1k^\minus} \cdots
          \gop^\ell_{12} c_0(A^{-\prime}_1,A^\minus_1)
          \gop^\ell_{11} c \cr
  & \qquad
          \gop^\ell_{11} c_0(A^{+\prime}_1,A^\plus_1)
          \gop^\ell_{21} \cdots
          \gop^\ell_{k^\plus1}
          c_0(A^{+\prime}_{k^\plus},A^\plus_{k^\plus}) \cr
            }
$$
where $c_0(A,A')$ is any element of $\eufJ_{0,1,1}(A,A')$ and $\ell\ge0$.
It follows that the induced isomorphisms
$$\eqalign{
  & H_*(\eufJ^*_{g,k^\minus,k^\plus}
      ((A^\minus_1,\ldots,A^\minus_{k^\minus})
       ,(A^\plus_1,\ldots,A^\plus_{k^\minus})
      )
     ) \cr
  & \qquad  \to
  H_*(\eufJ^*_{g,k^\minus,k^\plus}
      ((A^{-\prime}_1,\ldots,A^{-\prime}_{k^\minus})
       ,(A^{+\prime}_1,\ldots,A^{+\prime}_{k^\minus})
      )
     ) \cr
        }
$$
are given by
$$\eqalign{
  x
    \mapsto {}
  & \Theta_{0,1,1}(A^{-\prime}_{k^\minus},A^\minus_{k^\minus})
    \gop_{1k^\minus}
    \cdots
          \gop_{12} \Theta_{0,1,1}(A^{-\prime}_1,A^\minus_1)
          \gop_{11} x \cr
  & \qquad
          \gop_{11} \Theta_{0,1,1}(A^{+\prime}_1,A^\plus_1)
          \gop_{21} \cdots
    \gop_{k^\plus1} \Theta_{0,1,1}(A^{+\prime}_{k^\plus},A^\plus_{k^\plus}).
    \cr
            }
$$
It follows, by repeated use of Thm.~4.3.4 that the diagram
$$\hss
  \comdia{
    H_*(\eufJ^*_{g,k^\minus,k^\plus}
        ((A^\minus_1,\ldots,A^\minus_{k^\minus})
         ,(A^\plus_1,\ldots,A^\plus_{k^\minus})
        )
       )
  & \mapright{Q}
  & \bigotimes_{j=1}^{k^\minus} CF^*(M,A^\minus_j)
    \otimes
    \bigotimes_{i=1}^{k^\plus} CF_*(M,A^\plus_i) \cr
  \mapdown{}
  &
  & \mapdown{} \cr
  H_*(\eufJ^*_{g,k^\minus,k^\plus}
      ((A^{-\prime}_1,\ldots,A^{-\prime}_{k^\minus})
       ,(A^{+\prime}_1,\ldots,A^{+\prime}_{k^\minus})
      )
     )
  & \mapright{Q}
  & \bigotimes_{j=1}^{k^\minus} CF^*(M,A^{-\prime}_j)
    \otimes
    \bigotimes_{i=1}^{k^\plus} CF_*(M,A^{+\prime}_i) \cr
         }
  \hss
$$
commutes.
Thus we can view $Q$ as a homomorphism
$$Q:H_*(\eufJ^*_{g,k^\minus,k^\plus})
    \to H(CF^*(M)^{\otimes k^\minus}
        \otimes CF_*(M)^{\otimes k^\plus}).
$$

%%%%%%%%%%%%%%%%%%%%%%%%%%%%%%%%%%%%%%%%%%%%%%%%%%%%%%%%%%%%%

\references

\ref[D]
{\sl S.K. Donaldson,}
Orientation of Yang-Mills moduli spaces,
J. Diff. Geom. {\bf26} (1987), 397--428.

\ref[DK]
{\sl S.K. Donaldson and P.B. Kronheimer,}
``The geometry of four-manifolds,''
Oxford University Pr., Oxford, 1990.

\ref[EE]
{\sl C. J. Earle and J. Eells,}
A fiber bundle description of Teichm\"uller theory,
J. Diff. Geom. {\bf3} (1969), 19--43.

\ref[F1]
{\sl A. Floer,}
A relative Morse index for the symplectic action,
Commun. Pure Appl. Math. {\bf41} (1988), 393--407.

\ref[F2]
{\sl A. Floer,}
The unregularized gradient flow of the symplectic action,
Commun. Pure Appl. Math. {\bf41} (1988), 775--813.

\ref[F3]
{\sl A. Floer,}
Witten's complex and infinite dimensional Morse theory,
J. Diff. Geom. {\bf30} (1989), 207--221.

\ref[F4]
{\sl A. Floer,}
Symplectic fixed points and holomorphic spheres,
Commun. Math. Phys. {\bf120} (1989), 575--611.

\ref [FH]
{\sl A. Floer and H. Hofer,}
Coherent orientations for periodic orbit problems
in symplectic geometry,
Z. Math. {\bf212} (1993), 13--38.

\ref [FHS]
{\sl A. Floer, H. Hofer, and D. Salamon,}
Transversality in elliptic Morse theory for the symplectic action functional,
preprint, 1994.

\ref [Fu]
{\sl K. Fukaya,}
Morse homotopy and its quantization,
preprint, 1994.

\ref[HW]
{\sl P.J. Hilton and S. Wylie,}
``Homology theory,''
Cambridge Univ. Pr., Cambridge, 1960.

\ref[HS]
{\sl H. Hofer and D. Salamon,}
Floer homology and Novikov rings,
in ``A. Floer memorial volume,''
ed. H. Hofer, C. Taubes, A. Weinstein and E. Zehnder,
Birkh\"auser, Basel, to appear.

\ref[HZ]
{\sl H. Hofer and E. Zehnder,}
``Symplectic invariants and hamiltonian dynamics,''
Birkh\"auser, Basel, 1994.

\ref[J]
{\sl J. Jost,}
``Two-dimensional geometric variational problems,''
John Wiley, New York, 1991.

\ref[KM]
{\sl M. Kontsevich and Yu. Manin,}
Gromov-Witten classes, quantum cohomology, and enumerative geometry,
preprint, 1993. {\tt hep-th/9402177}

\ref[M]
{\sl D. McDuff,}
Elliptic methods in symplectic geometry,
Bull. Amer. Math. Soc. {\bf 23} (1990), 311--358.

\ref[MS1]
{\sl D. McDuff and D. Salamon,}
``$J$-holomorphic curves and quantum cohomology,''
AMS, Providence RI, 1994.

\ref[MS2]
{\sl D. McDuff and D. Salamon,}
``Introduction to symplectic topology,''
Oxford University Press, Oxford, to appear.

\ref[PW]
{\sl T. Parker and J. G. Wolfson,}
Pseudoholomorphic maps and bubble trees,
J. Geom. Anal. {\bf3} (1993), 63--98.

\ref[P]
{\sl S. Piunikhin,}
Quantum and Floer cohomology have the same ring structure,
preprint, 1994.

\ref[R]
{\sl J.~R\aa de,}
Remarks on a paper by A. Floer,
preprint, 1994.

\ref[Ru]
{\sl Y. Ruan,}
Topological sigma model and Donaldson type invariants in Gromov theory,
preprint.

\ref[RT1]
{\sl Y. Ruan and G. Tian,}
A mathematical theory of quantum cohomology,
preprint, 1994.

\ref[RT2]
{\sl Y. Ruan and G. Tian,}
Bott-type symplectic Floer cohomology and its multiplication structures,
in preparation.

\ref[S]
{\sl D. Salamon,}
Morse theory, the Conley index and Floer homology,
Bull. London Math. Soc. {\bf22} (1990), 113--140.

\ref[SZ]
{\sl D. Salamon and E. Zehnder,}
Morse theory for periodic solutions of Hamiltonian systems
and the Maslov index,
Commun. Pure Appl. Math. {\bf45} (1992), 1301--1360.

\ref[Sch1]
{\sl M. Schwarz,}
``Morse homology,''
Birkh\"auser, Basel, 1993.

\ref[Sch2]
{\sl M. Schwarz,}
Ph. D. thesis,
ETH, Z\"urich, in preparation.

\ref[T]
{\sl A. J. Tromba,}
``Teichm\"uller theory in Riemannian geometry,''
Birkh\"auser, Basel, 1992.

\ref[V]
{\sl C. Vafa,}
Topology mirrors and quantum rings,
in ``Essays on mirror manifolds'', ed. S.-T. Yau,
International Press, Hong Kong, 1992.

\ref[W]
{\sl E. Witten,}
Two dimensional gravity and intersection theory on moduli spaces,
Surveys in Diff. Geom. {\bf1} (1991), 243--310.

%%%%%%%%%%%%%%%%%%%%%%%%%%%%%%%%%%%%%%%%%%%%%%%%%%%%%%%

\endrule

\endlines{%
  Department of Mathematics, University of Texas at Austin,
  Austin, Texas 78713, U.S.A.\cr
  betz@math.utexas.edu\cr
  Department of Mathematics, Lund University, P.O.~Box 118,
  S-221\thinspace00 Lund, Sweden\cr
  rade@maths.lth.se\cr
  December 22, 1994\cr}

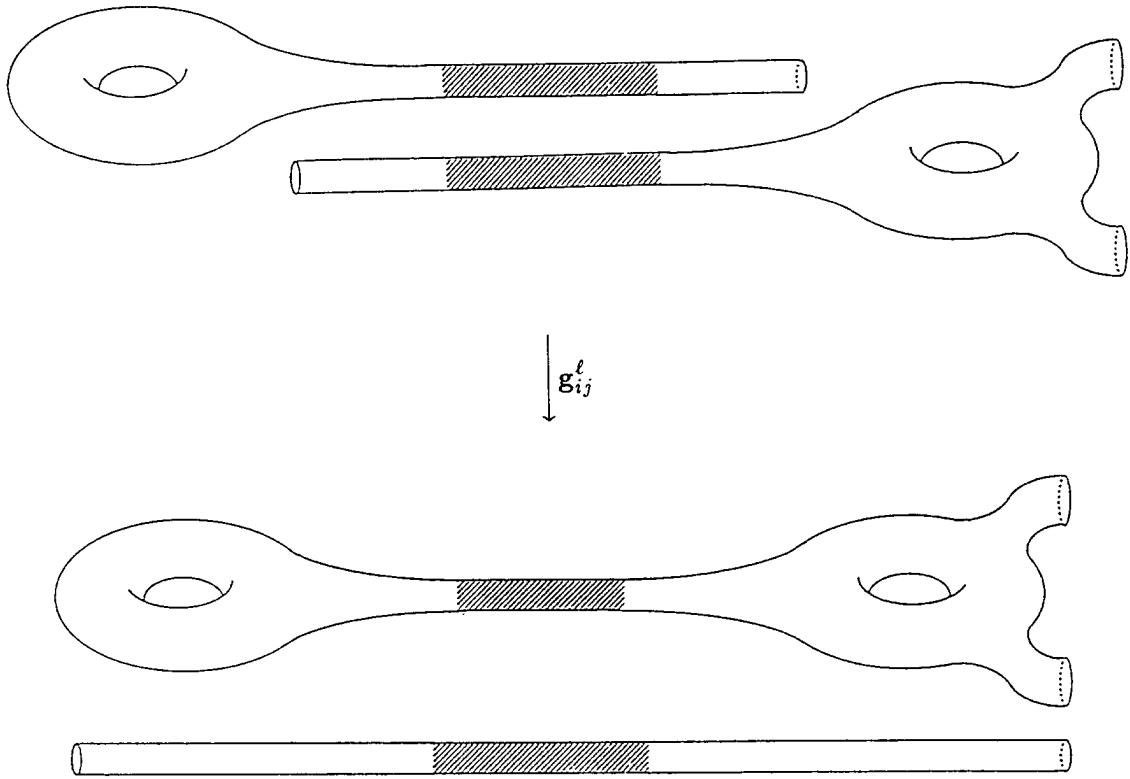
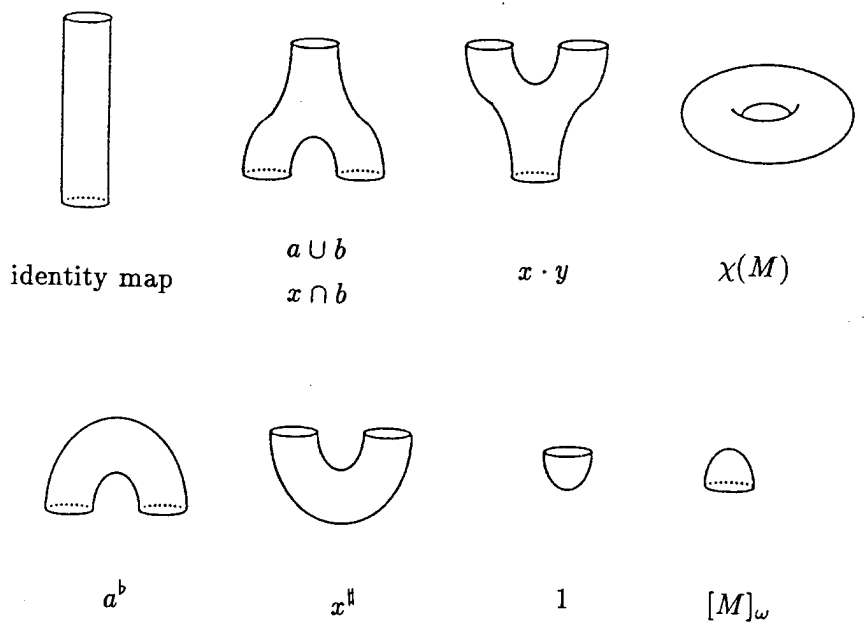
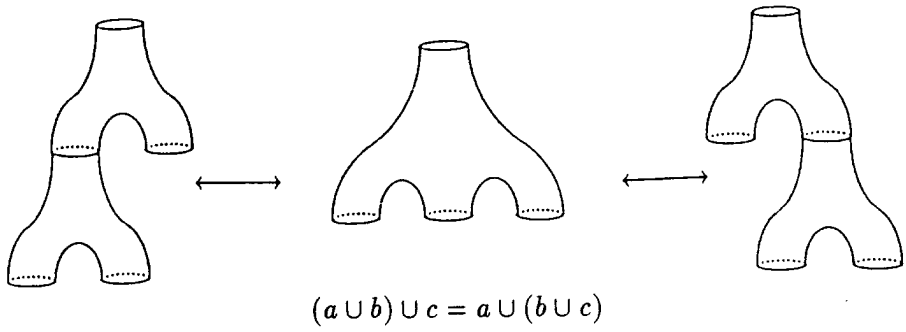
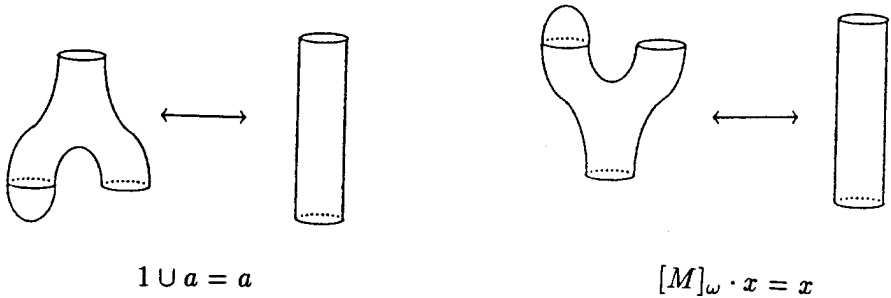
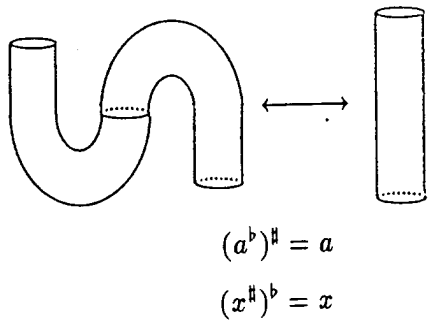
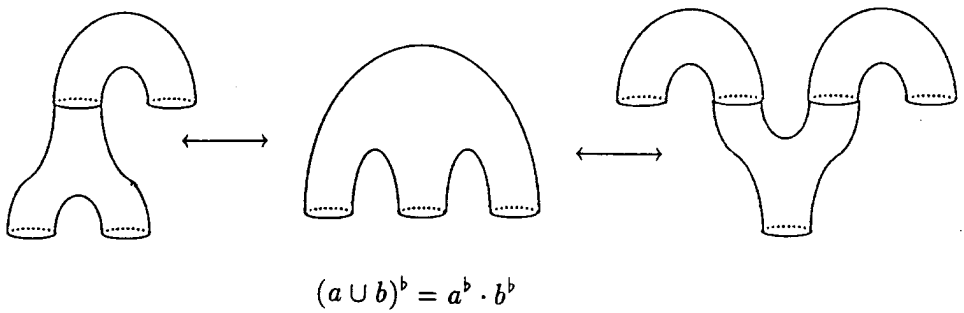
\bye